\documentclass[11pt]{amsart}
\usepackage{amssymb, amsmath, amsthm,ulem}
\usepackage[latin1]{inputenc}
\usepackage{graphicx}
\usepackage[hidelinks]{hyperref}
\usepackage{color}
\usepackage{ulem}
\usepackage{epstopdf}
\usepackage{cancel}
\usepackage{tikz}

\usepackage[a4paper, centering]{geometry}
\usepackage{color}
\usepackage{graphicx}
\usepackage{scalerel,stackengine}

\usepackage{mathtools}

\newtheorem{theorem}{Theorem}
\newtheorem{proposition}[theorem]{Proposition}

\newtheorem{lemma}[theorem]{Lemma}
\newtheorem{corollary}[theorem]{Corollary}

\theoremstyle{definition}

\theoremstyle{remark}
\newtheorem{remark}[theorem]{Remark}

\definecolor{verde}{RGB}{20,150,100}
\definecolor{purple}{RGB}{200,30,200}

\newcommand{\EEE}{\color{black}}

\stackMath
\newcommand\reallywidecheck[1]{%
\savestack{\tmpbox}{\stretchto{%
  \scaleto{%
    \scalerel*[\widthof{\ensuremath{#1}}]{\kern-.6pt\bigwedge\kern-.6pt}%
    {\rule[-\textheight/2]{1ex}{\textheight}}%
  }{\textheight}%
}{0.5ex}}%
\stackon[1pt]{#1}{\scalebox{-1}{\tmpbox}}%
}

\def \div {\mathop {\rm div}\nolimits}
\def \m {\mu _1}
\def \s {\sigma_1} 
\def\R{\mathbb{R}}

\def\N{\mathbb{N}}

\newcommand{\sm}{\setminus}

\newcommand{\Om}{\Omega}
\newcommand{\om}{\omega}

\newcommand{\sq}{\subseteq}

 \newcommand{\ka}{\mathcal K}
 
\newcommand{\cafive}{C_1}
\newcommand{\caseven}{C_2}

\def\ca0{K}

\def \d{\delta}
\def \p{\phi}
\def \e{\varepsilon}
\def \res{\mathop{\hbox{\vrule height 7pt width .5pt depth 0pt
\vrule height .5pt width 6pt depth 0pt}}\nolimits}

\DeclareMathOperator\sign{sign} 
\DeclareMathOperator\Per{Per}

\newcommand{\abs}[1]{{\left|#1\right|}}
\newcommand{\norma}[1]{{\left\Vert#1\right\Vert}}

\begin{document}
\title[]{
On Poincar\'e constants related to \\ isoperimetric problems in convex  bodies
} 

\author[]{Dorin Bucur, Ilaria Fragal\`a}

\begin{abstract}  For any convex set $\Om \subset \R ^N$, we provide a lower bound for the inverse of the Poincar\'e constant in 
$W ^ {1, 1}(\Om)$: it refines an inequality in terms of the diameter due to Acosta-Duran,
via the addition of an extra term giving account for the flatness of the domain. 
In dimension $N = 2$, we are able to make the extra term completely explicit, thus providing a new Bonnesen-type inequality for the Poincar\'e constant
 in terms of diameter and inradius. 
Such estimate is sharp, and it is asymptotically attained when the domain is the intersection of a ball with a strip bounded by parallel straight lines, symmetric about the centre of the ball.
As a key intermediate step, we prove that the ball maximizes the Poincar\'e constant in $W ^ {1, 1} (\Om)$,  among convex bodies $\Om$ of given constant width.    
\end{abstract}

\thanks{}

\date{\today}

\maketitle 
\section{Introduction} 
 
  Aim of this paper is  to     improve   the following optimal Poincar\'e inequality proved by  
Acosta-Duran in \cite{AcDu}:     if  $\Om$  is an open bounded convex set  in $\R^N$ of diameter $D _\Om$, it holds
 \begin{equation}\label{f:AD} \| u \| _{L ^ 1 (\Om)}  \leq \frac{ D_\Om}{2} \|\nabla u\| _{L ^ 1 (\Om)}  \qquad 
\forall u \in 
 W ^ { 1, 1}  (\Om) \text{ with } \int _\Om u = 0 
 \,.
 \end{equation}
The inequality is scale invariant and optimal in the sense that the constant $\frac 12$  in front of the diameter cannot be improved:  this is easily seen by taking 
$\Om$ as an interval in dimension one, or a long thin cylinder in arbitrary dimension, and  a sequence of
functions in $W ^ { 1, 1}  (\Om)$ with zero mean which converge  to a set  obtained bisecting $\Om$ by a hyperplane orthogonal to the axis through its midpoint.

Finding optimal bounds  for Poincar\'e constants in terms of the geometry of the domain is a classical topic, see    \cite{KN15, mazjabook}  and references therein.  Let us just mention that, some years after  the result by Acosta-Duran, nonlinear versions in   $W ^ {1,p}  ( \Om)$    of the classical  inequality proved by 
Payne-Weinberger  in   \cite{PW}  have been given in 
\cite{ENT, FNT12},  with the zero mean constraint appearing in \eqref{f:AD} replaced by  the  different constraint $\int _\Om |u | ^ { p-2} u  = 0$.
  For any $p \in ( 1,  +\infty)$, this 
is related to  the  minimization over $c \in \R$  of  the map $c \to \int_\Om|u-c|^p$, and hence  to the minimization of  
 the first nontrivial Neumann eigenvalue of the $p$-Laplace operator.     In the case $p=1$, this latter constraint reads $ \int_\Om  \sign  u  = 0$.

  Choosing either of the  constraints $ \int _\Om u = 0$ or $ \int_\Om  \sign  u  = 0$, 
Poincar\'e inequalities   in $W ^ { 1, 1} (\Om)$   are  known to be  intimately connected with   
isoperimetric inequalities:  with no attempt of completenes, see  
\cite{milman09} and the related papers \cite{KLS95, SZ99,K03,milman10, RV15}.     Incidentally   let us mention that an estimate in $W^{1,1}(\Om)$ different from \eqref{f:AD}, 
   has been given in  \cite{stein15}, relying on the isoperimetric inequality  by 
 Kannan-Lov\'asz-Simonovits 
\cite{KLS95}.

  Working more in general  in a weighted   Sobolev space $W ^ { 1, 1} _\phi (\Om)$ for some positive function $\p$ in $ L ^ 1 (\Om)$,   
    any of the  two distinct  values  
\begin{equation}\label{f:ms} 
\m (\Om, \p):= \!\!\!\! \inf_{\substack{u \in W ^ { 1, 1} _\p (\Om)\\   \int_\Om( \sign  u ) \phi = 0 }}  \frac {\int _\Om |\nabla u| \p } {\int_ {\Om} |u| \p } \, , \qquad
\s (\Om, \p) := 
\inf_{\substack{u \in W ^ { 1, 1}_\p  (\Om)\\   \int_\Om  u \p  = 0 }}  \frac {\int _\Om |\nabla u| \p} 
{\int_ \Om |u| \p}
 \end{equation} 
  is worth of interest:  each of them   gives   the inverse of a Poincar\'e constant in $W ^ { 1, 1}_\p  (\Om)$, and it  is related  to an optimal  domain  separation problem 
 (see Section \ref{sec:fence}).  In particular such fencing formulations reveal that 
$\m (\Om, \p) \geq \s (\Omega, \phi)$,  so that the smaller quantity
$\s (\Omega, \phi)$ is the natural object of investigation of lower bound inequalities.   
Let us point out that dealing with {\it weighted} constants is not just for the sake of generality,  but it  is 
a natural by-product of the proof  method via dimension reduction conceived by Payne-Weinberger, and exploited also by Acosta-Duran.   
Letting $\s (\Om):= \s (\Om, 1)$,  their inequality \eqref{f:AD} can be  reformulated as    \begin{equation}\label{f:ad} 
\s (\Om) \geq  \frac{2}{D_\Om } \,.
\end{equation} 

  In this paper we improve the inequality \eqref{f:ad} by   adding at its right hand side  some extra term 
which depends on the geometry of the domain, precisely on its flatness. 
Our first  result is  a sharp quantitative version of the inequality \eqref{f:ad} which holds in any space dimension
(Theorem \ref{t:1gap}), the quantitative term giving account for the flatness of the domain. 
  Our second main result  (Theorem \ref{t:newexplicit}) holds in dimension  $N = 2$, and for the constant weight $\phi = 1$, but it has the advantage of bringing an exact, sharp constant, in front of the flatness term, thus providing a Bonnesen-type inequality for $\sigma _ 1 ( \Om)$.  
The two results are proved with completely different methods.

\smallskip 

The former statement reads:

 \begin{theorem}\label{t:1gap}
Let  $N\ge 2$ and $m \in \N \setminus \{ 0 \}$. There exists a  positive  constant   $ C_0 = C_0(N , m) $  such that, 
 for every open bounded convex domain $\Omega$ in $\R ^N$ with diameter $D_\Om$
and John ellipsoid of semi-axes  $a_1 \geq \dots \geq a_N$,  and 
every  $(\frac {1}{m})$-concave function $\phi : \Om \to (0,+\infty)$ 
we have
	\begin{equation}
		\label{quantitative3dp}
		 \s (\Omega, \phi)\ge \frac{2}{ D_\Om}    +        C_0   \frac{a_2^2}{ D _\Om ^{3} }.
		\end{equation} 
		\end{theorem}

  The proof  of Theorem \ref{t:1gap} is based on   the approach introduced in \cite{ABF1, ABF2} and it 
 demands 
some new key  arguments with respect to the following
 nonlinear Poincar\'e inequality proved for any $p \in ( 1, + \infty)$ in \cite{ABF2},  for some positive constant $   \mathcal K _0  = \mathcal K _0 (p, N,m)  $: 
 \begin{equation}\label{f:ABF2}
\mu_p (\Omega, \phi):= \inf_{\substack{u\in W^{1, p}(\Om,\phi)\sm \{0\}\\ \int_\Omega \phi  |u | ^ { p-2} u =0}} \frac{\displaystyle{\int_\Omega \abs{\nabla u}^p} \phi}{\displaystyle{\int_\Omega |u |^p \phi }} 
 \ge \Big  (\frac{\pi_p}{ D_\Om} \Big ) ^ p   +      \mathcal K _0    \frac{a_2^2}{ D _\Om ^{p+2} }\,, \qquad \pi _p:= 2 \pi \frac{(p-1) ^ {\frac 1 p}}{p \sin  \big (\frac \pi p \big ) }  \,. 
\end{equation}  

We postpone to Section \ref{sec:teo1} the detailed presentation of  the proof of Theorem \ref{t:1gap}. Below  we add some more precise comment on 
the relationship between the inequalities  \eqref{quantitative3dp} and  \eqref{f:ABF2}.

\begin{remark}    We warn the reader  that, even the weaker inequality obtained by replacing $\s (\Omega, \phi)$ by $\m (\Omega, \phi)$ in \eqref{quantitative3dp}, 
cannot be obtained 
by passing to the limit in \eqref{f:ABF2} as one might reasonably guess. 
 Indeed,  letting
$\mu _ p (\Om):= \mu _p (\Om, 1)$, in the limit  as $p \to 1 ^+$ we have that  
$ (\mu _ p (\Om) ) ^ {\frac{1}{p} }  \to \m  (\Om):= \mu _1 (\Om, 1)$ (see \cite{Gaj, GajGar}).  
It is  also easily seen that $\pi _ p  \to 2$. 
 But the point is that the explicit value of the constant  $\mathcal K _0$,
as given   for $N = 2$  and $\phi = 1$ in \cite[Theorem 1]{ABF2},  
is infinitesimal (see the Appendix).   The reason is related to the origin of   
the  estimate for $\mathcal K _0$ when $p\in (1, +\infty)$: it  relies on explicit $L^\infty$-bounds for the eigenfunctions, coming  from the elliptic character of the equation, which degenerates as $p \to 1^ +$.  

Marginally let us observe that, in the other limiting case, namely when $p \to + \infty$, a quantitative  inequality for $\mu _ p (\Om)$ 
similar to \eqref{f:ABF2} 
cannot hold. In fact, as $p \to + \infty$,  we have that 
$ (\mu _ p (\Om) ) ^ {\frac{1}{p} }  \to \mu _\infty (\Om)$
(being $\mu _\infty (\Om)$ the first eigenvalue of  a viscosity eigenvalue problem for the $\infty$-Laplacian),  
$\pi _p \to 2$,  and
$\mu _\infty (\Om) = \frac{2}{D_\Om}$  (see \cite[Theorem 1]{EKNT15}).   \end{remark}    

\medskip
  Let us now present our second main result; it holds  in the $2$-dimensional setting,   
which is usually more favourable for the explicit solution of shape optimization problems, 
as it occurs for instance for the famous Blaschke-Lebsegue theorem \cite{B15}, and in many more cases, see e.g. 
the recent work \cite{HL24}.   

Prior to that, we state 
an intermediate inequality of independent interest, which  is an optimal improvement of  \eqref{f:ad}
when the class of competitors is restricted to 
 convex bodies of constant width (see \cite{MMO}):

\begin{theorem}\label{t:constwidth}
 For every open convex set $\Om \subset \R ^2$ having constant width, it holds 
\begin{equation}\label{f:teo2}\sigma _ 1  (\Om) \geq  \frac 8 \pi \frac{1}{D _\Om}\,,
\end{equation} 
with equality sign if and only if $\Om$ is a ball. 
\end{theorem} 

 Relying on Theorem \ref{t:constwidth}, we obtain   the following new  Bonnesen-type   estimate
 in which all the constants are sharp:   

 \begin{theorem}\label{t:newexplicit} 
 For every open bounded convex sets $\Om\subset \R ^2$, it holds
 \begin{equation}\label{f:teo3} \sigma _ 1 (\Om)  >    \frac{2}{D_\Om} + \frac{4}{3} \frac{\rho _\Om ^2  }{D_\Om ^ 3 }\,,
 \end{equation} 
with equality sign attained asymptotically when $\Om$ is the intersection of a ball with a strip bounded by two parallel straight lines symmetric with respect to the center of the ball, at infinitesimal distance. 
  \end{theorem}
 
 \smallskip 
  
 We emphasize that Theorem \ref{t:constwidth} is necessary for the proof of Theorem   \ref{t:newexplicit}.  It might seem quite surprising that 
the inequality \eqref{f:teo2} for bodies of constant width, which is attained at balls,
  does interplay with the  inequality   \eqref{f:teo3}, which is asymptotically attained along domains degenerating into a line segment, thus  quite far from having constant width!    
 Roughly, the origin of this interplay is that our  proof of Theorem \ref{t:newexplicit} is based on the study,  for a fixed $r \in ( 1, \frac 1 2)$,  of the shape optimization problem
  \begin{equation}\label{f:aimshow0}
\min  \Big \{ \frac{\sigma _ 1 (\Om) -2} {r ^2}  \ :\  \text{ $\Om \subset \R ^ 2$  convex with $D_\Om = 1$, $\rho _\Om \geq  r$ } \Big \}        \,. 
\end{equation} 
We stress that optimal shapes for such problem do not obey any a priori expectation: actually,  the fact that optimal sets for quantitative inequalities are not necessarily the same as for the corresponding non-quantitative ones,  has  been already observed in the literature, for instance in the case of the classical isoperimetric inequality (see e.g.\ \cite{CiLe12, AFN11, BCH17}).   In the case under study,  
by virtue of the shortest fence formulation of $\sigma _ 1 (\Om)$, it turns out that solutions to \eqref{f:aimshow0} have a peculiar structure, namely they are obtained 
as the intersection of a convex body of constant width with a strip bounded by two incident or parallel  straight lines. 
In particular, such intersection may happen to be a  full body of constant width, whence the role of Theorem \ref{t:constwidth}. 
We refer to Sections \ref{sec:teo2} and \ref{sec:teo3} for the detailed arguments.

  \smallskip
  Below we provide some comments  about Theorem \ref{t:constwidth}.  
  
 \begin{remark} (i) In  \eqref{f:teo3}, the flatness term is quantified by the inradius and not by the second  John semiaxis, as in the general statement of Theorem \ref{t:1gap}. In fact,  these quantities are equivalent in dimension $2$ and, by the way,    other equivalent quantities could have been used,   such as the width or the depth  orthogonal to a diameter.  On the other hand,  in order to determine explicitly the flatness term,  to our advice 
 the choice of the inradius  is the most  favourable or,  to some extent, the only possible one. 
Actually, when it is replaced by     any of the afore mentioned quantities, it may even occur that  the   corresponding sharp constant is no longer attained along a sequence of degenerating domains,
 leading to a shape optimization problem solved by  
  some non predictable  geometry.  

\smallskip
(ii)  By inspection of the proof of Theorem \ref{t:constwidth}, proceeding  along the same line with  minor changes
 (just replacing the cost functional $ \frac{\sigma _ 1 (\Om) -2} {r ^2} $  in \eqref{f:aimshow0}   by  $\sigma _ 1 (\Om)$), one obtains 
an alternative proof of the inequality  \eqref{f:ad} by Acosta-Duran in dimension $N = 2$
 . 
 
 \smallskip 
(iii)  Still by inspection of the proof, one can easily derive the validity,  for every open bounded convex set $\Om \subset \R ^ 2$ with $D _\Om = 1$, 
 of following  result    stronger than \eqref{f:teo3}  : 
 $$\sigma _ 1 ( \Om)   \ge  \frac{ 8 \rho _\Om  }
 { 2 \rho _\Om \sqrt { 1 - 4 \rho _\Om ^ 2} + \arcsin ( 2 \rho _\Om) 
  } >    2  + \frac{4}{3}  \rho _\Om ^2   \,,$$  
 the first inequality being attained  at   a  disc truncated by parallel lines at distance $\rho_\Om$ from its center.
    We point out that a finer analysis of the second inequality   may lead to  successive  refinements of our Bonnesen-type inequality, with optimal powers and coefficients,  such as, e.g., 
 $$ \sigma _ 1 (\Om)  >    \frac{2}{D_\Om} + \frac{4}{3} \frac{\rho _\Om ^2  }{D_\Om ^ 3 }  +   \frac{76}{45} \frac{\rho _\Om ^4  }{D_\Om ^ 5 }  \,,$$
 and
$$\sigma _ 1 (\Om)  >    \frac{2}{D_\Om} + \frac{4}{3} \frac{\rho _\Om ^2  }{D_\Om ^ 3 }  +   \frac{76}{45} \frac{\rho _\Om ^4  }{D_\Om ^ 5 } + 
 \frac{2648}{945}  \frac{\rho _\Om ^6  }{D_\Om ^ 7 } \,.$$
We leave the  detailed proofs of these inequalities as open questions for the interested reader.  
   \end{remark} 
 
 \smallskip 
The paper is organized as follows.

\smallskip 
	  In  Section \ref{sec:fence} we  provide  some useful characterizations of the two different (inverse) Poincar\'e constants introduced in \eqref{f:ms}, showing in particular that each of  them  is related to a shortest fence problem (see Propositions \ref{p:k} and \ref{p:ad}).  As a consequence, we obtain the inequality $\m (\Om, \p) \geq \s (\Omega, \phi)$  (see Corollary \ref{c:comparison}).  
In absence of a comprehensive reference   in the literature,  the proof of all these results is given in Appendix
(where  for the sake of completeness we also enclose a short proof of the above claim concerning
	the fact that the value given in  \cite{ABF2} for the constant $\mathcal K _0$ in \eqref{f:ABF2} is infinitesimal as $p \to 1^+$).

 \smallskip
 In Sections \ref{sec:teo1}, \ref{sec:teo2}, and \ref{sec:teo3}, 
 we prove  respectively Theorem  \ref{t:1gap}, Theorem \ref{t:constwidth}, and Theorem \ref{t:newexplicit}.

 	\smallskip

	\section{The geometric   side    of optimal Poincar\'e constants in $W ^{ 1, 1}$}\label{sec:fence}

 In order to present the shortest fence formulation of weighted  Poincar\'e constants in $W ^ { 1,1}$, we need to introduce a few notation. 
First of all we observe that 
  the natural space where  the minimization problems in \eqref{f:ms} 
can be  relaxed is a weighted $BV$ space. 
	Let us recall that, for any open bounded set  $\Om$  in $\R^N$ and any positive function $\p$ in $L ^ 1 (\Om)$, 
the  space $BV _\p (\Om)$ can be defined as the space $BV _\mu$ of functions with bounded variation with respect to a positive
measure introduced in \cite{BBF99}, by choosing   $\mu _\p := \p \mathcal L ^N  $ (tacitly meaning that $\p$ is extended to zero out of $\Om$). 
Specifically, we denote by   $L ^ 1  _\phi (\Om)$ the usual weighted Lebesgue space, and 
we set 
  $$
BV _\p (\Om):= \{ u \in L ^1 _\phi  (\Om) \ :\ |D _{\phi } u| (\Om)  < +
\infty \}\ , 
$$ 
 where
 $$
 \begin{aligned}
& |D_\p u| (\Om):= \sup \left \{  -  \langle u, \div (\sigma  \mu_\p) \rangle
 \ : \sigma \in X  _\p  \, , \, |\sigma | 
 \leq 1 \text{ a.e. on } \Om \right \}\,,
  \\  \noalign{\medskip} 
 & X  _\p    :=  \{ \sigma \in L ^ \infty  _{\mu _\p}  (\R^N ; \R ^N) \
:\ \div (\sigma \mu _\p )\in L ^\infty  _ {\mu _\phi} ( \R^N )  \}\,.
 \end{aligned} 
$$ 
Notice that here $\div (\sigma  \mu_\p)$ is meant as the  distribution in  $\mathcal D' (\R^N; \R ^N)$   defined by 
$ \langle \div (\sigma \mu_\p) , \psi \rangle
 := - \int _{\R^N}  \sigma \cdot \nabla \psi  \, d \mu _\phi 
 $  for all  $\psi \in {\mathcal D} (\R^N)$.
When $\phi \equiv 1$ and $\Om$ has Lipschitz boundary, the space $BV _\phi (\Om)$ defined above agrees with the classical space $BV (\Om)$  \cite[Proposition 3.1]{BBF99}.

If the characteristic function $\chi _ E$ of some Lebesgue measurable set $E$ belongs to  $BV _\p (\Om)$, we say that
$E$ has finite $\p$-perimeter,   
  and we  set
 $$\Per _\p ( E, \Om) := |D _\p \chi _ E | (\Om)  \,.$$ 
 Below
we adopt for brevity the notation   
$$|E| _\p:=  \int _ {  E  \cap \Om} \p \qquad \text{ for every measurable set $E \subset \R ^N$}\,.$$

The shortest fence formulation of the constants $\m (\Om, \p)$ and $\s (\Om, \p)$   in \eqref{f:ms} reads:

  \begin{proposition}\label{p:k} 
 Let $\Om$ be an open bounded set in $\R^N$, and  let $\p:\Om \to (0, + \infty)$ be a positive  $( \frac{1}{m})$-concave function, for some $m \in \N \setminus \{ 0 \}$.  Let $\m (\Om, \p)$ be defined as in \eqref{f:ms}, and   let 
 $$
 \overline \m  (\Om, \p) :=\!\!\!\!  \min _ {\substack{u \in BV _\p (\Om; \R^+)\\  | \{ u >0  \} | _\p  \leq  \frac{|\Om| _\p }{2} }} 
 \frac {|D_\p u| (\Om)  } {\int_ \Om |u| \p } \, ,   
  \qquad \m ^* (\Om, \p) := \!\!\!\!  \min _ {\substack{E\subset  \Om \\  |E| _\p  \leq \frac{|\Om| \p }{2}}}  
\frac {\Per_\p(E, \Om) } {|E|_\p}
  \,.$$ 
It holds $$\m  (\Om, \p) = \overline \m  (\Om, \p)  = \m  ^* (\Om, \p) \,;$$ moreover, if  $u$ is optimal for $\overline \m (\Om,\p) $, almost all its level sets $\{u > t \}$ are optimal for $\m  ^*  (\Om,\p) $.  

 \end{proposition}

 \begin{proposition}\label{p:ad} 
 Let $\Om$ be an open bounded set in $\R^N$, and  let $\p:\Om \to (0, + \infty)$ be a positive  $( \frac{1}{m})$-concave function, for some $m \in \N \setminus \{ 0 \}$.  Let $\s (\Om, \p)$ be defined as in \eqref{f:ms}, and   let  

 $$\begin{aligned} & 
\overline \s (\Om, \p)  := \min _ {\substack{u \in BV_\p (\Om)\\   \int_\Om  u \p = 0  }} 
 \frac { |D _\p u| (\Om) } {\int_ \Om |u| \p }\, \qquad \s ^*(\Om, \p)  :=  \min _ {{E\subset \Om}}  
\frac { |\Om |_\p  \Per_\p  (E, \Om) } {2 |E|_\p   |\Om \setminus E|_\p }   \\ 
&\widetilde \s (\Om, \p):= 
\min _ {\substack{A, B \subseteq \Om, \ a, b > 0 \\   A\cap B = \emptyset, \ a|A|_\p  = b |B|_\p } }  
\frac {a \Per_\p  (A, \Om) + b\Per_\p  (B, \Om)   } {a |A|_\p  + b |B|_\p }  \,.\end{aligned} 
 $$ 
 It holds
 $$\s (\Om, \p) = \overline \s  (\Om, \p)  = \s^*  (\Om, \p)  = \widetilde \s (\Om, \p)\,;$$  
 moreover, if  $u$ is optimal for $\overline \s (\Om,\p) $, almost all its level sets $\{u > t \}$ are optimal for $\s  ^*  (\Om,\p) $.  
 \end{proposition} 
 
 \bigskip 
 
 By combining Propositions \ref{p:k} and \ref{p:ad} we obtain 
 
 \begin{corollary}\label{c:comparison}  
We have $ \m (\Om, \p) \geq \s  ( \Om , \p )$, with equality sign in case $\s^ * (\Om,\p)$ is attained at a set      $E$    with $|E|_\p =
 \frac{|\Om|_\p }{2}$. 
 \proof 
We have 
$$ \m    ( \Om , \p ) =  \m  ^*  ( \Om , \p ) = \min _ {\substack{E\subset \Om \\  |E|_\p  \leq \frac{|\Om|_\p}{2}}}  
  \frac {  \Per _\p (E, \Om) }  { |E|_\p} \geq  \min _ {E\subset \Om }  
  \frac { |\Om |_\p  \Per_\p (E, \Om) }  {2 |E|_\p  |\Om \setminus E|_\p}  = \s  ^*  ( \Om , \p )  =   \s  ( \Om , \p )  \,.$$ 
If there exists an optimal set $E$  for $\s  ^*  ( \Om , \p ) $ with $|E| _\p = \frac{|\Om|_\p}{2}$, we have
$\s  ^*  ( \Om , \p ) =  \frac {  \Per_\p (E, \Om) }  { | E |_\p }  \geq \m  ^*  ( \Om , \p ) $, so that $\s ^* ( \Om , \p ) = \m^* (\Om, \p)$ , and
 $E$ is optimal also for $\m^* (\Om, \p)$.    \qed
\end{corollary}

\bigskip
 When $\p = 1$, the shortest fence problems $ \m ^* (\Om):= \m ^* (\Om, 1)$ and $\s^* (\Om):= \s ^* (\Om, 1)$   have been investigated,   
  especially in the planar setting,    in 
 \cite{cianchi89BUMI, cianchi89, EFKNT}, and  references therein. Let us also mention the more recent paper \cite{BogOud}, where  longest minimal length partitions are considered. Some related issues are recalled below. 

 \begin{remark}
 \begin{itemize} 
 \item[(i)] [The case of balls]
  For any ball $B$ in $\R^N$,  $\s ^ *  (B )$ is attained at a half-ball (see \cite[Theorem 1]{cianchi89}), and hence by Corollary \ref{c:comparison} it holds   
 $$\m ^*  (B )  =  \s ^ *  (B ) \,.$$     
\item[(ii)] [The case of triangles]   For an equilateral triangle $T$ in $
\R^2$, we have the strict inequality
 $$    \m ^ *  (T)>\s ^*  (T)   \,.$$ 
 Indeed, by  \cite[Theorem 5]{cianchi89BUMI}, $\m ^* (T)$ is attained at a circular sector $S_r$ with centre at a   corner of $T$ and radius $r'$ such that
 $|S_r| = \frac{|T|}{2}$. If $T$ has unit sides, it is readily checked that  $r'= \sqrt{ \frac {3 \sqrt 3} {4\pi}} $ and hence
 $\m ^* (T) = 4 \sqrt { \frac{\pi}{3 \sqrt 3} }$. On the other hand, we have
$\s ^* ( T) \leq\min _ { r\in [0,1]} \frac { |\Om | \Per (S_r, \Om) }  {2 |S_r| |\Om \setminus S_r|}$, and 
such minimum is attained at $r'' = \sqrt  {\frac {3} {2 \sqrt 3 \pi}}$, yielding the inequality  $\s^* (T) \leq 3 ^ {\frac{3}{4}} \sqrt {\frac{\pi}{2}} <  4 \sqrt { \frac{\pi}{3 \sqrt 3} }
$.   
 
 \item[(iii)] [Longest shortest fence] 
For any open bounded convex  set $\Om\subset \R ^2$, if $B\subset \R ^2$   is the ball with the same area as $\Om$, it holds 
 $$\m  ^ * (\Om) \leq \m  ^ * (B) \, , \qquad \s ^ * (\Om) \leq \s ^ * (B)  \,.$$ 
 Indeed, the first inequality has been proved in \cite[Corollary 2]{EFKNT}. The second inequality is satisfied because, by Corollary \ref{c:comparison} and statement (i) above, we have
 $\s ^ * (\Om) \leq \m ^ *(\Om) \leq \m ^ * (B) = \s ^ * (B)$. 

 \smallskip 
 \item[(iv)] [Planar  isoperimetric sets] 
For any open bounded convex  set $\Om\subset \R ^2$, $\m ^* (\Om )$ is attained at a convex set $E \subset \Om$ with $|E|  = \frac{|\Om|}{2}$ (see \cite[Theorem 3]{cianchi89BUMI}).   
 
  \end{itemize}  
  \end{remark}

 \medskip

\section{Proof of Theorem \ref{t:1gap}}\label{sec:teo1}

Before starting with the proof, we wish to give a short 
highlight of its
main ingredients, 
to explain the way they imply the result, and to outline how  they are presented in the next   sub-sections.    

The proof  relies on two key arguments.   The first 
is a   quantitative one dimensional inequality which 
is specifically related to the functional setting, here $W^{1,1}_\phi(\Om)$,
and holds for weights $\phi$ of particular feature.
The second one is  an area estimate for optimal sets in an isoperimetric problem inside $\Om$,  which enables us to control  the area proportion of each cell, 
in a convex partition that we perform following the approach introduced in \cite{ABF1, ABF2}. 
In order to describe more in detail these two  arguments, we need to  sketch the global strategy of the proof.
We refer to   \cite{ABF1, ABF2} for 
all the details that here are omitted so  to avoid repetitions.

\smallskip 
We first observe that it is enough to prove the result in dimension $N = 2$. The passage to higher dimensions can be done with minor modifications by following the  same line as in \cite[Section 6]{ABF1}. 
Moreover, it is not restrictive to  assume   that $\Om$ and $\p$ are smooth, and that $D _\Om = 1$.  
To fix the ideas, we can also work in  a system of coordinates such that a diameter is the horizontal segment $ [- \frac 1 2, \frac 1 2 ] \times \{0\} $.  
Since  $N=2$,  $a_2$ is comparable to the width  of $\Om$, or equivalently to  the maximum $\eta$ of the lengths of vertical sections of $\Om$, that we can assume small enough. If the inequality is proved for small $\eta$, then  its validity for arbitrary $\eta$ follows by the same argument as in \cite[Lemma 24]{ABF2}. 

Given a function $u \in W ^ { 1, 1} _\p (\Om)$, we consider a modified notion of partition \`a la Payne-Weinberger,  that we call  {\it $\phi$-weighted $L ^ 1$-equipartition of $u$}, 
specifically a family $ \{ \Om _ 1, \dots, \Om _n \}$ of mutually disjoint convex sets, with  $\Om = \Om _1 \cup \dots \cup  \Om _n $, such that 
$$\int _{\Om _i }  u  \phi= 0 \qquad \text{ and }   \qquad \int _{\Om _i }   | u  |   
\phi   = \frac{1}{n} \int _{\Om  } |  u |     \phi  \qquad \forall i = 1, \dots, n\,. $$ 
 They key feature of such a partition is that, if for a fixed  proportion of cells (say $ \theta  n$, 
 with $\theta \in (0, 1)$), it holds $ \frac { \int _{\Om_i} 
 |\nabla u| \p } {\int_ {\Om  _i} |u| \p }  \geq 2 + \delta$ for some $\delta >0$,  
then 
 \begin{equation}\label{f:pezzetto}
\frac { \int _\Om  |\nabla u| \p } {\int_ \Om |u| \p }    \geq  \frac{1}{n}  \sum _{i = 1} ^n  \frac { \int _{\Om_i} 
 |\nabla u| \p } {\int_ {\Om  _i} |u| \p } \geq \frac{1}{n} \big [ (n- \theta n ) 2 + \theta n  ( 2 + \delta)\big ]  =   2 + \theta \delta  \,.
 \end{equation}

Then the proof method consists in getting under control  the collective behaviour of the cells.  Such control is reached via the study of a list of binary cases,  
each  one leading either to the desired result, or to a distinction into two more sub-cases, 
up to the last sub-case when we arrive at a contradiction.
The list of binary cases works as follows:

\begin{itemize} 
\item{} {\it Case 1:} For a fixed   proportion of cells, the diameter $D_i$ of $\Om _i$  is  ``small'' (i.e., $1- D_i$ is sufficiently large in terms of $\eta$).  In this case the result follows directly from \eqref{f:pezzetto}  by applying Acosta-Duran inequality.

\smallskip \item{} For a  fixed   proportion of cells, $D_i$ is ``large''. In this case the cells $\Om _i$ can be ordered vertically in a strip upon the fixed diameter 
and the intersection of $\partial \Om _i$ with such a strip
is a polygonal line. We then distinguish two sub-cases: 
\smallskip

\item{} {\it Case 2.1:}   For a   fixed  proportion of cells, $\partial \Om _i$ has some vertex in the strip. In this case we are back to a situation of the same type as in Case 1. 

\smallskip 

\item{} {\it Case 2.2:}  For a  fixed  proportion of cells,  $\partial \Om _i$ has no vertex in the strip.  Such cells   $ \Om _i$ have an affine profile $h _i$  inside the strip, and we distinguish two further sub-cases:

\smallskip 
\item{}  {\it Case 2.2.1:}  For a   fixed  proportion of cells,   the energy of $u$ inside $\Om _i$ is far above $2$.
Then  the result follows from directly from \eqref{f:pezzetto} simply by applying Acosta-Duran inequality. 

\smallskip 

\item{}   {\it Case 2.2.2:}  If, for a   fixed  proportion of cells,  the energy of $u$ inside $\Om _i$ is controlled from above, then for those cells, 
  since they tend to collapse to one-dimensional segments,  we obtain   $$ \s (I _ {D_i}, h_i \phi )  \leq  \frac{5}{2} \qquad \text{ and } \qquad  \frac{ \int _{\Om _i}  |\nabla u|    \phi }{\int _ {\Om _i}  |  u|    \phi    }    
 \geq \s  ( I _ {D_i}, h_i  \phi ) + o ( 1)\,, 
$$ 
where $I _ {D_i} =[- \frac { D_i }{ 2},  \frac { D_i }{ 2  }]$, and $o ( 1)$ is an infinitesimal as the number of cells tends to $+ \infty$. 

\noindent Then, to finish the proof, we need to set up  a first key ingredient, which is an ad-hoc {\it quantitative one-dimensional inequality}, of the kind 
\begin{equation}\label{f:quanti1} \s (I _ {D_i}, h_i \phi )  \leq  \frac{5}{2} \quad  \Rightarrow \quad  \s (I _ {D_i}, h_i \phi) \geq 2  +   \cafive   \min _{[ - \frac{D_i}{2} + \delta_0, \frac{D_i}{2}  - \delta_0]} \Big ( \frac{h_i'}{h_i} \Big ) ^ 2   \,.\end{equation}  
\noindent If \eqref{f:quanti1} is in force, we can proceed and conclude  by distinguishing two final sub-cases: 

\smallskip 
\item{}  {\it Case 2.2.2.1}:  If the extra-term in \eqref{f:quanti1} is large in terms  of $\eta$, 
the result follows from \eqref{f:pezzetto} and \eqref{f:quanti1}.

\smallskip 
\item{}  {\it Case 2.2.2.2}:  If the extra-term in \eqref{f:quanti1} is small in term of $\eta$, 
we reach a contradiction, 
coming  essentially from the fact that the pile of cells   would be too high:    heuristically,   the shapes of the cells would be close to rectangles with controlled sides, 
so that the diameter of the pile can  be estimated, leading to a value strictly above $1$.  

\noindent This geometric argument works provided
one has,   as a  second key ingredient,  a   {\it lower bound for the area proportion of each cell},  in the sense that 
\begin{equation}\label{f:esticelle}
|\Om _i|\geq \caseven \frac{ |\Om|}{n} \,.
\end{equation}    
 
\end{itemize}

 In conclusion, the two main arguments  needed to obtain Theorem \ref{t:1gap} are 
the inequalities  \eqref{f:quanti1} and  \eqref{f:esticelle}.  
Both of them   demand  a completely different approach from their counterparts used in  \cite{ABF2} 
to deal with nonlinear Poincar\'e inequalities in $W^{1,p}$, for $p\in (1, +\infty)$.

\smallskip 
 The remaining of this section is organized as follows.

	\smallskip
	In Section \ref{sec:area},  thanks to the reformulation of $\s (\Omega, \phi)$ as an optimal profile problem, and to 
the fact that we are able to gain a control on the area  of  an optimal isoperimetric set, 
 we obtain the lower bound \eqref{f:esticelle}  for the area of each cell  in the convex partition (see Proposition \ref{p:consotto}  and Remark \ref{r:consottoarea}). 
 
\smallskip
In Section \ref{sec:1d}, we obtain 	
	the quantitative one-dimensional inequality \eqref{f:quanti1} (see Proposition \ref{p:1d}). Loosely speaking, the reason 
why we cannot follow the proof line used for $p \in (1, + \infty)$ in \cite[Section 4.1]{ABF2}, 
	is that, for $p =1$, we do not have at our disposal a pointwise form of the Euler-Lagrange equation for solutions to the  variational problems in \eqref{f:ms} (namely, for their relaxation in  spaces of $BV$ type). 
Thus we adopt an ad-hoc approach, whose starting point is  the elementary inequality pointed out by Acosta-Duran 
\begin{equation}\label{f:nuce}
\int _0 ^ 1 |u (x ) | \rho ( x) \, dx \leq 2  \int _0 ^ 1  |u' ( x)|  \p  ( x) J _\rho  ( x)  \, dx \,, \qquad \forall u \in W ^ { 1, 1} _\rho (0, 1) \hbox{  with } \int _0 ^ 1 u (x) \rho ( x) = 0\,, 
 \end{equation}
 where 
\begin{equation}\label{d:Jrho}
J _{\rho} (x)   :=  \frac{ \Big (\int _0 ^ x \rho ( t) \, dt  \Big ) \Big (\int _x ^ 1 \rho ( t ) \, dt  \Big )    }{\Big (  \int_0 ^ 1 \rho ( t) \, dt \Big ) \rho ( x)   } 
\end{equation}
In view of  \eqref{f:nuce}, in order to obtain a one-dimensional Poincar\'e inequality on $W ^ { 1, 1} _\rho (0, 1)$, 
it is enough to prove a pointwise upper bound on the function $J _{\rho}(x)$ on the interval $(0, 1)$. 
In fact, the inequality  \eqref{f:AD} was based on the upper bound
  $J _{\rho} (x) \leq \frac{1}{4}$, which was stated for power concave weights in   \cite[Lemma 3.1]{AcDu}. 
The heart of the matter in order to refine Acosta-Duran result consists  in showing that, for weights $\rho$ of the same form as in \eqref{f:quanti1},  the upper bound $\frac{1}{4}$ can be successfully improved. This purely one-dimensional target is achieved in Lemma \ref{l:1refinement}.

\bigskip
\subsection{Lower bound for the area proportion of cells}\label{sec:area} 

In this section, relying on the optimal profile formulation $\s^* (\Om, \p)$   of the eigenvalue $\s (\Om, \p)$, we achieve the lower bound \eqref{f:esticelle}
 on the area proportion of the cells.

We start by noticing that, thanks to Proposition \ref{p:ad}, in the proof strategy described at the beginning of this section, we can replace the Rayleigh quotient 
$\frac{\int _\Om |\nabla u| \p} {\int _\Om |u| \p}$ of a generic function in $W ^ { 1, 1} _\p (\Om)$ with $\int _\Om ({\rm sign} u ) \p = 0$, by the Rayleigh quotient  
$$\frac{| D\overline u | _ \p (\Om)} {\int _\Om |\overline u| \phi}\,$$   
where  
  \begin{equation}\label{f:ue}  \overline u:= \chi _{\Om \setminus E} - \frac{|\Om \setminus E| _\p}{| E |_\p } \chi _ { E} \,,
  \end{equation} 
being $E \subseteq \Om$ an  optimal set for $\sigma _ 1 ^*  (\Om, \p)$. 

In turn, by Proposition \ref{p:BV}, there exists  a sequence $u _h  \subset \mathcal D ( \R ^N)$  such that $u _ h \to  \overline u  $ in $L ^ 1  _\p (\Om)$  
and $\lim _h \int _\Om |\nabla u _h |  \phi =  |D _\p  \overline u| (\Om) $; moreover, 
since $\int _\Om \overline u \p = 0$, up to a translation we can assume that $\int _\Om u_h \p = 0$.   
Thus, we are reduced to control the area of the cells of a 
 $\phi$-weighted $L ^ 1$-equipartitions of an element of this sequence, that we   shall    call for brevity a {\it smooth admissible approximation of $\overline u$}. 
 As detailed in Remark \ref{r:consottoarea}, this goal  is reached as a consequence of the next result.  

\begin{proposition}\label{p:consotto}
 Let $\Om$ be an open bounded convex set in $\R^N$ with diameter $D _\Om = 1$,  and let $\p$ be a positive function on $\Om$, which is $( \frac{1}{m})$-concave for some $m \in \N \setminus \{ 0 \}$.  Then there exists a   constant $\gamma = \gamma (N, m ) \in (0, \frac{1}{2})$ such that, if 
 $\s  ^* (\Om, \p)    \leq 3$ and $E$ is an   optimal set for $\s  ^* (\Om, \p)$ with $|E |_\p \leq \frac{1}{2} |\Om |_\p$, it holds   
\begin{equation}\label{f:areae} |E | _\p \geq  \gamma  |\Om | _\p\,.
\end{equation} 
\end{proposition}

\proof We first prove \eqref{f:areae}  when $\p \equiv 1$. Letting 
$$ u:= \chi _{\Om \setminus E} - \frac{|\Om \setminus E| }{| E | } \chi _ { E} \,,
$$
the assumption $\s  ^* (\Om)    \leq 3$ gives 
\begin{equation}\label{f:star}  3   \geq \frac { |D  u| (\Om) } {\int_ \Om |u|  } 
= \frac{|\Om |    \Per  (E, \Om)  }{ 2 |E|  |\Om \setminus E |   }  \geq \frac{  \Per  (E, \Om)  }{  2  | E |  }   
\end{equation}
	  Let us denote by $a_1 \geq a_2 \geq \dots \geq a_N $ the John semi-axes of  $\Om$, and let us assume without loss of generality that 
	$a_1= 1$.  We consider the transformation 
	$T(x)= X$ defined by 
	$$
	X_1=x_1, \quad X_2= \frac{x_2}{a_2}  , \quad \dots\quad,   \quad X_N= \frac{x_N}{a_N} \,,
	$$
	and we set $A: = T (\Om)$. 
	We define on $A$  the function $v : = \chi _ E \circ T ^ { -1} = \chi _{T ( E)}$  i.e., 
		$$ \
	v (X_1,\dots , X_N) := \chi _ E (X_1, a_2  X_2, \dots , a_N X_N )\,.$$	
	
	We observe that
\begin{equation}\label{f:changevar} \frac{  |D \chi _E | (\Om)   }{  \int _\Om \chi _ E     }  \geq \frac{  |D \chi _{T (E)} | (A)   }{ \int _A \chi _ {T(E)}   } =  \frac{ \int _A |Dv|   }{\int _A |v|  }  \,.
\end{equation} 

Indeed, by the change of variables $x = T ^ {-1} ( X)$,  for any $\varphi \in \mathcal D ( \R ^N)$,  letting $\psi = \varphi \circ T ^ { -1}$, it holds 
$$\frac{\int _\Om |\nabla \varphi| } {\int _\Om | \varphi|  } =  \frac  {    \int _A\Big (  \sum _{i = 1} ^N \frac{1}{a_i ^2} |\nabla_i \psi|  ^2 \Big ) ^ {\frac{1}{2}}   }{ \int _A |\psi|}   \geq \frac {\int _A  |\nabla \psi|  }{ \int _A |\psi|}   \,.  $$ 
Hence, 
$$\begin{aligned} 
\frac{  |D \chi _E | (\Om)   }{  \int _\Om \chi _ E     } & = \inf \Big \{ \liminf \frac{\int _\Om |\nabla \varphi _h | } {\int _\Om | \varphi _h|  }  \ :\  \{ \varphi _h\} \subset \mathcal D ( \R ^N) \,, \ \varphi _h \to \chi _E  \text{ in } L ^ 1 ( \R ^N) \Big \}
\\ 
& \geq \inf \Big \{ \liminf \frac{\int _A |\nabla \psi _h | } {\int _A | \psi _h|  }  \ :\  \{ \psi _h\} \subset \mathcal D ( \R ^N) \,, \ \psi _h \to \chi _A  \text{ in } L ^ 1 ( \R ^N) \Big \} =  \frac{ \int _A |Dv|   }{\int _A |v|  } \,.
\end{aligned} $$

We conclude that 
$$  \begin{aligned} 
6 \geq \frac{  \Per  (E, \Om)  }{  | E |  }  &     \geq   \frac{ \int _A |Dv|   }{\int _A |v|  }  \geq C_N
\frac{  \Big ( {\int _A |v|  \Big ) ^ {\frac{N-1}{N}}    }}{\int _A |v|  }  =  C_N
  \Big ( {\int _A |v|  \Big ) ^ { - \frac{1}{N}}    } 
\\ \noalign{\medskip} 
& \geq  C'_N
  \Big ( \frac{1}{|A|} {\int _A |v|  \Big ) ^ { - \frac{1}{N}}    } =  C'_N  \Big ( \frac{   | A |  }{  |T(E) |  }   \Big ) ^ {  \frac{1}{N}}   
  =  C'_N  \Big ( \frac{  |\Om  |   }{   | E |  }   \Big ) ^ { \frac{1}{N}}    \,,
  \end{aligned}    
 $$ 
 where we have used  \eqref{f:star} for the first inequality and  \eqref{f:changevar} for the second one,  while the third 
 and fourth inequalities hold   for some dimensional constants $C_N$ and $C'_N$,   respectively because the 
immersion constant of the continuous embedding $ BV ( A) \subset L ^  {\frac{N-1}{N}} (A)$ and the volume of $A$   are purely dimensional. This is a consequence of the fact that,  by construction, the outradius and inradius  of $A$ are bounded respectively from above and from below. 
The proof is thus achieved for $\phi \equiv 1$ (and in arbitrary space dimension).  

The validity of the inequality \eqref{f:areae} for a general  weight $\phi$  follows by applying the statement with $\phi = 1$ to
the ``inflated'' convex set $\widetilde \Om \subset \R ^ { N + m}$  defined by 
	 $$
	 \widetilde{\Omega } = \left\{ (x,y) \in \R^{N}\times \R^{m} : \, x \in \Omega ,\,  \norma{y}_{\R^m} <   \omega _m ^ { - \frac{1}{m}}  \p^{1/m}(x)\right\}\,.
	 $$ 
We omit the detailed argument, since it can be easily reproduced as a variant in $BV$ of   Proposition 8 in \cite{ABF2}.  
 \qed 

\begin{remark}\label{r:consottoarea} 
 From Proposition \ref{p:consotto}, we readily obtain the validity of the area estimate  \eqref{f:esticelle} for 
  the cells of a $\phi$-weighted $L ^1$-equipartition of $u _h$, being $u_h$  
 a smooth admissible approsimation  of the function $\overline u$ defined in \eqref{f:ue}
  (with the terminology introduced just above the statement of Proposition \ref{p:consotto}). 
 
Indeed, notice first that,  as a consequence of the inequality \eqref{f:areae},  the function $\overline u$ defined in \eqref{f:ue}    satisfies  the $L ^ \infty$ estimate 
$$1  \leq \overline u ( x)  \leq   \frac{|\Om \setminus E| _\p}{| E |_\p } \leq \frac{1 - \gamma}{\gamma}  \qquad \forall x \in \Om \,.
$$
In turn, this implies that the smooth admissible approximation $u _h$ satisfys  for $h$ sufficiently large
 $$\frac { 1}{2}   \leq u  _h ( x)  \leq    \frac{2(1 - \gamma)}{\gamma}  \qquad \forall x \in \Om \,.
$$

Now, from the upper bound on $u_h$, we have 
$$
		 \frac{\int _{\Om  } |  u_h | \p  }{n\int _{\Omega _i }\phi }= \frac{\int _{\Om _i }   | u _h |   
\phi   }{\int _{\Omega _i} \phi } \leq     \frac{2(1 - \gamma)}{\gamma}  \,.$$ 
	Hence, using also the lower bound on $u _h$, we have 
	$$|\Om _ i | \| \p \| _{L^\infty (\Om)} \geq \int _{\Omega _i} \phi  \geq  \frac{   \gamma  }{2n ( 1 - \gamma) }\int _{\Om  } |  u_h | \p   \geq 
	  \frac{ \gamma  }{4n  ( 1 - \gamma)  } { |\Om | _\p    }  \geq  
	  \frac{ \gamma   }{n ( 1 - \gamma)  }   \frac{1}{2 ^ { N+m+2} } |\Om|  \| \p \| _{L^\infty (\Om)} \,.
	 $$ 
Here the last equality is satisfied  assuming with no loss of generality that $0 \in \overline \Om$ is a maximum point for $\phi$, so that by $\big ( \frac 1m)$-concavity  it holds 
$\phi (x) \ge  \frac {1}{2^m} \|\phi\|_{L^\infty(\Om)}$ for every  $x \in \frac 12 \Om $.  

We conclude that the estimate \eqref{f:esticelle} is fulfilled by taking 
$\caseven = \caseven (N, m):=  \frac{\gamma   }{n ( 1 - \gamma)  }   \frac{1}{2 ^ { N+m+2} } $.  
\end{remark}
\EEE

\subsection{One-dimensional refined Poincar\'e inequality}\label{sec:1d}  

  In this section we prove the quantitative one-dimensional inequality \eqref{f:quanti1}, precisely in the following form:  

\begin{proposition}\label{p:1d} 
There exist constants 
 $\cafive = \cafive ( m)>0$   and $ \delta_0 = \delta_0   ( m)  \in (0, \frac{1}{2})$  such that, 
if $\rho = \phi h$ is a positive weight  on $I = (0, 1)$,  with  $\| \rho \| _\infty = 1$,   $\phi$ positive 
$(\frac{1}{m})$-concave for some $m \in \N \setminus \{ 0 \}$,  and $h$ affine on $[\delta_0, 1-\delta_0]$, it holds 
\begin{equation}\label{f:1dq}
\s ( I , \rho) \leq    \frac 5 2   \ \Rightarrow \
\s ( I , \rho)  \geq   2   +  \cafive  \min _{[ \delta_0, 1 - \delta_0  ]} \Big ( \frac{h'}{h} \Big ) ^ 2 \,.
\end{equation}
\end{proposition} 
 
\bigskip
\begin{remark}   The above inequality scales as follows. If $I_d=(0,d)$ and $\rho$ satisfies the assumptions of 
Proposition \ref{p:1d} on $I _d$, by applying \eqref{f:1dq}  
 to  the weight   $R(x)=\rho (d\cdot x)$, 
one gets 
\begin{equation}\label{f:1dqbis} 
\s (I_d, \rho ) \leq  \frac{ 3}{d}    \ \Rightarrow \ \s ( I_d , \rho)  \geq   
 \frac{ 2}{d}   +  \cafive d  \min _{[d  \delta_0   , d - d \delta_0  ]} \Big ( \frac{h'}{h} \Big ) ^ 2 \,.
\end{equation} 
  In particular, for $d \in [ \frac{2}{3} , 1]$,     from  
\eqref{f:1dqbis} one gets 
\begin{equation}\label{f:1dqbis.1} 
  \s (I_d, \rho ) \leq  3  \ \Rightarrow \     \s ( I_d , \rho )  \geq     2   +  \cafive\frac 23 \min _{[ \delta_0  , d -  \delta_0  ]} \Big ( \frac{h'}{h} \Big ) ^ 2. \end{equation} 
\end{remark}

\bigskip

The proof of Proposition \ref{p:1d} is based on the next two lemmas. 
In the first one below we give refinements of the inequality \eqref{d:Jrho} 
by Acosta-Duran.   In particular, we shall exploit the inequality stated in item (ii). 
Although we are not going to use it, 
for convenience of the reader and for possible future refinements, we also give in item (iii) an inequality which is of different nature from (ii)
 (not being a consequence of it, neither the converse).

\begin{lemma}\label{l:1refinement} 
Let $\rho$ be a positive weight on $[0, 1]$,  and let 
$J _\rho (x)$ be defined as in  \eqref{d:Jrho}.   
\begin{itemize}
\item[(i)] If $\rho$ is log-concave, it holds 
$$ J _ {\rho}  (x) \leq x ( 1-x)  \qquad \forall x \in (0, 1)$$
with equality sign   if and only if either $\rho$ is constant, or $x = 0$, or $x = 1$.

\smallskip 
\item[(ii)]  For every $\delta \in (0, \frac 14)$ and every  $m \in \N \setminus \{ 0 \}$, 
there exists
$  \Lambda_0 = \Lambda  _0 ( \d, m) >0$ such that, 
if $\rho = \phi h$ is a positive weight   with  $\| \rho \| _\infty = 1$, being   $\phi$ positive 
$(\frac{1}{m})$-concave and $h$ affine on $[\delta, 1 - \delta]$,  
 it holds 
 $$
\begin{aligned} J _ {\rho}  (x) & 
\leq x ( 1-x) - \Lambda  _0 \min _{[\d, 1 - \d]} \Big ( \frac{h'}{h} \Big ) ^ 2  \qquad \forall x \in ( 2\delta, 1 - 2\delta) 
\,.
\end{aligned} 
 $$
 
 \smallskip
\item[(iii)]   
For every $\delta \in (0, \frac 1 2 )$ and every  $m \in \N \setminus \{ 0 \}$, 
there exists  a positive constant  $\Lambda_ 0' = \Lambda _0 ' (\d, m)$ such that, 
if $\rho$ is  a  positive weight  with $\| \rho \| _\infty = 1$ which is $(\frac{1}{m})$-concave,  it holds 
 $$ J _ {\rho}  (x) \leq \frac{1}{4} - \Lambda_0 '  \Big ( \frac{\rho' ( x)} {\rho (x) }  \Big ) ^2  \qquad \forall x \in ( \delta, 1 - \delta) \,.$$ 
\end{itemize}

 \end{lemma} 

\proof  (i)  Let $x \in (0, 1)$ be fixed.  We  remark that, since
$$ J _ {\rho}  (x) = \frac{1}{\rho (x)}  \cdot \frac{1}{\frac{1}{\int_0 ^ x \rho} + \frac{1}{\int _x ^ 1 \rho}}  \, ,$$ 
we have that  $  J _ {\rho}  (x) \leq J _ {\tilde \rho}  (x)$
for any positive log-concave function $\tilde \rho$  such that 
$\tilde \rho(t)  \geq \rho(t)$ for every  $t \in (0, 1)$, with equality at $t = x$.

Next we observe that, letting  $g := \log \rho$, by assumption it holds $g ( t)  \leq g ( x)+ g' ( x) ( t-x)$  $\forall\, t \in (0, 1)$. Hence  
setting
\begin{equation}\label{f:ab} a := g' (x) = \frac{\rho' (x)} {\rho (x)}\, , \qquad   b :=  e ^ {  g(x) -  x g' (x) }  =   \rho ( x) e ^ { 
- a x   }  \, , 
\end{equation}
the weight $\rho$ satisfies the estimate \begin{equation} \label{f:113a}
\rho ( t) \leq b  e ^ { a t } \qquad \forall   t \in ( 0, 1)  \,. 
\end{equation}
 Applying the above remark  we see that 
\begin{equation}\label{f:parziale}
J _\rho (x) \leq  J _ {b e ^ { at}}  (x)  \,, 
\end{equation} 
with equality if and only if $\rho ( t) \equiv b e ^{ at}$. 
Hence it is enough to prove the inequality 
$$ J _ {b e ^ { at}}  (x)   \leq x ( 1-x) \,. $$ 
For $a = 0$, namely when $\rho$ is constant,  the above inequality becomes an equality for every $x \in [0, 1]$.  

Thus we can assume that $ a\neq 0$,  in which case we have 
$$ \frac{ \Big (\int _0 ^ x e ^ { at}  \, dt  \Big ) \Big (\int _x ^ 1 e ^ { at}   \, dt  \Big )    }{\Big (  \int_0 ^ 1 e ^ { at}   \, dt \Big ) e ^ { ax}  }  = 
\frac{( e ^ { ax} -1) ( e ^ a - e ^ { ax})}{a (e ^ a -1) e ^ { ax} }$$ 
So we have to show that, 
for $a \neq 0$, it holds
$$E _ a ( x):= \frac{( e ^ { ax} -1) ( e ^ a - e ^ { ax})}{a (e ^ a -1) e ^ { ax} }  - x ( 1-x) \leq 0 \,, $$
with equality if and only if $x = 0 $ or $x = 1$. 
We have $E_a ( 0 ) = E _ a ( 1) = 0$, and 
$$
E _ a'  ( x)= \frac{ e ^ { a ( 1-x)}- e ^ { ax}    }{ e ^ a -1 }   - 1 + 2x \,.
$$  
We notice  in particular that 
$E _ a'  ( x) =-  E _ a'  (1- x)$,
so that we are reduced to prove that
$$ E _ a'  ( x) \geq 0 \quad \forall x \in \big ( \frac{1}{2}, 1 \big ) \,.$$  
Indeed, this implies that the maximum of $E_a$ on the interval $(0, 1)$ equals  $E_a ( 0 ) = E _ a ( 1) = 0$ (and it is  attained only at $x = 0 $ and at $x = 1$).  

Since $E_{a} = E _{-a}$ we can assume with no loss of generality that $a>0$. Then, in terms of the variable $y := x - \frac{1}{2}$, the above inequality is equivalent to
$$H _ a (y):= e ^ { - a y} - e ^ { ay} + 2y \big ( e ^ { \frac{a}{2} } - e ^ { - \frac{a}{2}}  \big )  \geq 0 \qquad \forall y \in \big  ( 0 , \frac{1}{2} \big  ) \,.$$ 
Since $a>0$, and 
$$
H'' _ a (y ) = a^ 2  (e ^ { -a y} -  e ^ { ay}  )  \,,
$$  
we have that $H_a$ is concave, and hence it attains minimum at $0$ or at $\frac{1}{2}$, where it is equal to zero. 
Accordingly, we have $E' _a \geq 0$ on $\big ( \frac{1}{2}, 1 \big )$ as required. 
 
\medskip 
(ii) Let $x \in (2 \d, 1 - 2 \d)$ be fixed, and let $a,b$ be the corresponding quantities defined by \eqref{f:ab}. 
From \eqref{f:113a}, for  every $\e >0$ such that $x + \e < 1$ and $x - \e >0$,  we have
$$\begin{aligned}
& \int_x ^ 1   \rho \leq  \int _ x ^ 1 b  e ^ { a t } + \int _{ x + \frac{\e}{2}} ^ {x + \e } \big ( \rho - b  e ^ { a t } \big ) 
\\ 
& \int _0 ^ x   \rho \leq  \int _ 0 ^ x  b e ^ { a t } + \int _{ x -  \e }^ {x  - \frac{\e}{2} } \big ( \rho - b e ^ { a t } \big )  \, , 
 \end{aligned}
 $$ 
and hence 
\begin{equation}\label{f:pro} J _ \rho ( x) \leq \frac{1}{\rho ( x)} 
\left  \{   
\frac{1}{ \int _ x ^ 1 b  e ^ { a t } 
 \Big  [ 1 +  
\frac {\int _{ x + \frac{\e}{2}} ^ {x + \e } \big ( \rho - b  e ^ { a t } \big ) } {\int _ x ^ 1 b  e ^ { a t } } 
 \Big ]   } 
 + 
 \frac{1}{\int _ 0 ^ x b e ^ { a t } 
 \Big  [ 1 +  
\frac {\int _{ x - \e } ^ {x - \frac{\e}{2}  } \big ( \rho - b  e ^ { a t } \big ) } {\int _ 0 ^ x b  e ^ { a t } } 
 \Big ]   } 
 \right \}  ^ { -1} \,.
 \end{equation} 
We claim that,  
setting for brevity 
$$  c_\delta :=    \min _{[\delta, 1 - \delta]  }  \Big ( \frac{h'}{h} \Big ) ^ 2 \,, $$   
for a positive  constant $\Lambda_1= \Lambda_1( \delta, m)$ it holds
\begin{equation}\label{f:claim}   
\frac {\int _{ x + \frac{\e}{2}} ^ {x + \e } \big ( \rho - b e ^ { a t } \big ) } {\int _ x ^ 1 b  e ^ { a t } } 
  \leq  - \Lambda_1  c_\delta \, , \qquad \frac {\int _{ x - \e } ^ {x - \frac{\e}{2}  } \big ( \rho - b  e ^ { a t } \big ) } {\int _ 0 ^ x b  e ^ { a t } }  \leq - \Lambda_1   c_\delta\,.
  \end{equation} 
Once proved this claim,  the conclusion readily follows. Indeed, we have \begin{equation} \label{f:jump} J _ \rho ( x)  \leq \frac{1}{\rho ( x)}  
 \left \{ \frac{1}{ \int _ x ^ 1 b  e ^ { a t } } + \frac{1}{\int _ 0 ^ x b  e ^ { a t }  } \right \} ^ { -1}  ( 1 - \Lambda _1 c_\delta)  
 \leq x ( 1-x )  ( 1 - \Lambda _1 c_\delta)  \leq x ( 1-x) - \Lambda _0  c _\delta  \,,  
  \end{equation}  
 where the first inequality follows from \eqref{f:pro} and \eqref{f:claim}, the second inequality follows from part (i) of the statement already proved, and the last inequality holds by taking  $\Lambda _0  =  ( 1- 2\d) ^ 2 \Lambda _1  $.  

It remains to prove \eqref{f:claim}. Let us focus on the first of the two inequalities (the other one being analogous), that we rewrite as 
\begin{equation}\label{f:claim1}
\frac {\int _{ x + \frac{\e}{2}} ^ {x + \e } \big ( b e ^ { a t }  - \rho \big ) } {\int _ x ^ 1 b  e ^ { a t } } 
  \geq  \Lambda  _1  c_\delta\,.
  \end{equation} 
 By  the current assumptions on $\rho$,  
we have 
$$g''  = (\log \phi) '' +  (\log h ) ''   \leq  (\log h ) '' = -  \Big ( \frac{h'}{h} \Big ) ^ 2  \,$$ 
Therefore, 
for every $t \in (\d, 1 - \d)$, we have,  for some $\xi = \xi _{ x, t }$ in  the interval with endpoints $x$ and $t$ 
$$ \begin{aligned}
g ( t)  & = g ( x)+ g' ( x) ( t-x) + \frac{1}{2} g'' (\xi)( t-x ) ^ 2 
\\
& \leq  g ( x)+ g' ( x) ( t-x) - c_\delta  ( t-x ) ^ 2 \,.
\end{aligned} 
$$ 
Hence  the weight $\rho$ satisfies the estimate 
\begin{equation} \label{f:113b}
 \rho ( t) \leq b  e ^ { a  t } e ^ { - \frac{c_\delta }{2} ( t-x) ^ 2 } \qquad \forall  t \in (\delta, 1 - \delta )\,,
  \end{equation}
Since  $x \in ( 2 \d, 1 - 2 \d)$, choosing $\e = \d$ we have that   $( x+ \frac{\e}{2}, x + \e) \subset ( \delta, 1 - \delta) $. Then   integrating the above inequality
over $( x+ \frac{\d}{2}, x + \d)  $ we obtain 
 $$\begin{aligned} 
  \frac {\int _{ x + \frac{\d}{2}} ^ {x + \d } \big (  b e ^ { a t } - \rho \big ) } {\int _ x ^ 1 b  e ^ { a t } }  & \geq  
  \frac {\int _{ x + \frac{\d}{2}} ^ {x + \d }   b e ^ { a t }    \big ( 1-e ^ { - \frac{c_\delta }{2} ( t-x) ^ 2 }   \big ) } {\int _ x ^ 1 b  e ^ { a t } }  
    =   \frac{a}{e ^ a - e ^ { ax}  }  
   {\int _{ x + \frac{\d}{2}} ^ {x + \d }    e ^ { a t }    \big ( 1-e ^ { - \frac{c_\delta }{2} ( t-x) ^ 2 }   \big ) }  
   \\ \noalign{\bigskip}
&   \geq    \frac{a e ^ { ax} } {e ^ a - e ^ { ax}  }  
   {\int _{ x + \frac{\d}{2}} ^ {x + \d }     \big ( 1-e ^ { - \frac{c_\delta }{2} ( t-x) ^ 2 }   \big ) }
   \geq    \frac{\d }{2}  \frac{a e ^ { ax} } {e ^ a - e ^ { ax}  }  
     \big ( 1-e ^ { -  {c_\delta } \frac{\d  ^ 2}{8}}   \big )    
     \\ \noalign{\bigskip}  
&     \geq 
      \frac{\d }{2}  \frac{a e ^ { 2 \d a} } {e ^ a - e ^ { 2 \d a}  }  
     \big ( 1-e ^ { -  {c_\delta } \frac{\d  ^ 2}{8}}   \big )   \geq \Lambda_1  c_\d \,,
         \end{aligned}
   $$ 
where the last inequality is satisfied  for some positive constant $\Lambda _1= \Lambda _1( \d, m)$. This follows from the fact that
there exists   positive constants
$\Lambda_2 = \Lambda_2 ( \d)$  and $\Lambda_3 = \Lambda _3 (\d, m)$ such that 
 \begin{equation}\label{f:dueineq}  \big ( 1-e ^ { -  {c_\delta } \frac{\d  ^ 2}{8}}   \big )  \geq \Lambda_2 c_\d   \qquad 
 \text{ and } \qquad \frac{\d }{2}  \frac{a e ^ { 2 \d a} } {e ^ a - e ^ { 2 \d a}  }   \geq \Lambda_3  \,.\end{equation} 
 To show the first inequality it is enough to observe that the function 
 $$\varphi( c_\d):=   \frac{ \big ( 1-e ^ { -  {c_\delta } \frac{\d  ^ 2}{8}}   \big ) }{c_\d}$$ 
 admits a positive finite limit as $c _ \d \to 0 ^+$, and hence it has a positive minimum when $c_\d$ varies in its range 
 of possible values (which go from $(c_\d ) ^{min}= 0$  to some  $ (c_\d ) ^{max} >0$).

To show the second inequality notice that,  by the assumptions made on $\phi$ and $h$, the function $\rho$ turns out to be $\big ( \frac{1}{m+1} \big)$-concave, so that $\rho$ and $\rho'$ satisfy 
$$\min \{ x, 1-x \} ^ {m+1} \leq \rho ( x) \leq  1 \,,  \qquad  - \frac{m+1}{x} \leq \rho ' ( x) \leq \frac{m+1}{x}  \qquad \forall x \in (0, 1) \,.$$

Then, since our fixed point $x$ belongs to $(2\d, 1 - 2 \d)$,  we have,  
  for some positive constant $\Lambda_4 =\Lambda _4 ({\d, m})$  
$$|a| = \Big |  \frac{\rho' (x)} {\rho (x)}  \Big |  \leq  \Lambda_4 ({\d, m}) \,,$$ 
 which implies the validity of the second inequality in \eqref{f:dueineq} for some positive constant $\Lambda _ 3 =\Lambda _ 3 (\d, m)$. 

\bigskip
(iii)   By the inequality \eqref{f:parziale} proved above (for $a$ as in \eqref{f:ab}),  we know that 
 $$ J _\rho ( x) \leq J _ {e ^ { at}}  (x) \,.  $$ 
   We are thus reduced to show that, 
   for every $\d \in (0, 1)$
there exists a positive constant $\Lambda_ 0' = \Lambda _0' (\d, m)$ such that
     $$J _ {e ^ { at}}  (x)  \leq \frac{1}{4} - \Lambda_0' a^2  \qquad \forall x \in ( \d, 1- \d)   \,.
 $$ 
 
For $a = 0$, the above inequality is  trivially satisfied because we have
$$J _ {e ^ { at}}  (x) = J _ 1 ( x) = x ( 1-x) \leq \frac{1}{4}  \qquad \forall x \in [0, 1]\,.$$
For $a \neq 0$, it can be rewritten as
\begin{equation}\label{f:goal} F (a, x):=
 \frac{e ^ a - e ^ { ax} - e ^ { a ( 1-x)} + 1 }{a (e ^ a -1) }  \leq \frac{1}{4} - \Lambda_0 '  a ^ 2  \,.
 \end{equation} 
Via elementary computations  one gets that
$$F  \big ( 0, \frac{1}{2} \big ) = \frac{1}{4}\,, \qquad \nabla F \big  ( 0, \frac{1}{2} \big )  = 0 \,, \qquad 
\nabla^ 2  F \big  ( 0, \frac{1}{2} \big  ) \xi \cdot \xi <0 \quad \forall \xi \in  \R ^ 2 \setminus \{ 0 \} \,.$$ 
Hence the inequality \eqref{f:goal} holds true, with a suitable $\Lambda_0' >0$, for  $(a, x)$ in a neighbourhood of $ \big ( 0, \frac{1}{2} \big )$. 

Assume now that statement (iii) is false for some fixed $\d \in (0, 1)$.  Then along a  sequence $(a_n , x_n)$, 
with $x_n \in (\d, 1 - \d)$, 
 we have $F (a_n , x_n ) \to \frac{1}{4}$.  

Similarly as in the proof of statement (ii) above, the power-concavity of $\rho$ and the assumption $x _n \in (\d, 1 - \d)$ implies that
the sequence  $|a_n| =  |\frac{\rho' (x_n) } {\rho ( x_n) } |$  is bounded from above (by a constant depending only on $\d$ and $m$). 
Hence,   
possibly passing to a subsequence, we get  a pair $(\overline a, \overline x)$ 
such that
$F ( \overline a, \overline x) = \frac{1}{4}$. 

Recalling that $F ( \overline a, \overline x) \leq \overline x ( 1 - \overline x)$, we infer that $\overline x = \frac{1}{2} $. 

Since we have already shown the validity of the inequality \eqref{f:goal}  for  $(a, x)$ in a neighbourhood of $ \big ( 0, \frac{1}{2} \big )$, we infer that $\overline a \neq 0$. 
But  the equality conditions in statement (i) already proved ensure in particular the validity of the {\it strict} inequality $$F  ( \overline a, \frac {1}{2}) - \frac{1}{4} = J _ { e ^ {\overline a t} } (x) - x ( 1-x) \Big | _{ x = \frac{1}{2}} < 0 \qquad \text{  for any  } \overline a \neq 0 \,,$$ 
so we have reached a contradiction. \qed  

\bigskip

The second lemma needed for the proof of Proposition \ref{p:1d} reads:

\begin{lemma}\label{l:nondeg}  There exist  constant $\d _ 0 \in (0, \frac 1 2)$ and  $  \lambda _0 >0  $ such that,
for any positive  weight    $\rho$ on $I= (0, 1)$  which is $(\frac{1}{m})$-concave
 for some $m \in \N \setminus \{0 \}$, if 
$u \in  W ^ { 1, 1} _\rho  (I)$ is  such that
$$\int_0 ^ 1 u \rho = 0  \qquad  \text{ and } \qquad  \frac {\int _{0} ^  { 1}  |u'| \rho } {\int_ 0 ^ 1  |u| \rho } \leq 3 \,, $$
 then 
 $$ 
 \frac {\int _{\d_0} ^  { 1 - \d_0}   |u'| \rho   } {\int_ 0 ^ 1  |u| \rho } \geq \lambda  _0 \,.
$$ 
  \end{lemma}

\proof Let $\d \in (0, \frac 1 4 ]$, whose precise value will be fixed later. 
Assume by contradiction that there exists a sequence  of  positive $(\frac{1}{m})$-concave weights $\rho _h$ on $I$, and a sequence of functions
$u _h \in  W ^ { 1, 1} _{\rm loc} (I)$  such that 
$$\int_0 ^ 1 u \rho_h = 0 \, , \qquad \frac {\int _{0} ^  { 1}  |u' _h | \rho _h  } {\int_ 0 ^ 1  |u_h | \rho _h } \leq 3\,, \qquad  
  \lim _h \frac {\int _{\d } ^  { 1 - \d}   |u _h'| \rho  _h   } {\int_ 0 ^ 1  |u _h | \rho  _h} =0  \,.
$$  
 Without loss of generality, we can assume that $\| \rho _h \| _\infty = 1$ and $\int _0 ^ 1 |u _h | \rho _h = 1$.

 	Let  ${\Omega_h} \sq \R^{1+m}$  be the open convex set defined by
	 $$
	 {\Omega_h } = \left\{ (x,y) \in I \times \R^{m} : \, x \in I ,\,  \norma{y}_{\R^m} <   \om _m ^ { -\frac{1}{m}}  \rho^{ \frac 1 m}(x)\right\}\,,
	 $$
	 and consider the sequence of functions 
	 $$
	U_h(x,y):= u_h(x) \qquad \forall (x, y ) \in  \Om_h \,.$$
 The functions $U _h$ belong to  $W ^ { 1, 1} (\Om _h)$ and satisfy 
$$\int_{\Om _h}   U _h = 0 \, , \qquad  \int_ { \Om _h}  |U _h |   = 1  \, , \qquad  \int _{\Om _h} |  \nabla     U _h |    \leq 3\,, \qquad    
  \lim _h  \int _ {\Om _h \cap ( (\d, 1 - \d) \times \R ^m ) }  \!\!\!\!\!\!\!  |  \nabla    U _h|   =0  \,.
$$   
We observe that, thanks to the power concavity of $\rho _h$,  the  inradii and the outradii of the domains $\Om _h$ are respectively uniformly bounded respectively from below and from above, and hence the trace operators from $W ^ { 1, 1} (\Om _h)$ to $L ^ 1 (\partial \Om _h)$ have uniformly bounded norms. 
Combined with the fact that $\int _{\Om _h} | \nabla \EEE  U _h |    \leq  3$, this ensures that 
the extensions of the functions $U _h$ to zero outside $\Om _h$ (that we still denote by $U _h$) are uniformly bounded in $BV ( \R ^ {1 + m})$. 

Therefore, up to subsequences, $U _h$ converge in $L ^ 1 ( \R ^ { 1 +m})$ to a function $U\in BV ( \R ^ {1 + m})$ such that, denoting by 
$\Om$ the Hausdorff limit (up to subsequences) of the convex domains $\Om _h$,
$$\int_{\Om }   U  = 0 \, , \qquad  \int_ { \Om }  |U  |   = 1  \, , \qquad  \int _{\Om } |D U  |    \leq 3\,, \qquad    
 \int _ {\Om  \cap ( (\d, 1 - \d) \times \R ^m ) }  \!\!\!\!\!\!\!  |D U |   =0  \,.
$$   

In particular, the latter condition implies that there exists a constant $\lambda \in \R$ such that
$U (x, y ) = \lambda$ for every $(x, y) \in \Om  \cap ( (\d, 1 - \d) \times \R ^m )$. 
Up to changing $U$ by $-U$, we shall assume with no loss of generality that $\lambda \geq 0$. 
In the sequel, we set for brevity 
$$\Om _\ell:= \Om \cap ( (0, \d) \times \R ^m)  \, , \qquad  \Om _c:= 
\Om  \cap ( (\d, 1 - \d) \times \R ^m ) \, , \qquad  \Om _r:=   \Om \cap ( (1- \d, 1) \times \R ^m)  \,.$$ 
We observe that the Lebesgue measure of $\Om _c$ satisfies
$$\om _m \geq |\Om _ c| \geq \beta    _{m, \d} \geq\beta   _{m, \frac{1}{4}} \,,$$  
where the first inequality holds   because $\rho _h \leq 1$, and the second one holds (for some positive constant $ \beta  = \beta  _ {m, \d}  $ monotone decreasing in $\d$) by the  power concavity of $\rho _h$. 

We denote by  $ {\Om _\ell} ^ *$ the reflection of $\Om _\ell$ about the hyperplane $\{ \d \} \times \R ^m$,  and we set
$$\om : =  {\Om _\ell}  \cup  {\Om _\ell ^ *}  \, , \qquad \Gamma :=  \Om \cap (\{ \d \} \times \R ^m ) = \partial \Om _\ell \cap \Om  \,.$$

We now take $\e  >0$, whose precise value will be fixed later, and we distinguish two cases. 

\medskip
{\it Case 1:} $\lambda \geq \e $.  
Since
$$ \begin{aligned}
 0 = \lambda |\Om _ c| + \int _{\Om _\ell} U +  \int _{\Om _r} U & \geq   \e  \beta  _{m, \frac{1}{4}} + \int _{\Om _\ell} U +  \int _{\Om _r} U
 \\ 
 &   \geq   \e  \beta  _{m, \frac{1}{4}} + \int _{\Om _\ell} U ^ - +  \int _{\Om _r} U^-   \EEE  \,, 
 \end{aligned} 
  $$ 
   where $U ^ - :=  \min \{ U, 0 \}$,  
we can assume with no loss of generality that 
$$ \int _{\Om _\ell}  -   U^-      \geq   \frac{1}{2} \e   \beta _{m, \frac{1}{4}} \,.$$ 

We 
consider the function $v $ defined  on $\om$ by 
$$ v (x, y):= \begin{cases} 
U  ^- ( x, y) & \text{ if } (x, y) \in \Om _\ell 
\\ 
 -U ^- (2 \d -  x, y) & \text{ if } (x, y) \in \Om _\ell   ^* \,,
\end{cases}
$$ 

By construction, $v$ has zero mean on the 
the   one dimensional   horizontal slices of the set $\om$  (which have length less than or equal to $2 \d$, the horizontal width of $\om$), so that the Poincar\'e inequality    in $L ^ 1$   applied on such slices yields
 $$\int _\om |Dv| \geq \frac{1}{ \d} \int _\om |v|\,.$$

We have 
$$\int _\om |v| =  2 \int _{ \Om  _\ell}  - U ^ -   \geq \e  \beta  _{m, \frac{1}{4}} \,;$$ 
on the other hand, since $v = a$  on  the left side of  $\Gamma$ for some constant $a \leq 0$,  we have
 $$\int _\om |D v|   \leq   2 \Big [ \int _{\Om _\ell}  |DU ^- |   -  a \mathcal H ^ {  m  } (\Gamma) \Big ]  
=  2  \Big [ \int _{\Om _\ell}  |DU ^- | +\int _{\Gamma}  |DU ^- |  \Big ]  \leq 2  \Big [ \int _{\Om _\ell}  |DU | +\int _{\Gamma}  |DU  |  \Big ]    \leq 6 \,,
  $$ 
where the first inequality follows from \cite[Theorem 3.84]{AFP}.
Hence  
\begin{equation}\label{f:assurdo1} 
 6 \geq   \frac{1}{ \d} \e  \beta _{m, \frac{1}{4}}\,. 
 \end{equation}

\medskip
{\it Case 2:} $0 \leq \lambda < \e $.  
Since
$$ 1 = \lambda |\Om _ c| + \int _{\Om _\ell} |U |  +  \int _{\Om _r}|  U |  <   \e   \om _m  + \int _{\Om _\ell} |U| +  \int _{\Om _r} |U |\,,  $$ 
we can assume with no loss of generality that 
$$ \int _{\Om _\ell}  | U | \geq   \frac{1}{2} ( 1 - \e   \om _m) \,.$$

We 
consider the function $w $ defined  on $\om$ by 
$$ w (x, y):= \begin{cases} 
U   ( x, y) & \text{ if } (x, y) \in \Om _\ell 
\\ 
 -U  (2 \d -  x, y) & \text{ if } (x, y) \in \Om _\ell   ^* \,.
\end{cases}
$$ 

Similarly as in Case 1,  by Poincar\'e inequality   in $L ^ 1$  applied on the  one dimensional horizontal   slices of $\om$ we obtain 
 $$\int _\om |Dw| \geq \frac{1}{ \d} \int _\om |w|\,,$$

We have 
$$\int _\om |w| =  2 \int _{ \Om  _\ell}   | U |   \geq  ( 1 - \e   \om _m)   \,;$$ 
on the other hand, since $w = b$  on  the left side of  $\Gamma$ for some constant $b \in \R$,  we have
 $$ \begin{aligned}
 \int _\om |D w|   &    \leq  2 \Big [ \int _{\Om _\ell}  |DU |   + |b |   \mathcal H ^ {   m  } (\Gamma) \Big ]  
=  2  \Big [ \int _{\Om _\ell}  |DU  | + ( |b | - \lambda) \mathcal H ^ {   m  } (\Gamma)   \Big ]  + 2 \lambda \mathcal H ^ {   m   } (\Gamma)  \\  
\noalign{\medskip} 
& \leq 2  \Big [ \int _{\Om _\ell}  |DU  | + ( |b | - \lambda) \mathcal H ^ {   m  } (\Gamma)   \Big ]  + 2 \e  \om _m  \\ 
\noalign{\medskip} 
&  \leq 2    \Big [ \int _{\Om _\ell}  |DU  | +   |b  - \lambda | \mathcal H ^ {  m  } (\Gamma)   \Big ]  + 2 \e  \om _m \\ 
\noalign{\medskip}  
&
=  2  \Big [ \int _{\Om _\ell}  |DU  | + \int _{\Gamma}  |DU |   \Big ]  + 2 \e  \om _m  \leq  6  + 2 \e  \om _m \,.
 \end{aligned} 
  $$ 
Here,  similarly as in Case 1,  we have used   \cite[Theorem 3.84]{AFP} for the first inequality,   and the estimate $\mathcal H ^ m (\Gamma) \leq \omega _m$ in the second inequality;    moreover, in the last line we have used the equality 
$  |b  - \lambda | \mathcal H ^ { m  } (\Gamma)  = \int _{\Gamma}  |DU |$, which holds by \cite[Theorem 1, Section 5.4]{EG}. 
We conclude that  
\begin{equation}\label{f:assurdo2} 
 6 + 2 \e  \om _m  \geq   \frac{1}{ \d} ( 1 - \e   \om _m)    \,. 
 \end{equation}

We now fix the constants $\delta$ and $\e $ so that both \eqref{f:assurdo1} and \eqref{f:assurdo2} lead to a contradiction. 
 
The inequality \eqref{f:assurdo2} is false provided $\delta \leq \frac{1}{8}$ and $\e  < \frac{1}{5 \om _m }$. 

We fix $\e  = \frac{1}{6 \om _m}$. For such value of $\e$, inequality \eqref{f:assurdo1} is false provided 
$\delta < \frac{1}{36 \om _m} \beta  _{m, \frac{1}{4}}$. 
Then our proof is achieved by taking
$$\delta_0 := \frac{1}{2}  \min \big \{ \frac{1}{8} , \frac{1}{36 \om _m} \beta _{m, \frac{1}{4}} \big \}\,.$$

 \qed
 
\bigskip   
\underbar {Proof of Proposition \ref{p:1d}}.  Let $\{u _h \}  \in W ^ {1, 1 }  _\rho (0, 1)$, with $\int_ 0 ^ 1  |u _h| \rho = 0$, be a minimizing sequence for $\s (I, \rho)$, and let $\delta _0  \in (0, \frac 1 2)$ be given by Lemma \ref{l:nondeg}. 
By inequality \eqref{f:nuce} and Lemma \ref{l:1refinement} we have 
  $$\begin{aligned} \int_ 0 ^ 1  |u _h| \rho  & \leq 2 \int _{0} ^  { 1}   |u _h '| \rho  J _\rho  = 2 \Big [ 
 \int _ 0 ^ {\delta_0}    |u _h '| \rho  J _\rho   + 
 \int _  { 1 - \delta _0}    ^1   |u _h '| \rho  J _\rho +   \int _  {\delta_0} ^ { 1 - \delta _0}     |u _h '| \rho  J _\rho    \\ 
 & \leq 2 \Big \{  
  \frac{1}{4} \int _ 0 ^ {\delta_0}    |u _h '| \rho    + 
   \frac{1}{4} \int _  { 1 - \delta _0}    ^1   |u _h '| \rho  
 + \Big [ \frac{1}{4}  - \Lambda   _0 \min _{[\d_0, 1 - \d_0]} \Big ( \frac{h'}{h} \Big ) ^ 2  \Big ] 
  \int _  {\delta_0} ^ { 1 - \delta _0}     |u _h '| \rho   \Big \} 
   \\ 
 &  =  \frac{1}{2} \int _ 0 ^ 1    |u _h '| \rho  - 2\Lambda  _0  \Big [ \min _{[\d_0, 1 - \d_0]} \Big ( \frac{h'}{h} \Big ) ^ 2  \Big ]   \int _  {\delta_0} ^ { 1 - \delta _0}     |u _h '| \rho \,. 
   \end{aligned} 
$$ 
Then,  by Lemma \ref{l:nondeg},  we conclude that 
$$\begin{aligned} \frac {\int _{0} ^  { 1}   |u _h '| \rho  } {\int_ 0 ^ 1  |u _h| \rho }   & \geq  2 + 4 \Lambda _0 \Big [  \min _{[\d_0, 1 - \d_0]} \Big ( \frac{h'}{h} \Big ) ^ 2   \Big ] 
\frac{ 
 \int _  {\delta_0} ^ { 1 - \delta _0}     |u _h '| \rho}  {\int_ 0 ^ 1  |u _h| \rho }
 \\ 
  & \geq  2 + 4  \Lambda  _0 \lambda _0  \Big [  \min _{[\d_0, 1 - \d_0]} \Big ( \frac{h'}{h} \Big ) ^ 2   \Big ] \,.  \end{aligned} 
  $$ 
  So Proposition \ref{p:1d} is proved, with $C_1 =  4 \Lambda _0 \lambda _0$. 
\qed

\section{Proof of Theorem \ref{t:constwidth}}\label{sec:teo2} 
  Throughout this section (and the next one), we are going to exploit the equality $\sigma _ 1 (\Om) = \sigma _ 1 ^* (\Om)$ given by Proposition \ref{p:ad}. 
For brevity,  and since the two quantities coincide, we always write $\sigma _ 1 (\Om)$  in place of $\sigma ^* (\Om)$.   
 
Let us start by recalling  that, for every open bounded convex set $\Om$ in $\R ^2$, $\sigma_1 (\Om)$ is attained at some set $E$. Note that a minimizer $E$ is a connected isoperimetric region, so that the arc $\partial E \cap \Om$ has constant mean curvature, being either a line segment or an arc of circle.  Moreover, the contact points between $\partial E \cap \Om$ are points
 where $\partial \Om$ is differentiable. We refer to \cite[Theorems 2 and 3]{cianchi89BUMI} for the analogue properties in  case of $\mu_1 (\Om)$. The arguments extend naturally to  $\sigma_1 (\Om)$. Next we observe that, setting
 $$\mathcal W_1 := \text{ the family of open  convex sets $\Om$ in $\R ^2$ of constant width equal to $1$},  $$ 
 the statement is equivalent to
 $$\sigma _ 1  (\Om) \geq  \sigma _1 (B _ { \frac 1 2 } ) = \frac 8 \pi\,,$$ 
 with equality sign if and only if $\Om = B _ {\frac 1 2 }$.  
Notice that   
a solution to   the minimization problem 
 \begin{equation}\label{f:mincw} 
 \min \Big \{ \sigma _ 1 (\Om) \ : \   \Om  \in \mathcal W _ 1  \Big \} \,, 
 \end{equation}  
exists   because  (in the Hausdorff topology)  the map $\overline \Om \mapsto \sigma _1 (\Om)$ is continuous and the class of competitors is compact. 
    Then we are going to argue by contradiction.  
    
    \smallskip 
 Our proof is built upon the  three lemmas:
 in Lemma \ref{l:symmetry} we obtain some preliminary geometric information about a minimizer $\Om$ which is assumed to be different from $B _ \frac{1}{2}$; 
  then in Lemma \ref{l:nor0} and Lemma \ref{l:blaschke}, 
we extract some   additional  information on such a minimizer via the use of different kinds of perturbations.  
 \smallskip 
 
 In this respect we point out that  perturbing sets of constant width in order to extract some information from optimality is a classical question. We refer to \cite{Ha02} for an early work on the variational approach to the Blaschke-Lebesgue theorem, and to 
  \cite{Bo24nv, HL21, HL24} and  references therein  for some more recent related works. 
  Our strategy    is, up to our knowledge, new   and
of  independent interest for its possible application in other optimization problems regarding convex bodies of constant width. 
 
   The  key novelty is a delicate singular perturbation argument for convex bodies of constant width, that we present  as a separate result
to highlight its generality, see Proposition \ref{p:cruciall}. This  enables us to exclude,  in Lemma \ref{l:nor0},  
the  presence  in $\mathbb S ^ 1$ of  directions  where the radius of curvature $r_\Om$ of $\Om$ vanishes in the Lebesgue sense. We recall that, for an arbitrary body of constant width, 
this set may be even a dense set, see \cite[Section 5]{kallay}.

   Then,  relying in essential  way on Lemma \ref{l:nor0},  and  using Minkowski-type perturbations somehow in the spirit of \cite{Ha02},   in Lemma \ref{l:blaschke} we arrive at writing the boundary of $\Om$ as the union of $4$ arcs whose radii of curvature are explicitly written in terms of the Lebesgue measures of an optimal set $E$ for $\sigma _1 (\Om)$ and of its complement.    Eventually,  we conclude that such configuration of $4$ arcs   cannot occur, except for the ball,  because 
    we reach the contradiction that   $r_\Om$ would not be balanced.

 \smallskip 
Let us start with the first lemma:

 \begin{lemma}\label{l:symmetry} 
Assume that problem \eqref{f:mincw}  is not solved by $B _ \frac{1}{2}$. Then it admits a minimizer $\Om$ such that,  if   $E$ is an optimal set for $\sigma _ 1 (\Om)$, and  $S$ is the strip bounded by two tangent lines to $\partial \Om$ at the endpoints of $E$,   we have that: 
 \begin{itemize}
 \item[(i)] $\Om$ is symmetric about the bisector of $S$;
 \smallskip
 
 \item[(ii)] $\partial E \cap \Om$ is   an arc of circle of strictly positive curvature.  
  \smallskip

  \end{itemize}
 \end{lemma}

 \proof (i)  Let  $\Om$  be a solution to problem \eqref{f:mincw}. 
Let $\Om _S$ be  the Steiner symmetrization  of $\Om$ with respect to the bisector of $S$.
Then there exists  a convex set  $\Om_*$  of constant width equal to $1$ which contains $\Om_S$ and is symmetric about the bisector of  
$S$   (see \cite[Theorem 1]{Rogers71}).  We claim that  \begin{equation}\label{f:doppiaug0}
\sigma_ 1 (\Om)= \sigma _ 1 (\Om  _S) = \sigma _ 1 (\Om  _*)\,, 
  \end{equation}  
so that $\Om_*$  is a solution to problem \eqref{f:mincw}, symmetric about the bisector of $S$.

\smallskip
To prove the first equality in \eqref{f:doppiaug0},  we observe that 
the arc $\gamma: = \partial E \cap \Om$ splits the set $\Om  _S$ into  a set  $E'$ 
and its complement $\Om  _S \setminus  E'$, with 
$| E '| = |E|$ and  $|\Om   _S \setminus E' | = |\Om  \setminus E|$. 
Hence we have
$$\sigma _ 1 (\Om  _S) \leq \frac{1}{2} {{\rm Per} ( E', \Om _S)}  \Big ( { \frac{1}{| E ' | } + \frac{1}{|\Om _S  \setminus  E ' | } }  \Big )
= \frac{1}{2} {{\rm Per} ( E, \Om )} \Big ( { \frac{1}{| E| } + \frac{1}{|\Om  \setminus E | } }   \Big ) = \sigma _ 1 (\Om) \,. $$ 
 By the properties of Steiner symmetrization, we have $D _{\Om  _S} \leq D _{\Om}   ( = 1)$. Then, 
up to replacing $\Om _S$ by its dilation by a factor larger than $1$, we may assume that $D _{\Om  _S }= 1$. 
Therefore, the inequality $\sigma _ 1 (\Om  _S) \leq \sigma _ 1 (\Om)$, combined with the optimality of $\Om$, implies that 
$\sigma_ 1 (\Om)= \sigma _ 1 (\Om  _S)$, and that $E'$ is an optimal set for $\sigma _ 1 (\Om _S)$. 
To prove the second equality in \eqref{f:doppiaug0},  we proceed in a similar way, namely we observe that
the arc $\gamma$ splits the set $\Om  _*$ into  a set  $E''$ 
and its complement $\Om  _* \setminus  E''$, with 
$| E ''| \geq |E '|$ and  $|\Om   _* \setminus E'' | \geq   |\Om   _S \setminus E' |$, so that  
$$ \sigma _ 1 (\Om _*) \leq \frac{1}{2} {{\rm Per} ( E'', \Om _*)}  \Big ( { \frac{1}{| E '' | } + \frac{1}{|\Om _*  \setminus  E '' | } }  \Big ) \leq \frac{1}{2} {{\rm Per} ( E', \Om _S)}  \Big ( { \frac{1}{| E ' | } + \frac{1}{|\Om _S  \setminus  E ' | } }  \Big ) = 
\sigma _ 1 (\Om _S) \,. $$  
From the optimality of $\Om _S$, we infer that $\sigma_ 1 (\Om_*)= \sigma _ 1 (\Om  _S)$, and that $E''$ is optimal for $\sigma_ 1 (\Om_*)$. 
\smallskip

\smallskip
(ii) Let  $\Om$  be a solution to problem \eqref{f:mincw}, symmetric about the bisector of $S$, as found in the previous item.  Assume by contradiction that $\partial E \cap \Om$  is contained into a straight line $\sigma$, letting $\Om _\sigma$ 
be  the Steiner symmetrization  of $\Om$ with respect to $\sigma$ (possibly dilated so that it has diameter $1$), 
by arguing as in the proof of item (i) we obtain that a 
convex set of constant width equal to $1$ which contains $\Om_\sigma$ is still a solution  to problem \eqref{f:mincw}.  But such a set, having constant width and two orthogonal axes of symmetry, is necessarily a ball.
Indeed, we notice  that all diameters have to pass through the origin, which is a center of symmetry
(if by contradiction a diameter does not pass through the origin, its symmetric segments with respect to the two orthogonal axes of symmetry would be two parallel diameters, which is impossible).
 Since all the diameter have the same length, the set is a ball, 
against the assumption that problem \eqref{f:mincw} is not solved by $B _ \frac{1}{2}$.  
  \qed

 \bigskip\smallskip
Now, in order to pursue our proof by contradiction, we are going to obtain some information on a solution $\Om$ to problem
 \eqref{f:mincw} of the form found in Lemma \ref{l:symmetry}.  

Given such a minimizer, we call
 {\it secured points} of $\partial \Om$  the family of $4$ points given by the two endpoints of the arc $\partial E \cap \Om$, and their diametral points in $\Om$.  (Recall that, in any set of constant width, each boundary point is an extremum of a diameter. We say that two points are diametral  if they are extrema of the same diameter.)

We denote by $\nu _\Om (\cdot)$  the (possibly multivalued)  Gauss map of $\partial \Om$.
  
We consider the subset of $\mathbb S ^ 1$ defined as the spherical image through $\nu _\Om$ 
of $\partial \Om$ deprived of the 
secured points, namely: 
 \begin{equation}\label{f:erre} 
 \mathcal R (\Om)  : = \big \{ \xi  \in \mathbb S ^ 1 \ :\ \xi  \in  \nu _ \Om ( x) \text{ for some  $x\in \partial \Om$ which is not a secured point of $\partial \Om$}  \big \} \,,
 \end{equation} 

 \smallskip 
Recall that, since $\Om$ has constant width, its surface area measure 
is absolutely continuous with respect to $\mathcal H ^1 \res \mathbb S ^ 1$ (see \cite[Theorem 3.7]{HugII})    
and its density is given by the curvature function 
$r_\Om := h _ \Om + h _\Om ''$, where $h _\Om $ denotes the support function of $\Om$, which is of class $C ^ {1,1}$ 
 (see \cite[Corollary 2.6]{Howard}).  

Then we decompose $\mathcal R(\Om)$  as  the following disjoint union:  
   $$\mathcal R(\Om)  =  \mathcal R _0 (\Om) \cup \mathcal R _1(\Om)   \cup \mathcal R _*(\Om)  \,,$$
 where 
 \begin{equation}\label{f:raggi}
 \begin{aligned}
 \ & \mathcal R _0(\Om)   = \{ \xi \in \mathcal R(\Om)  \ :\ r _\Om (\xi ) =0 \}  \\ 
 & \mathcal R _1 (\Om)  = \{ \xi \in \mathcal R(\Om)  \ :\ r _\Om (\xi ) =1 \}  \\
  & \mathcal R _*(\Om)   = \{ \xi \in \mathcal R(\Om)  \ :\ r _\Om (\xi ) \in ( 0, 1) \} \,.
 \end{aligned} 
\end{equation}

In the definition of $\mathcal R _0(\Om) $, the  equality $r_\Om (\xi ) = 0$ must be intended as follows: $\xi$ is a Lebesgue point for $r_\Om$, such that the approximate limit of $r_\Om$ at $\xi$ equals $0$ (and similarly in  definitions of $\mathcal R _1(\Om) $ and $\mathcal R _*(\Om) $). 
\smallskip

 \begin{lemma}\label{l:nor0}   
 Assume that problem \eqref{f:mincw}  is not solved by $B _ \frac{1}{2}$. 
  Let $\Om$ be a solution to the minimization problem \eqref{f:mincw} as given by Lemma \ref{l:symmetry}, and let 
 $\mathcal R _0(\Om) $ and $\mathcal R _ 1 (\Om)$ be the corresponding sets of directions $\xi\in \mathbb S ^ 1$ defined by \eqref{f:raggi}. Then
$$  \mathcal R _0(\Om)   =   \mathcal R _1(\Om)     = \emptyset  \,.$$ 
\end{lemma}

\bigskip 

Since the  spherical image of  the  family of secured points is  a symmetric closed  
subset of $\mathbb S ^ 1$ given by  $4$ distinct points,  as an immediate consequence of Lemma \ref{l:nor0}  
we obtain: 

\begin{corollary} Under the same assumptions of Lemma \ref{l:nor0}, it holds 
 \begin{equation}\label{f:r*}
\mathcal R  (\Om) = \mathcal R _* (\Om)=  \bigcup _{i = 1, 2}  \big (\gamma _i \cup \gamma ^s _i  \big ) \,, 
\end{equation} 
 where each $(\gamma _i , \gamma _i ^ s)$, for $i = 1, 2$, is a pair of mutually disjoint and symmetric open arcs of $\mathbb S ^ 1$ --  that we call {\it free arcs} -- whose endpoints belong to the spherical image of some secured point of $\partial \Om$. 
\end{corollary} 

\bigskip

 The proof of Lemma \ref{l:nor0} is based on the following

\begin{proposition}\label{p:cruciall} 
Let $\Om \in \mathcal W_1$.  If there exists  a Lebesgue point  $\xi _0 \in \mathbb S ^ 1$ of $r_\Om$ such that $r_\Om (\xi_0)= 0$, letting $x_0 \in \partial \Om$ be such that $\xi_0 \in \nu _\Om  (x_0)$, and  $x_1$ be its diametral point, it is possible to construct a one-parameter family of convex bodies $\Om _ \e$ belonging to $\mathcal W_1$, such that 
  $\Om \setminus \Om _ \e$ and $\Om _\e \setminus \Om$ are regions of infinitesimal area neighbouring $x_0$ and $x_1$ respectively, with 
\begin{equation}\label{f:singlim}  \lim _ { \e \to 0 }  \frac {| \Om _ \e \setminus \Om | } {| \Om  \setminus \Om_ \e |  }  = + \infty \,.
\end{equation}
In particular, we have $|\Om _ \e| > |\Om|$ for $\e$ small enough. 
\end{proposition}

\bigskip Let us first show how  Lemma \ref{l:nor0} is readily obtained from Proposition \ref{p:cruciall}, and then we turn back to the proof of the lemma.

\bigskip 
\underbar {Proof of Lemma \ref{l:nor0}}.  Assume by contradiction that $\mathcal R _0 (\Om) \neq \emptyset$.  Then there exists 
a Lebesgue point  $\xi _0 \in \mathbb S ^ 1$ of $r_\Om$ such that $r_\Om (\xi_0)= 0$,  with $\xi _0 \in \nu _\Om (x_0)$ for some 
$x_0 \in \partial  \Om$ which is not a secured point of $\partial \Om$.  Let $x_1$ denote the diametral point  of $x_0$ in $\partial \Om$. 
We consider a perturbation $\Om _ \e$ of $\Om$ in $\mathcal W _1$ as given by Proposition \ref{p:cruciall}. It induces in a natural way the perturbation $E_\e$ of $E$, which is obtained replacing  the portion of  $\partial E$ lying in $\partial \Om$ by
its perturbation lying in $\partial \Om _\e $. 
Since by construction the arc $\gamma = \partial E \cap \Om$ is left fixed by the perturbation, we have 
${\rm Per } (E _\e, \Om _\e)=  {\rm Per } (E  , \Om  )$ for every $\e$. In particular, in order to contradict the optimality of $\Om$
as a solution to 
the minimization problem \eqref{f:mincw} 
 it is enough to show that, for $\e$ small enough, it holds
\begin{equation}\label{f:divoluta}
{ \frac{1}{|E _\e| } + \frac{1}{|\Om_\e \setminus E_\e| } }  < { \frac{1}{|E | } + \frac{1}{|\Om \setminus E| } } \,. 
\end{equation} 
since this implies $\sigma _ 1 (\Om _\e) < \sigma _ 1 (\Om)$. 

 Up to exchanging $E$ with its complement, let us assume that  $E$ is convex. 
Then it may occur that either both the points $x_0$ and $x_1$ belong  to $\partial (\Om \setminus E)$, or that 
one of them belongs to $\partial E$ and the other one to $\partial (\Om \setminus E)$. 
In the first case, the required inequality \eqref{f:divoluta} is satisfied because
$$
\begin{cases} 
|\Om _ \e \setminus E _ \e | = |\Om  \setminus E  | & 
\\  \noalign{\medskip}
 |E _ \e | = | \Om_ \e | - |\Om \setminus E  | \ \Rightarrow \ |E _\e | - |E| = 
|  \Om _ \e | - |\Om|  >0 &
\,, \end{cases} $$  
where the last inequality holds for $\e$ sufficiently small by Proposition \ref{p:cruciall}. 

Assume now for definiteness that $x_0 \in \Om \setminus E$ and $x_1 \in \partial E$ (the case when  $x_1 \in \Om \setminus E$ and $x_0 \in \partial E$    is analogous).  Then we have 
$$
\begin{cases} 
 |E _ \e | - |E| =   | \Om _ \e \setminus \Om| 
& \\  \noalign{\medskip}
 |\Om _ \e \setminus E _ \e | - |\Om  \setminus E  |=   - | \Om   \setminus \Om_\e |  \,, 
&
\end{cases} $$   
We deduce that also in this case  the required inequality \eqref{f:divoluta} is satisfied because
$$ \begin{aligned}
{ \frac{1}{|E _\e| } + \frac{1}{|\Om_\e \setminus E_\e| } }   & = 
 \frac{1}{|E | }   \Big ({1+ \frac{| \Om _ \e \setminus \Om| }{|E|}}  \Big )   ^ { -1}+ 
  \frac{1}{|\Om \setminus E | } \Big ( {1- \frac{| \Om  \setminus \Om _\e| }{|\Om \setminus E|}}  \Big )  ^ { -1}
  \\ \noalign{\bigskip} 
  & \sim   \frac{1}{|E | }   \Big ({1-  \frac{| \Om _ \e \setminus \Om| }{|E|}}  \Big )   + 
  \frac{1}{|\Om \setminus E | } \Big ( {1+ \frac{| \Om  \setminus \Om _\e| }{|\Om \setminus E|}}  \Big )  
   \\ \noalign{\bigskip} 
  & =  \frac{1}{|E | }  + \frac{1}{|\Om \setminus E | }  -  \frac{| \Om _ \e \setminus \Om| }{|E| ^2} 
  + 
  \frac{| \Om  \setminus \Om _\e| }{|\Om \setminus E| ^2}   <    \frac{1}{|E | }  + \frac{1}{|\Om \setminus E | } 
  \end{aligned} 
 $$ 
where the last inequality holds for $\e$ sufficiently small by Proposition \ref{p:cruciall}.  \qed

\bigskip

\underbar {Proof of Proposition \ref{p:cruciall}}. 

\bigskip\bigskip
\bigskip\bigskip

\begin{figure} [h] 
\centering   
\def\svgwidth{7cm}   
\begingroup%
  \makeatletter%
  \providecommand\color[2][]{%
    \errmessage{(Inkscape) Color is used for the text in Inkscape, but the package 'color.sty' is not loaded}%
    \renewcommand\color[2][]{}%
  }%
  \providecommand\transparent[1]{%
    \errmessage{(Inkscape) Transparency is used (non-zero) for the text in Inkscape, but the package 'transparent.sty' is not loaded}%
    \renewcommand\transparent[1]{}%
  }%
  \providecommand\rotatebox[2]{#2}%
  \ifx\svgwidth\undefined%
    \setlength{\unitlength}{265.53bp}%
    \ifx\svgscale\undefined%
      \relax%
    \else%
      \setlength{\unitlength}{\unitlength * \real{\svgscale}}%
    \fi%
  \else%
    \setlength{\unitlength}{\svgwidth}%
  \fi%
  \global\let\svgwidth\undefined%
  \global\let\svgscale\undefined%
  \makeatother%
  \begin{picture}( 4 , .92)%
    \put(0.45, 0){
     \includegraphics[height=8cm]{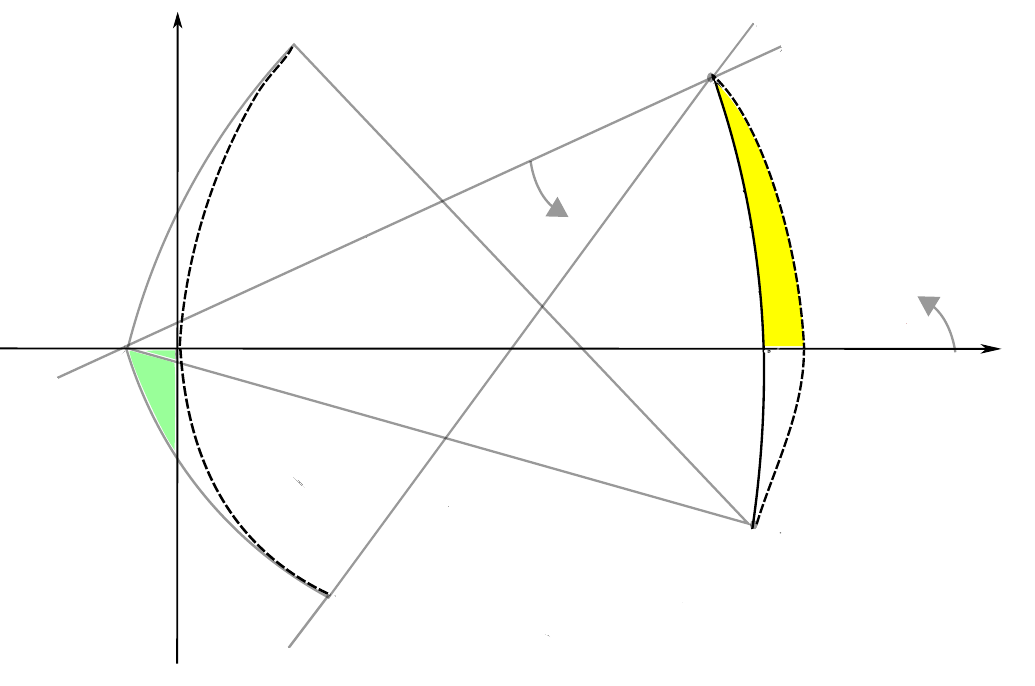}} %
         \put(0.52, 0.59){\color[rgb]{0,0,0}\makebox(0,0)[lb]{\smash{$( - \e, 0)$}}}%
  \put(1.82, 0.59){\color[rgb]{0,0,0}\makebox(0,0)[lb]{\smash{$( 1, 0)$}}}%
  \put(1.62, 1.05){\color[rgb]{0,0,0}\makebox(0,0)[lb]{\smash{$p _ \e$}}}%
  \put(1.72, 0.21){\color[rgb]{0,0,0}\makebox(0,0)[lb]{\smash{$\tilde p _ \e$}}}%
    \put(2.1, 0.6){\color[rgb]{0,0,0}\makebox(0,0)[lb]{\smash{$\theta$}}}%
    \put(1.3, 0.81){\color[rgb]{0,0,0}\makebox(0,0)[lb]{\smash{$\omega$}}}%
    \put(0.54, 0.81){\color[rgb]{0,0,0}\makebox(0,0)[lb]{\smash{$\partial B (\tilde p _ \e; 1)$}}}%
    \put(0.54, 0.33){\color[rgb]{0,0,0}\makebox(0,0)[lb]{\smash{$\partial B (p _ \e; 1)$}}}%
    \put(1.4, 0.63){\color[rgb]{0,0,0}\makebox(0,0)[lb]{\smash{$\partial B ((- \e, 0); 1)$}}}%
      \put(1.79, 1.04){\color[rgb]{0,0,0}\makebox(0,0)[lb]{\smash{$\theta = \theta _ \e$}}}%
      \put(1.74, 1.12){\color[rgb]{0,0,0}\makebox(0,0)[lb]{\smash{$\theta = \Theta _ \e$}}}%
      \put(0.98, 0.06){\color[rgb]{0,0,0}\makebox(0,0)[lb]{\smash{$\theta = \Theta _ \e + \pi$}}}%
           \put(0.38, 0.46){\color[rgb]{0,0,0}\makebox(0,0)[lb]{\smash{$\theta = \theta _ \e + \pi$}}}%

 \end{picture}%
\endgroup%
\caption{The geometry of Proposition \ref{p:cruciall}: the dotted  curves  represent portions of  $\partial \Om$; 
the coulored regions above and below $\vec x$ represent respectively $\Delta _\e ^+$ and $\Delta _\e ^ -$. }
\label{fig:1}   
\end{figure}

Up to a rotation, let us assume that $\xi _0 = (1,0)$. Moreover, let us fix a system of coordinates with the origin at $x_1$,  and the $\vec x$ axis along the diameter joining $x_0$ and $x_1$, so that  $x_0 =  (1, 0)$. In a neighbourhood of  $x_0$, the curve $\partial \Om$ can be parametrized by arc-length 
as   $\gamma _ \Om ( t) =  x_0 + \int _0 ^ t \gamma _\Om ' ( s) \, ds$. 
Since  as already recalled above the surface area measure of $\Om$ is absolutely continuous with respect to $\mathcal H ^ 1 \res \mathbb S ^ 1$, 
we can perform a change of variables  through the Gauss map of $\Om$, in which the Jacobian is given by the curvature function $r_\Om$ (see e.g. \cite[Proposition 1]{CF24}). 
Taking into account that $\gamma _\Om '  = \nu _ \Om ^ \perp$, we obtain, for $\xi \in \nu _\Om ( t) $ in a neighbourhood of $\xi _0$ in $\mathbb S ^ 1$,
$$\gamma _ \Om ( \xi ) =  x_0   + \int _ {\xi _0} ^  \xi  \nu _\Om   ^ \perp \,  r _\Om (\nu _\Om ) \, d \nu_\Om  \,. $$  
 We are going to exploit the above  equality by components,   by writing 
$\xi  = (\cos \theta, \sin \theta)$ for $\theta$ in a neighbourhood of $0$. Setting for brevity
$r ( t):= r _\Om ( \cos t, \sin t)$, 
we have
 \begin{equation}\label{f:parangle} 
 \gamma _\Om ( \theta ) =( \gamma _\Om^1 ( \theta ), \gamma _\Om ^ 2 ( \theta)  ) =  \Big ( 1 - \int _0 ^  \theta    r (t) \,  \sin t  \, d t \,,   \int _0 ^  \theta    r (t) \,  \cos t \, d t \Big )\,.
 \end{equation} 
 
 We are now ready to construct the perturbed body $\Om _ \e$. 
Let $B ( (- \e, 0); 1)$ be the ball of centre $(- \e, 0)$ and radius $1$, parametrized by
 $$f_\e (\theta)= (f_\e ^ 1 (\theta) , f _ \e ^ 2 ( \theta)) : = ( -\e + \cos \theta , \sin \theta )\, , \qquad  \theta \in \big ( 0, 2 \pi)\,.$$ 
 
 Since $\Om\in \mathcal W_1$, it lies inside the ball of radius $1$ centred at the origin, the unique contact point being $x_0 = ( 1, 0)$. It follows that, 
 for $\e$ sufficiently small, $\partial \Om$ intersects $\partial B ( (- \e, 0); 1)$ at  some point   $p _ \e$ of the form 
 $$p_\e = f _ \e (\theta _\e) \, , \qquad \text { with } \theta _ \e \in \big ( 0, \frac \pi 2 \big )\,.$$ 
Moreover, due to the upper bound $\| r _\Om \| _\infty \leq 1$, it is easily seen that such a point is unique.  Similarly, there is a unique intersection point  between 
 $\partial \Om$ and  $\partial B ( (- \e, 0); 1)$ 
of the form 
$\tilde p_\e  = f _ \e (  \tilde \theta _\e) $, with $ \tilde \theta _\e \in \big ( - \frac \pi 2 , 0 \big ) $. 

As $\e \to 0 ^+$, the points $p_\e$ and $\tilde p _ \e$ move continuously along $\partial \Om$  towards $x_0$. 
Thus, taking into account that
 $\partial \Om$ is twice differentiable up to a $\mathcal H ^ 1$-negligible set, up to diminishing the value of $\e$, we can assume that $\partial \Om$ is twice differentiable at the points $p_\e$ and $\tilde p _ \e$. 
 
 This means that the angles $\Theta _\e$ and $\widetilde \Theta _ \e$  such that, using the parametrization \eqref{f:parangle}  of $\gamma _\Om$, 
 we have respectively $p _ \e = \gamma _\Om ( \Theta _ \e) $ and $ \tilde p _ \e = \gamma _\Om ( \widetilde \Theta  _ \e)$, are uniquely determined.  
 Notice also  for later use that $\Theta _\e$ satisfies by construction the following system
 \begin{equation}\label{f:systhetae} 
 \begin{cases}
 1 - \int _0 ^ { \Theta _\e}     r (t) \,  \sin t  \, d t  =  -\e + \cos \theta _\e  & 
 \\ 
 \noalign{\bigskip} 
    \int _0 ^ {\Theta _\e}     r (t) \,  \cos t \, d t =  \sin \theta _\e   \,.  & 
   \end{cases}
\end{equation}
 Let $B ( p _ \e; 1)$ and  $B ( \tilde p  _ \e; 1)$ be the balls of radius $1$ and centre $p _ \e$ and $\tilde p _\e $,  parametrized respectively by
 $$\begin{aligned}
 & g_\e (\theta) = (g_\e  ^1 (\theta) , g_\e ^2(\theta) ): = p _ \e  +  (\cos \theta , \sin \theta )\, , \qquad \theta  \in \big ( 0, 2 \pi)\,,
 \\ 
 & \tilde  g_\e (\theta)=  (\tilde g_\e  ^1 (\theta) , \tilde g_\e ^2(\theta) )  : = \tilde p _ \e  +  (\cos \theta , \sin \theta )\, , \qquad \theta  \in \big ( 0, 2 \pi)\,.
   \end{aligned}
$$ 
 
We consider the  perturbation $\Om _ \e$ of $\Om$ which is is constructed as follows:
 \begin{itemize}
 \item{} Replace  the portion of $\partial \Om$ corresponding to 
 $\gamma _\Om ( \theta)$ for $\theta \in [ \widetilde \Theta _ \e, \Theta _ \e]$ by  the arc of circle lying in $\partial B ( (- \e, 0); 1)$ 
 corresponding to $f _ \e (\theta)$ for $\theta \in [ \widetilde \theta _ \e, \theta _ \e ] $. 
 
 \smallskip 
 \item{}  Replace  the portion of $\partial \Om$ corresponding to 
 $\gamma _\Om ( \theta)$ for $\theta \in [  \pi , \Theta _ \e + \pi]$ by  the arc of circle lying in $\partial B ( p _ \e; 1)$ 
 corresponding to $g _ \e (\theta)$ for $\theta \in [\theta _ \e + \pi  ,  \Theta _ \e + \pi ] $.  

 \smallskip 
 \item{}  Replace  the portion of $\partial \Om$ corresponding to 
 $\gamma _\Om ( \theta)$ for $\theta \in [ \widetilde \Theta _ \e + \pi , \pi,  ]$ by  the arc of circle lying in $\partial B ( p _ \e; 1)$ 
 corresponding to $\tilde g _ \e (\theta)$ for $\theta \in [ \widetilde \Theta _ \e + \pi ,  \tilde \theta _ \e + \pi ] $.

 \end{itemize}
 
  By construction, we have that  $\Om _ \e$ belong to $\mathcal W_1$, with 
$\Om \setminus \Om _ \e$ contained into a small neighbourhood of $x_0$, and $\Om_ \e \setminus \Om$  contained into a small neighbourhood of $x_1$.
It remains to check that the equality \eqref{f:singlim} is satisfied.  We claim that such condition holds thanks to the following two separate asymptotics:
\begin{eqnarray}  
&\displaystyle \lim _ { \e \to 0 }  \frac {|  ( \Om _ \e \setminus \Om  )  \cap \{ (x, y) \ : \ y \leq 0   \} | }  {| (\Om  \setminus \Om_ \e ) \cap \{ (x, y) \ : \ y \geq 0   \}  |  }  = + \infty \,,
& \label{f:singlim1}  
\\ \noalign{\bigskip} 
& \displaystyle \lim _ { \e \to 0 }  \frac {|  ( \Om _ \e \setminus \Om  )  \cap \{ (x, y) \ : \ y \geq 0   \} | }  {| (\Om  \setminus \Om_ \e ) \cap \{ (x, y) \ : \ y \leq 0   \}  |  }  = + \infty \,.
& \label{f:singlim2}  
 \end{eqnarray} 
We focus attention on the proof of \eqref{f:singlim1}, the proof of \eqref{f:singlim2} being analogous. 

\smallskip 
We first prove the following claim: it holds
\begin{equation}\label{f:claimtheta}
\theta _ \e  ^ 2 = o ( \e ) \qquad \text{ in the limit as } \e \to 0^+\,. 
\end{equation} 
Indeed, we observe that the map $\e \mapsto \Theta _\e$ is monotone nonincreasing, and hence it admits a finite limit as $\e \to 0 ^ +$, that we denote by $\Theta _0$. We are going to prove the validity of \eqref{f:claimtheta} in the two distinct situations when $\Theta _0 >0$ and $\Theta _0 = 0$. 
Below, all the asymptotic equivalences and hierarchies (for which we adopt the usual notation $\sim$  and small $o$), are  tacitly meant in the limit as $\e \to 0 ^ +$. 

Assume first that $\Theta _0 >0$. Then,  for every $\theta \in [0, \Theta _0]$,  the direction $(\cos \theta, \sin \theta)$  lies in $\nu _\Om ( x_0)$, and hence  $r ( \theta) = 0$.  By using  \eqref{f:systhetae}, we see  that

$$
 \begin{cases}
 \sin \Theta _0 \int _{ \Theta _0}  ^ { \Theta _\e}     r (t)  \, d t  \sim 
 \int _{ \Theta _0}  ^ { \Theta _\e}     r (t) \,  \sin t  \, d t =  \int _0 ^ { \Theta _\e}     r (t) \,  \sin t  \, d t  =  1 +\e - \cos \theta _\e   = \e + \frac{\theta _\e ^ 2 }{2} + o ( \theta _ \e ^2 )& 
 \\ 
 \noalign{\bigskip} 
  \cos \Theta _0 \int _{ \Theta _0}  ^ { \Theta _\e}     r (t)  \, d t  \sim  
     \int _{ \Theta _0}  ^ {\Theta _\e}     r (t) \,  \cos t \, d t =   \int _{ 0}  ^ {\Theta _\e}     r (t) \,  \cos t \, d t =  \sin \theta _\e   = \theta _ \e + o ( \theta _ \e) \,,  & 
   \end{cases}
$$
which shows that
$$ \e + \frac{\theta _\e ^ 2 }{2}  \sim \theta _ \e \,,$$
yielding \eqref{f:claimtheta}. 

Assume now that $\Theta _0 =0$.  Then, from the assumption that $\xi_0$ is a  Lebesgue point  for $r_\Om$,  such that $r _\Om (\xi _0) = 0$,   we have
\begin{equation}\label{f:Lebesgue}
 \delta ( t)  := \frac{1}{t} \int _0 ^ { t } r ( \theta) \, d\theta  \,  \to \,   0  \qquad \text{ as } t \to 0 ^+ \,.
\end{equation} 
By \eqref{f:Lebesgue} and the assumption $\Theta _0 = 0$, we have
$$\delta  _ \e:= \sup _{ t \in [0, \Theta _ \e]} \delta ( t) \, \to \,  0 \qquad \text{ as } \e \to 0 ^ +\,. $$ 
Hence
$$\begin{aligned}
\int _0 ^ { \Theta _\e }  r ( t) \, t \, dt & =  \int _0 ^ { \Theta _\e }  \frac{r ( t) \,  }{\delta ( t ) } \int _0 ^ { t } r ( \theta) \, d\theta  \, dt  \geq 
\frac{1}{\delta_\e} \int _0 ^ { \Theta _\e }  {r ( t)  } \int _0 ^ { t } r ( \theta) \, d\theta  \, dt    \\ 
\noalign{\medskip}
& = \frac{1}{2\delta_\e} \int _0 ^ { \Theta _\e } \Big[\Big (  \int _0 ^ { t } r ( \theta) \, d\theta  \Big ) ^ 2 \Big ] '  \, dt 
=  \frac{1}{2\delta_\e} \Big (  \int _0 ^ { \Theta _\e} r ( \theta) \, d\theta  \Big ) ^ 2  \,. 
\end{aligned} 
$$ 
This shows that
\begin{equation}\label{f:teclemma}
\lim _ {\e \to 0 ^ + } \frac{  \Big ( \int _0 ^ { \Theta _\e } r ( t) \, dt \Big ) ^ 2  }{  \int _0 ^ { \Theta _\e }  r ( t) \, t \, dt }  =  0  \,.
\end{equation}
 Now, in order to obtain  \eqref{f:claimtheta}, it is enough to insert into  \eqref{f:teclemma} 
the following asymptotic expansions, which are deduced from 
 \eqref{f:systhetae}:
$$
 \begin{cases}
 \int _{ 0}  ^ { \Theta _\e}     r (t) \,  t  \, d t \sim  \int _0 ^ { \Theta _\e}     r (t) \,  \sin t  \, d t  =  1 +\e - \cos \theta _\e   = \e + \frac{\theta _\e ^ 2 }{2} + o ( \theta _ \e ^2 )& 
 \\ 
 \noalign{\bigskip} 
     \int _{0}  ^ {\Theta _\e}     r (t) \,   d t \sim   \int _{ 0}  ^ {\Theta _\e}     r (t) \,  \cos t \, d t =  \sin \theta _\e   = \theta _ \e + o ( \theta _ \e) \,.  & 
   \end{cases}
$$ 
Having thus achieved the proof of  claim \eqref{f:claimtheta}, we are ready to show \eqref{f:singlim1}. 
Let us show that the areas at the numerator and the denominator in \eqref{f:singlim1} can be estimated respectively as: 
\begin{eqnarray}
&  {| (\Om  \setminus \Om_ \e ) \cap \{ (x, y) \ : \ y \geq 0   \}  |  }    \leq \e \sin \theta _ \e + o ( \e \theta _ \e) & \label{f:show1}  
\\ \noalign{\medskip} 
& {|  ( \Om _ \e \setminus \Om  )  \cap \{ (x, y) \ : \ y \leq 0   \} | }  \geq \e ^ {\frac{3}{2}}  + o \big (  \e ^ {\frac{3}{2}} \big ) \,. & \label{f:show2}   
\end{eqnarray}
The proof will be achieved by combining these estimates  with claim \eqref{f:claimtheta}. 

Let us prove \eqref{f:show1}. Denoting for brevity  by 
$\Delta _ \e^ +$ the set  $(\Om  \setminus \Om_ \e ) \cap \{ (x, y) \ : \ y \geq 0   \} $, and by $\nu _ \e ^+$ its unit outer normal, we compute the area of $\Delta _ \e ^ + $ as the normal flow of the vector field $(x, 0)$ through the boundary: 
$$\begin{aligned} |\Delta _ \e ^ + | 
& 
= \int _{\partial \Delta _ \e ^ + } (x, 0) \cdot \nu _ \e ^ +  = \int _0  ^ { \Theta _\e  } \gamma _\Om ^ {1} (\theta)    \big ( \gamma _\Om ^ {2 } \big ) '  ( \theta  )  \, d \theta  -  \int _0  ^ { \theta _\e  } f_\e  ^ {1} (\theta)    \big (f_\e  ^ {2 } \big ) '  ( \theta  )  \, d \theta 
\\ 
&=  \int _0 ^ { \Theta _\e}  \Big ( 1 - \int _0 ^ { \theta }     r (t) \,  \sin t  \, d t  \Big ) r ( \theta) \cos (\theta) \, d \theta - \int _0 ^ {\theta _\e} 
 (-\e + \cos \theta  ) \cos \theta \, d \theta
\\ 
&=  \int _0 ^ { \Theta _\e}   r ( \theta) \cos (\theta) \, d \theta
- \int _0 ^ { \Theta _\e}  r ( \theta) \cos (\theta)  \int _0 ^ { \theta }     r (t) \,  \sin t  \, d t  \, d \theta 
 + \e \sin \theta _\e - \frac{\theta _\e}{2} - \frac{\sin ( 2 \theta _\e)} {4} 
\\  &=  \sin \theta _ \e - \int _0 ^ { \Theta _\e}  r ( \theta) \cos (\theta)  \int _0 ^ { \theta }     r (t) \,  \sin t  \, d t  \, d \theta 
 + \e \sin \theta _\e - \frac{\theta _\e}{2} - \frac{\sin ( 2 \theta _\e)} {4}  \,,
\end{aligned}
$$  
where in the last equality we have used the second equation in \eqref{f:systhetae}.  Hence
$$\begin{aligned} |\Delta _ \e ^ + |  & \leq \sin \theta _ \e 
 + \e \sin \theta _\e - \frac{\theta _\e}{2} - \frac{\sin ( 2 \theta _\e)} {4}   
 \\  & \leq \theta _ \e + \e \sin \theta _ \e  - \frac{\theta _\e}{2} - \frac{   \theta _\e} {2}     + \frac{  \theta _\e ^ 3 } {3}  + o (\theta _\e ^ 3)
 \\ 
 &
 = \e \sin \theta _\e + o (\e \theta _\e)\,,
 \end{aligned}
 $$ 
where in the last equality we have used  \eqref{f:claimtheta}. 

Let us prove \eqref{f:show2}.  We estimate from below  $|(\Om_\e   \setminus \Om ) \cap \{ (x, y) \ : \ y \leq 0   \} |$ 
by the area of the set $\Delta _ \e ^ -:= (\Om_\e   \setminus \Om ) \cap \{ (x, y) \ : \ y \leq 0, \ x\leq 0   \} $.  
The proof of \eqref{f:show2} will be achieved by showing that
\begin{equation}\label{f:partegiu}
|\Delta _ \e ^ - |  = \e ^ {\frac{3}{2}}  + o \big (  \e ^ {\frac{3}{2}} \big ) \,.
\end{equation}  
Denoting by 
$\nu _ \e ^-$ the unit outer normal to $\Delta _\e ^ -$, again we can compute the area of $\Delta _ \e ^ - $ as 
the normal flow of the vector field $(x, 0)$ through the boundary. Moreover,  the normal component of the vector field  $(x, 0)$ 
vanishes on $\partial \Delta _ \e ^ -$ except for its portion lying in $\partial B (p _ \e; 1)$. So we have
 $$ |\Delta _ \e ^ -  |
= \int _{\partial \Delta _ \e ^ - } (x, 0) \cdot \nu _ \e ^ -  =
\int _{\partial \Delta _ \e ^ -  \cap \partial B ( p _ \e, 1)} (x, 0) \cdot \nu _ \e ^ - \,,
$$ 
From  the construction of $\Om _\e$, we see that  $\partial \Delta _ \e ^ -  \cap \partial B ( p _ \e, 1)$ is parametrized by the function 
   $$ g_\e (\theta) = ( p _ \e  +  \cos \theta ,  p _ \e + \sin \theta )=  ( - \e + \cos \theta _\e  + \cos \theta\,,\  
   \sin \theta _\e + \sin \theta)
 \, , \qquad \theta \in [ \theta _ \e + \pi , \overline \theta _ \e + \pi] \,,$$ 
 where the limiting angle $\overline \theta _\e$ is determined by the condition
\begin{equation}\label{f:barthe}
  - \e + \cos \theta _\e  + \cos (\overline \theta  _\e  + \pi ) = 0\,.
  \end{equation} 
Equivalently,  in terms of the angle $\om:= \theta - ( \theta _\e + \pi )$, a parametrization of  $\partial \Delta _ \e ^ -  \cap \partial B ( p _ \e, 1)$  is given by 
   $$ \begin{aligned} g_\e (\om + \pi + \theta _\e)
 &   = ( - \e + \cos \theta _\e  + \cos (\om + \pi + \theta _\e)\,,\ 
   \sin \theta _\e + \sin (\om + \pi + \theta _\e))
   \, , \qquad \om \in [   0, \om _\e  ] 
   \\
   &   = ( - \e + \cos \theta _\e  - \cos (\om + \theta _\e)\,,\ 
   \sin \theta _\e - \sin (\om  + \theta _\e)) \,,
      \end{aligned} 
 $$  
where the limiting angle $\om _\e$  (equal to $\overline \theta _ \e - \theta _\e)$ is 
determined by the condition 
 \begin{equation}\label{f:barthe2}
  - \e + \cos \theta _\e  -  \cos (\om _ \e + \theta  _\e   ) = 0\,.
  \end{equation} 
Hence, if $\om _\e$ is defined by the above equation, we can  write $|\Delta _ \e ^ -  |$  as 
  $$\begin{aligned} |\Delta _ \e ^ -  | 
  & = \int _0 ^ {\om _\e}  ( - \e + \cos \theta _\e  - \cos (\om + \theta _\e) ) (- \cos  (\om  + \theta _\e))  ) \, d \om
  \\ \noalign{\medskip}
  & =   (  \e - \cos \theta _\e   ) \big ( \sin (\om  _\e  + \theta _\e)  - \sin  ( \theta _\e)  \big )  + \frac{\om _\e }{2} + \frac{1}{4} \big ( \sin ( 2  (\om _\e + \theta _\e)) 
   - \sin ( 2   \theta _\e)  \big )
   \\ \noalign{\medskip}
  & =   -  \cos (\om _ \e + \theta  _\e   )  \big ( \sin (\om  _\e  + \theta _\e)  - \sin  ( \theta _\e)  \big )  + \frac{\om _\e }{2} + \frac{1}{4} \big ( \sin ( 2  (\om _\e + \theta _\e)) 
   - \sin ( 2   \theta _\e)  \big )\,,
    \end{aligned} 
    $$ 
    where in the last equality we have used \eqref{f:barthe2}. Since by  \eqref{f:barthe2},  $\om _\e$ is infinitesimal as $\e \to 0$, we have 
      $$ \begin{aligned}  |\Delta _ \e ^ -  |  
       =  & - \Big  ( 1 - \frac{1}{2}   (\om _ \e + \theta  _\e   ) ^ 2 + o ( (\om _ \e + \theta  _\e   ) ^ 2 ) \Big )
       \Big ( \omega _\e - \frac{1}{6}  (\om _ \e + \theta  _\e   ) ^ 3 + o ( (\om _ \e + \theta  _\e   ) ^ 3 )
        + \frac{1}{6}   \theta  _\e    ^ 3 + o (  \theta  _\e    ^ 3 )
        \Big )   \\ \noalign{\medskip} 
        & + \frac{\om _\e }{2} +  \frac{\om _\e }{2} -  \frac{1}{3}  (  \om _\e + \theta _\e) ^ 3 + o ( (\om _ \e + \theta  _\e   ) ^ 3 )
 + \frac{1}{3}  \theta _\e ^ 3  + o (  \theta  _\e    ^ 3 )
        \end{aligned} 
        $$ 
    
 Thus, in order to estimate  $|\Delta _ \e ^ -  |$, we need to precise the asymptotic behaviour of $\delta _\e$ as $\e \to 0$. 
Let us show that 
 \begin{equation}\label{f:asomega}
 \om_\e ^ 2 \sim \e \,. 
 \end{equation} 
 Indeed, starting from \eqref{f:barthe2} we obtain 
 $$- \e + 1 - \frac{\theta_\e ^ 2}{2} + o (\theta_\e ^ 2)- 1 + \frac{\theta_\e ^ 2}{2} + \frac{\omega_\e ^ 2}{2}   + \theta _\e \omega _\e + o ( (\omega _ \e + \theta_\e  ) ^ 2) ) = 0 $$ 
 which using \eqref{f:claimtheta} entails 
 $$- \e   + o (\e) + \frac{\omega_\e ^ 2}{2}   + \theta _\e \omega _\e + o ( (\omega _ \e + \theta_\e  ) ^ 2)  = 0 \,.$$  
Recalling that $\theta _\e = \e ^ { \frac 1 2}+ o ( \e ^ { \frac 1 2} )$,   we get
 $$
  - \e   + o (\e) + \frac{\omega_\e ^ 2}{2}   +  \big (   \e ^ { \frac 1 2}+ o ( \e ^ { \frac 1 2} ) \big )  \omega _\e + 
  o \big ( (\omega _ \e +  (  \e ^ { \frac 1 2}+ o ( \e ^ { \frac 1 2} )   ) )^ 2 \big )   = 0$$ 
or equivalently  
\begin{equation}\label{f:oe}  \omega _\e ^ 2 \big (  \frac{1}{2}   +   o ( 1) \big  ) +     \e ^ { \frac 1 2} \omega _\e  \big (  1 +  o ( 1    ) \big  ) = 
\e \big (  1 + o (1)\big )  \,.
 \end{equation}   
 For $\e$ small enough we have 
 \begin{equation}\label{f:2lati}  \frac{\omega _\e ^ 2  }{4} 
 \leq \omega _\e ^ 2 \big (  \frac{1}{2}   +   o ( 1) \big  ) +     \e ^ { \frac 1 2} \omega _\e  \big (  1 +  o ( 1    ) \big  )  \leq {\omega _\e ^ 2  }    + \frac{\e}{4} + 4 \omega _\e ^ 2 = 5 \omega _ \e ^ 2 + \frac{\e}{4}  \,.
 \end{equation} 
By combining \eqref{f:oe} respectively with the first and with the second inequality in \eqref{f:2lati}, we obtain 
$$\e \big (  1 + o (1)\big )  \geq 
\frac{\omega _\e ^ 2  }{4}  \qquad \text{ and } \qquad \frac 3 4  \e \big (  1 + o (1)\big )  \leq {5 \omega _\e ^ 2  } \,,$$
thus achieving the proof of \eqref{f:asomega}. 

In view of \eqref{f:claimtheta} and \eqref{f:asomega},  we have $\theta_\e = o (\omega _\e)$, so that
the expression found above for  $|\Delta _ \e ^ -|$ has the following asymptotic expansion: 
   $$ |\Delta _ \e ^ -  |  
       =   - \Big  ( 1 - \frac{1}{2}   \om _ \e ^2   +  o ( \omega _ \e ^ 2  ) \Big )
       \Big ( \omega _\e - \frac{1}{6}  \om _ \e ^ 3 + o ( \om _ \e  ^ 3 )
              \Big )   + {\om _\e }  -  \frac{1}{3}   \om _\e  ^3  + o ( \om _ \e  ^ 3 ) =    \frac{1}{3}   \om _\e  ^3  + o ( \om _ \e  ^ 3 )\,.
        $$ 
In view of \eqref{f:asomega}, this shows the validity of \eqref{f:partegiu}  and achieves our proof.
 \qed

\bigskip\bigskip 
Our next step is to show that, still working under the assumptions of Lemma \ref{l:nor0}, 
the value of $r_\Om$ on each of the $4$ free arcs appearing in the decomposition \eqref{f:r*} is uniquely determined. 
To that aim, we 
are going to label the free arcs   so that, assuming up to exchange $E$ with its complement that $E$ is convex,
\begin{eqnarray} 
&\displaystyle \xi \in \gamma_1  \  \Rightarrow \ \ \nu _ \Om   ^ { -1} (\xi)  \in    \partial E   \quad \text{ and } \quad  
 \nu _ \Om  ^ { -1} (-  \xi ) \in  \partial (\Om \setminus E)  
& \label{f:gammasym}  
\\  \noalign{\medskip} 
&\displaystyle \xi \in \gamma_2  \    \Rightarrow \ \  
\nu _ \Om   ^ { -1} ( \pm \xi  ) \in   \partial (\Om \setminus E) \,.  
& \label{f:gammasymbis} 
  \end{eqnarray}

\begin{lemma}\label{l:blaschke} Under the same assumption of Lemma \ref{l:nor0}, we have: 
 \begin{eqnarray} 
 &  r_\Om (\xi) =    \frac{|E| ^ 2} {|E| ^ 2 + |\Om \setminus E| ^ 2 }  \qquad \hbox{ on }  \gamma _1  & \label{f:determin1} 
\\   \noalign{\medskip}
 & r _\Om (\xi) =     \frac{|\Om \setminus E| ^ 2} {|\Om \setminus E| ^ 2 + |E| ^ 2 }     \qquad \hbox{ on } \gamma_1 ^s     \,.& \label{f:determin2} 
\\   \noalign{\medskip}
&  r_\Om  (\xi) = \frac 1 2 \qquad \hbox{ on }  \gamma _2  \cup \gamma _2 ^ s \,.& \label{f:determin3}  
 \end{eqnarray}
  \end{lemma}

 The proof of Lemma \ref{l:blaschke} is based on a  shape derivative argument, and demands preliminarily to recall some basic facts about bodies of constant width, and
to introduce some definitions. 

\smallskip 
 Let $\Om \in \mathcal W _1$. From the equality  $h _ \Om ( \xi) + h _ \Om ( -\xi ) = 1$ for all $\xi \in \mathbb S ^ 1$, 
it follows that the same equality holds for  $r _\Om$ as well. Recalling also that 
the surface area measure of $\Om$ is equilibrated, 
we infer that $r_\Om$ belongs to the following class of functions: 
$$\mathcal A _1 := \Big \{ \varphi : \mathbb S ^ 1 \to [0,1]\,,\quad \varphi  ( \xi) + \varphi(  -\xi ) = 1 \,, \ \int_{\mathbb S ^ 1 } \xi \varphi ( \xi) =  0 \Big \}\,.$$ 
Conversely, by the Minkowski theorem, for any $\varphi \in \mathcal A _1$, there exists a unique convex body $\Om _\varphi$ such that $r _{\Om_\varphi} = \varphi$
(notice that the condition that $\varphi$ is not concentrated on any half-circle, required  to have existence in the Minkowski theorem, is 
automatically satisfied by all functions in $\mathcal A _ 1$). 

It will be convenient to associate, with any $\varphi \in \mathcal A _1$, the function
$$\overline \varphi : = \varphi - \frac{1}{2} \,,$$
and to introduce the class
$$\begin{aligned} 
\overline {\mathcal A _1} &  := \Big \{ \psi: \mathbb S ^ 1 \to [ - \frac 1 2 , \frac 1 2 ]\,,\quad \psi ( \xi) + \psi (  -\xi ) = 0 \,, \ \int_{\mathbb S ^ 1 } \xi \psi ( \xi) =  0 \Big \}
\\ 
&  = \Big \{ \overline \varphi  :  \varphi \in \mathcal A _1 \Big \} \,.
\end{aligned} 
$$

Notice that, for every $\varphi \in \mathcal A _1$ we have
$$\int _ {\mathbb S ^ 1 } \varphi = \pi \, , \qquad 
\int _ {\mathbb S ^ 1 }  \overline \varphi = 0\,.$$ 

Let  $\Om \in \mathcal W _1$. We call a {\it Minkowski perturbation of $\Om$ in $\mathcal W _1$} a family of convex bodies in $\mathcal W _1$
which is written, for some function $\varphi \in \mathcal A _1$, as 
$$\Om _ \e [\varphi]= ( 1- \e ) \Om + \e \Om _\varphi\,.$$ 

Since we are in dimension $2$, Blaschke and Minkowsi addition agree, namely we have 
$$r _ { \Om _ \e[\varphi]} = ( 1 - \e ) r _\Om + \e \varphi = r _\Om  + \e (\varphi - r _\Om ) \,.$$ 
Moreover, the corresponding first variation of the volume functional  can be representad as an integral over $\mathbb S ^ 1$, precisely as 
\begin{equation}\label{f:recall}  
 \frac{d}{d \varepsilon}  |  \Om _ \e[\varphi] |  _ { \e = 0 }   = \int_{\mathbb S ^ 1} h _\Om ( \varphi - r _\Om) \,. 
 \end{equation}

\medskip
\underbar {Proof of Lemma \ref{l:blaschke}}. 
Let  $( \Om, E)$ be an optimal pair as in the assumptions.

We are going to consider Minkowski perturbations  of $\Om$ of the form $\Om _ \e = \Om _ \e [\varphi]$ of the form 
\begin{equation}\label{f:perturb}
\varphi - r _\Om =:  \overline \psi  \in \overline {\mathcal A _ 1}\, , \qquad  \overline{\{\overline \psi \neq 0 \}} \subset   \
 ( \gamma _i\cup \gamma _i ^ s ) \cap \big \{\xi \ :\  r _\Om (\xi) \in  ( \delta, 1 - \delta)  \big \}\,,
\end{equation}  
where, for $i= 1$ or $i = 2$,  $\gamma_i$ is a free arc in the decomposition \eqref{f:r*}, and $\delta$ is a fixed number in $(0, 1)$.

We shall later make the distinction between the cases $i =1$ and $i =2$. 

Notice in particular that, for $\e$ small enough, since $\overline \psi$ is supported on a set where
$r _\Om$ is detached from $0$ and $1$,  we have $\Om _ \e \in \mathcal W _1$. 

Any such perturbation $\Om _ \e [\varphi]$ induces in a natural way a perturbation   $ E _\e$  of $ E$, where 
$E _\e\subset \Om _\e$  is obtained by from $E$  by  replacing  the portion of  $\partial E$ lying in $\partial \Om$ by
its perturbation lying in $\partial \Om _\e[\varphi]$.

\smallskip
Notice in particular that the perturbation $\Om _\e$ of $\Om$ keeps fixed the endpoints of 
$\nu _\Om ^ { -1} (\gamma_i)$  and of $\nu _\Om ^ { -1} (\gamma_i^s)$, which are secured points of $\partial \Om$: 
this is due to the fact that 
both $\overline \psi \res \gamma$ and $\overline \psi \res \gamma^s$ have zero barycenter. Indeed, 
\begin{equation}\label{f:barparziali} 0 = \int _{\gamma_i}  \xi  \overline \psi (\xi)  +  \int _{\gamma_i ^s} \xi  \overline \psi (\xi) = 
2 \int _{\gamma _i} \xi \overline \psi (\xi) =  2 \int _{\gamma _i ^ s}  \xi \overline \psi (\xi)  \,,
\end{equation} 
where the first equality holds because $\overline \psi$ is supported into $\gamma _i \cup \gamma _i ^ s$ and has zero barycenter, and the second one because $\overline \psi$ is odd.  

\smallskip

Then, from the optimality of $\Om$ in the minimization problem \eqref{f:mincw}  and the optimality of $E$ for $\sigma _ 1 (\Om)$, 
the corresponding first order shape derivative of 
 $$F (\Om, E):= { \frac{1}{|E| } + \frac{1}{|\Om \setminus E| } }\, $$
 must be equal to zero. Such shape derivative is given by 
\begin{equation}\label{f:shapeF}
\begin{aligned} dF (\Om, E)& :=  \frac{d}{d \varepsilon}  F (\Om_\e,  E_\e)\Big | _ { \e = 0 } 
\\ 
&  =   - \frac{1}{|E|^2 } \frac{d}{d \varepsilon}  |  E_\e|  _ { \e = 0 }  
 - \frac{1}{|\Om \setminus E| ^2}   _ { \e = 0 }\frac{d}{d \varepsilon}  |  \Om _ \e \setminus E_\e|  \,.
  \end{aligned}
  \end{equation}

Let us compute $\frac{d}{d \varepsilon}  |  E_\e|  _ { \e = 0 } $ and  $\frac{d}{d \varepsilon}  |  \Om _ \e \setminus E_\e|$. 

To that aim, we  distinguish the cases $i = 1$ and $i = 2$. 
Assume first that $i = 1$. Let us  
consider 
 the convex bodies 
 $\Om _ \e ' $ and $\Om _ \e ''$ which  are defined respectively by
\begin{equation}\label{f:mink} 
r _{ \Om' _ \e} = r _\Om + \e \overline \psi \res \gamma _i \qquad \text{ and } \qquad r _{ \Om'' _ \e} = r _\Om + \e \overline \psi \res \gamma_i  ^s \,.
\end{equation} 
Notice that,  though they do not belong to $\mathcal W_1$,  for $\e$ small enough $\Om' _ \e$  and $\Om '' _\e$ are well-defined by Minkowski theorem, because  both
 $r _\Om + \e \overline \psi \res \gamma_1 $ and $ r _\Om + \e \overline \psi \res \gamma_1  ^s$ 
are not supported on a half circle (thanks to Lemma \ref{l:nor0}), and have zero barycenter (cf.\ \eqref{f:barparziali}). 

From \eqref{f:mink}, by using the known formula for the shape derivative of volume under Blaschke addition, we infer that 
\begin{equation}\label{f:shapeK}   \begin{cases}
\displaystyle \frac{d}{d \varepsilon}   |\Om' _ \e | _ { \e = 0 } =  \frac{d}{d \varepsilon}   |\Om' _ \e\cap S | _ { \e = 0 }  =   \int_{\mathbb S ^ 1} h _\Om ( \overline \psi \res \gamma_1  ) = \int _{\gamma_1 }  h _ \Om \overline \psi
& \\ \noalign{\medskip} 
\displaystyle  \frac{d}{d \varepsilon}   |\Om '' _ \e | _ { \e = 0 }
=   \frac{d}{d \varepsilon}   |\Om '' _ \e \cap S | _ { \e = 0 } =   \int_{\mathbb S ^ 1} h _\Om ( \overline \psi \res \gamma_1 ^s ) = \int _{\gamma_1  ^ s} h _ \Om \overline \psi
= -  \int _{\gamma_1}  (1 -h _\Om )  \overline \psi  \,,  & 
\end{cases} 
\end{equation}      
where the last equality holds because
$$ \int _{\gamma _1 ^ s}  h _\Om  (y) \overline \psi (y)  =  \int _{\gamma _1  ^ s}  [1 -h _\Om (-y)]  [- \overline \psi  (-y) ] = -  \int _{\gamma _1  }  [1 -h _\Om (z)]  \overline \psi  (z) \,,$$ 
Next we observe that, by \eqref{f:gammasym},  we have 
 \begin{equation}\label{f:scompo}
|E _\e | = | \Om' _ \e \cap S | - |\Om \setminus E| \qquad \text{ and }  \qquad |(\Om_\e \setminus E_\e) \cap S  | =  | \Om'' _ \e  | - | E|  \,. 
\end{equation} 
By \eqref{f:shapeK} and \eqref{f:scompo}, we infer  
\begin{equation}\label{f:claime}  \frac{d}{d \varepsilon}  |  E_\e|  _ { \e = 0 }   =  \int_{\gamma_1} h _\Om \overline \psi \qquad \text{ and } 
\qquad \frac{d}{d \varepsilon}  |  \Om _ \e \setminus E_\e| _ { \e = 0 }   =  \int_{\gamma_1^s} h _\Om \overline \psi  
\end{equation}
 
Inserting \eqref{f:claime} into \eqref{f:shapeF}, we get 
\begin{equation}\label{f:formdF}
\begin{aligned}  dF (\Om, E) 
& =  \int _ {\gamma _1}    \Big ( - \frac{h_\Om}{|E|^2 }  + \frac{ 1 - h _\Om }{|\Om \setminus E| ^2} \Big )  \,\overline \psi \\ 
& =  \frac{1}{|\Om \setminus E| ^2 } \int _ {\gamma _1}  \Big [ 1  - h _ \Om \Big ( \frac{ |E| ^ 2 + |\Om \setminus E| ^2 }{|E|^2}  \Big )    \Big  ] \, \overline \psi  
\\ 
& =  \frac{1}{| E| ^2 } \int _ {\gamma _1   ^s} \Big [ 1  - h _ \Om \Big ( \frac{ |E| ^ 2 + |\Om \setminus E| ^2 }{|\Om \setminus E|^2}  \Big )    \Big  ] \, \overline \psi  
\,.  \\ 
\end{aligned}
\end{equation} 
Writing for brevity $C (\Om, E):= \frac{|E| ^ 2} {|E| ^ 2 + |\Om \setminus E| ^ 2 }$, 
we have found so far that the equalities 
\begin{equation}\label{f:condimento}  0 =  \int_{\gamma_1} [ h _\Om - C ( \Om, E)]  \overline \psi  =     
\int_{\gamma_1 ^ s} [ h _\Om -  (1-C ( \Om, E))]  \overline \psi     
 \end{equation}
hold for every $\overline \psi \in \overline {\mathcal A _ 1}$ such that, for some $\delta \in (0, 1)$,   $\overline{\{\overline \psi \neq 0 \}} \subset   \
 ( \gamma _1\cup \gamma _1 ^ s ) \cap \big \{\xi \ :\  r _\Om (\xi) \in  ( \delta, 1 - \delta)   \big \}$.  
 Let us show how this implies  the equalities  \eqref{f:determin1}-\eqref{f:determin2}. 
Since the proof of the two equalities is analogous, 
let us focus on the proof \eqref{f:determin1}. 

We claim that 
\begin{equation}\label{f:gollazo}
\int _{\gamma_1} [ h _\Om - C ( \Om, E)]  {f} = 0  \qquad \forall f \in L ^ 2 ( \gamma_1) \ :\ \int _{\gamma _1} \xi  f (\xi) = 0 \,.
\end{equation} 
Indeed, setting
$$S _n := \Big \{ \xi \in \mathbb S ^ 1 \ :\ r _\Om (\xi) \in  \Big ( \frac 1 n, 1 - \frac 1 n \Big ) \Big \} \,,$$ 
 by Lemma \ref{l:nor0} we have that
 $$\mathcal H ^ 1 \Big (  \gamma _1  \setminus \bigcup _{n \in \N } S _n \Big  )=  0\,.$$ 
As a consequence, the sequence of functions 
$$f _n:= \big (   f   \vee  (-n) \wedge n \big  ) \chi _ {S _n }    \,, $$  
converges to $f$ in $L ^ 2 (\gamma_1)$.   In order to gain also the zero barycenter condition, we fix  an integer $\overline n >2$ such  that 
$\gamma _ 1 \cap S _{\overline n } \neq \emptyset$,
and we
select subsets $A, B \subset \gamma_1 \cap S _{\overline n }$  such that the matrix 
$$M:= \left ( \begin{matrix} 
\int_{A} \xi _1 & \int  _{B } \xi _ 1 
\\ \noalign{\medskip} 
\int_{A} \xi _2 & \int  _{B } \xi _ 2
\end{matrix}\right  ) 
$$  
  is invertible.  Then we modify the sequence $f_n$ into the new sequence
  $$\widetilde f _n = f _n + a_n \chi _{A} + b _ n \chi _B \,, $$ 
where $(a_n, b _n)$ is the unique solution to the following system: 
$$\begin{cases}
\int_{\gamma _1} \xi _ 1 f _n + a _n \int _  A \xi _1 + b _n \int _  B \xi _1 = 0 & 
\\    \noalign{\medskip} 
\int_{\gamma _1 } \xi _ 2 f _n + a _n \int _  A \xi _2 + b _n \int _  B \xi _2 = 0  \,.
 \end{cases}
$$ 
Notice that $a_n$ and $b _n$ are infinitesimal because $\int _{\gamma_1} \xi f_n \to 0$ (as a consequence of 
the assumption $\int_{\gamma _1} \xi f = 0$ and of the fact that
 $f_n\to f$ in $L ^ 2 (\gamma_1)$).  In particular, this ensures that $\widetilde f _n \to f$ in $L ^ 2 (\gamma_1)$. 
  
Now, for every $n\geq \overline n$,  we have that the odd extension to $\mathbb S ^ 1$ of  the function $\frac {\widetilde f _n} { 2 \| \widetilde f_n \| _ \infty }$ belongs to the class $\overline {\mathcal A _ 1}$, and it is supported in $ ( \gamma _ 1 \cup \gamma _ 1 ^ s )  \cap S _n$. 
It can therefore be chosen as a test function $\overline \psi$ in \eqref{f:condimento}, yielding
$$ 0 =  \int_{\gamma_1} [ h _\Om - C ( \Om, E)]  \widetilde f _ n = 0\,.$$ 
The claimed equality \eqref{f:gollazo} follows by passing to the limit as $n \to + \infty$,

\smallskip 
Now we observe that, since the vector equation  $ \int _{\gamma _1} \xi  f (\xi) = 0$ can be rephrased as the system $ \int _{\gamma _1} \xi _i f (\xi) = 0$ for $i = 1, 2$, 
 condition \eqref{f:gollazo} implies  that the function $h _\Om - C ( \Om, E)$ belongs to the subspace $M$ of $L ^ 2 (\gamma_1)$ spanned by the functions $\xi _1$ and $\xi _2$. 
By writing $\gamma_1$ under the form
 $$\gamma_1 = \big \{ \xi (\theta) = (\cos \theta, \sin \theta) \ : \ \theta \in (\alpha , \beta)  \big \}\,, $$
we infer  that there exists constants $A, B \in \R$ such that 
$$ h _\Om (\xi (\theta))  - C ( \Om, E)  = A \cos \theta + B \sin \theta \qquad \forall \theta \in (\alpha , \beta) \,.$$ 
Recalling that $r _\Om= h'' _\Om + h _\Om$,  we conclude that 
$$r_\Om ( \xi ( \theta) ) = C ( \Om, E) \qquad \forall \theta \in (\alpha , \beta) \,. $$

Assume now that $i =2$. In this case by \eqref{f:gammasymbis}  we have  $E _\e = E$, so that  $ \frac{d}{d \varepsilon}  |  E_\e|  _ { \e = 0 }$.
On the other hand, 
by \eqref{f:recall},  for all functions $\overline \psi \in \overline {\mathcal A_1}$ such that  
$\overline{\{\overline \psi \neq 0 \}} \subset \gamma_2  \cup \gamma_2 ^ s$,  we obtain   
\begin{equation}\label{f:claimebis}  \frac{d}{d \varepsilon}  | \Om _\e \setminus E_\e|  _ { \e = 0 }   =  \frac{d}{d \varepsilon}  | \Om _\e|  _ { \e = 0 }   = \int_{\gamma_2} h _\Om \overline \psi +  \int_{\gamma_2 ^s} h _\Om \overline \psi   =   \int_{\gamma_2} h _\Om \overline \psi  -  \int _{\gamma_2}  (1 -h _\Om )  \overline \psi  =  \int_{\gamma_2}  \big ( 2h _ \Om - 1 \big ) \overline \psi \,. 
\end{equation} 
Starting from this equality, and arguing as done above for  $i =1$, we obtain $r _\Om ( \xi) = \frac{1}{2}$ on $\gamma  _2 \cup \gamma _2  ^ s$.  

\qed

 \bigskip
 We are now ready to derive the following

 \bigskip 
\underbar{Conclusion of the proof of Theorem \ref{t:constwidth}}.   
 Assume that problem \eqref{f:mincw}  is not solved by $B _ \frac{1}{2}$. 
  Let $\Om$ be a solution to the minimization problem \eqref{f:mincw} as given by Lemma \ref{l:symmetry}. By Lemma  \ref{l:blaschke}, if we fix a coordinate system with the $\vec x$ axis along the bisector of $S$, there exists an angle $\theta \in  (0, \pi)$ such that, writing $\xi = (\cos t, \sin t)$ for $t \in (0, 2 \pi)$, and letting $r (t): = r _\Om ( \cos t, \sin t)$, we have 
     \begin{eqnarray} 
 &  r (t) =  \overline r :=   \frac{|E| ^ 2} {|E| ^ 2 + |\Om \setminus E| ^ 2 }  \qquad \hbox{ for }   t\in  (0, \theta) \cup  ( \pi , \pi + \theta)   & \label{f:determin10} 
\\   \noalign{\medskip}
 & r (t)   = 1 - \overline r  \frac{|\Om \setminus E| ^ 2} {|\Om \setminus E| ^ 2 + |E| ^ 2 }     \qquad \hbox{ for } t\in  (\pi - \theta, \pi) \cup ( 2\pi- \theta  ,  2 \pi )        \,.& \label{f:determin20} 
\\   \noalign{\medskip}
&  r (t)   =  \frac 1 2 \qquad \hbox{ on }  \qquad \hbox{ for } t\in  (\theta, \pi - \theta) \cup ( \pi + \theta , 2 \pi - \theta)   \,.& \label{f:determin30}  
 \end{eqnarray}
 Written in terms of the function $r (t)$, the condition $\int _{\mathbb S ^ 1 } \xi  r _\Om (\xi) = 0$ reads: 
 $$\int_0  ^ {2 \pi} r ( t) \cos t \,dt = \int _0 ^ {2 \pi} r ( t) \sin t  \, dt= 0 \,.$$  
 Inserting into this system the expression of $r ( t)$ on $(0, 2 \pi)$ given by \eqref{f:determin10}-\eqref{f:determin20}-\eqref{f:determin30}, this is equivalent to
 $$( 2 \overline r - 1) \sin \theta =  ( 2 \overline r - 1)  (1- \cos \theta)  = 0\,.$$
 Since $\theta \neq 0$, we conclude that $\overline r = \frac{1}{2}$, and hence $r ( t) = \frac{1}{2}$ on $(0 , 2 \pi)$, which implies $\Om = B _ {\frac 1 2 }$, contradiction. \qed
 
 \bigskip

 \section {Proof of Theorem \ref{t:newexplicit}}\label{sec:teo3}
We have to prove that, for every open convex set $\Om \subset \R ^2$ with $D _\Om = 1$, it holds
$$\frac{\sigma_ 1 (\Om) -2} {\rho _\Om ^ 2} \geq \frac 4 3 \,.$$ 
To that aim it is enough to show that, 
 for every  fixed $r \in (0, \frac 12 )$,   it holds
$$\frac{\sigma_ 1 (\Om) -2} {r ^2} \geq \frac 4 3 \qquad \forall \Om \in \mathcal K ( r) \,, $$ 
  where
  $$
 \mathcal K ( r):= \Big \{ \Om  \text{ open bounded convex sets in $\R ^2$ with $D_\Om = 1$, $\rho _\Om \geq  r$} \Big \}\,, 
  $$
 Thus, for a fixed $r \in ( 1, \frac 1 2)$, we focus our attention on 
  the minimization problem
\begin{equation}\label{f:defsigma}
\Sigma ( r):=  \min _{\Om \in \mathcal K ( r) }   \sigma_ 1 (\Om)  \,,
\end{equation}
with the aim of showing that
\begin{equation}\label{f:aimshow}
\frac{\Sigma ( r) -2} {r ^2} \geq \frac 4 3  \,. 
\end{equation}

 \smallskip 
Let us introduce the subclass of $\mathcal K ( r)$ given by 
$$  \mathcal X ( r):= \Big \{ \Om \in \mathcal K ( r)\, , \ \Om = K  \cap S, \text{ for some } K \in \mathcal W_1, \ S \in \mathcal S  \text{ such that  }  \gamma _{1, 2}  \cap \partial K \neq \emptyset \Big \}\,, $$ 
   where $S$ denotes the family of strips bounded either by 
two parallel straight lines $\gamma _{1, 2}$, or by two incident half-lines (still denoted  $\gamma _{1, 2}$), with the incidence point outside $K$. 
In the sequel,  when dealing with elements $\Om = K  \cap S$ of $\mathcal X ( r)$, 
we call them:
\begin{itemize}
\item[--] {\it proper}  if $\Om \neq K$ ;

\smallskip  
\item[--]   {\it improper}  if $\Om = K$.
\end{itemize}

As a preliminary result, in the next lemma we show  that the minimum $\Sigma ( r)$ defined in \eqref{f:defsigma} is attained on a set $\Om$  belonging to  $\mathcal X ( r)$, and in addition we derive some geometric information on such a minimizer $\Om$ in case it  is proper. 
Then the proof of Theorem \ref{t:newexplicit} will be given by analyzing separately the two cases $\Om$ improper  
(in which the conclusion will be reached relying on Theorem \ref{t:constwidth}) and  $\Om$ proper 
(in which we shall exploit the geometric information coming from Lemma \ref{l:cw}).

 \begin{lemma}\label{l:cw}   For every fixed  $r \in (0, \frac{1}{2})$, 
 $\Sigma ( r)$ admits a minimizer   belonging to $\mathcal X ( r)$.  
 Moreover, either it admits an improper minimizer, or it admits a proper minimizer $\Om= K \cap S$ with the following properties: 
 $K$  is symmetric about the bisector of $S$, 
  $\sigma _ 1 (\Om )$ is attained at a set $E$ such that 
 the endpoints of the arc $\partial E \cap \Om $ belong to $\partial S\setminus \partial K$, and the  
 four  points in
 $\partial K \cap  \partial S$  are two by two diametral in $K$.  

  \end{lemma} 
 
\proof Let $r \in (0, \frac{1}{2})$ be fixed, and let us show that
 $\Sigma ( r)$ admits a minimizer $\Om =  K \cap S  \in \mathcal X ( r)$.  Let  $\Om _0$ be an optimal set for $\Sigma ( r)$, and let $E  _0$ be an optimal set for $\sigma _ 1 (\Om  _0)$. 
Let $S$ be the strip determined by the two tangent lines $\gamma _{ 1, 2} $ to $\Om  _0$ at the two points 
of intersection between $\partial \Om  _0$ and the arc 
$\gamma:= \partial E  _0 \cap \Om  _0$. 

If ${\Om  _0}$ has constant width,   we are done simply by taking 
$K = {\Om  _0}$.

Otherwise, we consider  a convex body  $V \in \mathcal W _1$ such that $V \supset {\Om  _0}$  and we claim that we are done by taking $K = V$ and 
$\Om  = V \cap S$. 
Indeed, we notice that $\gamma $ splits $\Om $ into two sets $E $ and $\Om \setminus E$ such that
 $E  \supseteq E  _0$ and $(\Om  \setminus E ) \supseteq  
 (\Om  _0 \setminus E  _0)$. 
 Therefore, 
 $$\sigma _ 1 (\Om ) \leq  \frac{1}{2} {{\rm Per} (E, \Om)} \Big ( { \frac{1}{|E| } + \frac{1}{|\Om \setminus E | } } \Big )  \leq 
  \frac{1}{2} {{\rm Per} (E_0, \Om_0 )}  \Big ( { \frac{1}{|E_0| } + \frac{1}{|\Om \setminus E_0 | } }   \Big ) = 
 \sigma _ 1 (\Om _0)\,.$$  
 Hence $\Sigma ( r) = \sigma _ 1 (\Om )$, and $E$ is optimal for $\sigma _1 (\Om)$, with the endpoints of $\gamma = \partial E \cap \Om$ belonging to $\partial S$ by construction.

\medskip
 Assume now that  $\Sigma ( r)$ does not admit improper minimizers. In particular, the minimizer $\Om= K \cap S$  found above is proper.   
To obtain the symmetry of $K$ about the bisector of $S$, we proceed in a  similar   way as in the proof of Lemma  \ref{l:symmetry} (i), except that now we have to pay attention to the constraint on the inradius.
Let $K _S$ be  the Steiner symmetrization  of $K$ with respect to the bisector of $S$.
Let $K_*$ be a set of of constant width equal to $1$ which contains $K_S$ and is symmetric about the bisector of  
$S$.  We set 
  $$\Om  _S :=   K_S \cap S  \,,  \qquad   \Om  _* :=   K_* \cap S\,,$$
 and we claim that
  \begin{equation}\label{f:doppiaug}
  \Sigma ( r ) = \sigma _ 1 (\Om  _S) = \sigma _ 1 (\Om  _*)\,.
  \end{equation}  

\smallskip
To prove the first equality in \eqref{f:doppiaug}, we notice
the arc $\gamma$ splits the set $\Om  _S$ into  a set  $E'$ 
and its complement $\Om  _S \setminus  E'$, with 
 $| E '| = |E|$, and $|\Om   _S \setminus E' | = |\Om  \setminus E|$, so that 
$$\sigma _ 1 (\Om  _S) \leq \frac{1}{2} {{\rm Per} ( E', \Om _S)}  \Big ( { \frac{1}{| E ' | } + \frac{1}{|\Om _S  \setminus  E ' | } }  \Big )
= \frac{1}{2} {{\rm Per} ( E, \Om )} \Big ( { \frac{1}{| E| } + \frac{1}{|\Om  \setminus E | } }   \Big ) = \sigma _ 1 (\Om) = \Sigma ( r) \,. $$

By the properties of Steiner symmetrization, we have
\begin{equation}\label{f:inradius}
\rho _{\Om  _S } \geq \rho _{\Om} \,, \qquad D _{\Om  _S} \leq D _{\Om}   ( = 1)\,;
\end{equation}
up to replacing $\Om _S$ by its dilation by a factor larger than $1$, we may assume that $D _{\Om  _S }= 1$. 
Hence $\Om  _S$ is admissible for $\Sigma (\rho _{\Om} )$, so that
$$\Sigma (\rho _{\Om})  \leq \sigma _ 1 (\Om  _S)\,.$$ 
Since, by the increasing monotonicity of  the map $r \mapsto \Sigma ( r)$, it holds 
\begin{equation}\label{f:stesso} 
\Sigma ( r) \leq \Sigma (\rho _{\Om} )\leq  \sigma _ 1 (\Om ) = \Sigma (r)\,,
\end{equation} 
we conclude that 
$\Sigma (r)  = \sigma _ 1 (\Om  _S)$, and $E'$ is optimal for   $\sigma _ 1 (\Om  _S)$. 

To prove the optimality of $ \Om  _*$, we proceed similarly as done above to prove the optimality of $\Om$, starting from the optimal set 
$\Om  _S$  in place of the initial optimal set $\Om _0 $.
  Notice that the endpoints of the optimal arc $\gamma$ belong by construction to $\partial S$. In addition, they  cannot belong to $\partial K$, since otherwise
by removing the strip $S$ we would obtain an improper minimizer, against the assumption.

\smallskip

\medskip
Finally, let us consider a proper minimizer $\Om_*$ for $\Sigma ( r)$ as found above, 
namely $\Om _* = K _* \cap S$, with $K_* \in \mathcal W _1$, $S \in \mathcal S$, and $K_*$ symmetric  about the bisector of $S$. 
Moreover we know that the endpoints of the 
arc $\gamma = \partial E\cap \Om _*$, being $E$ an optimal set for $\sigma _ 1 (\Om_*)$,   belong to $\partial S$. 
Since $\Om_*$ is proper and $K_*$ is symmetric about the bisector of $S$, 
$\partial K _*\cap \partial S$  consists of $4$ points. 
Let us show that they are two by two diametral in $K_*$.   
Let $P$, $Q$ be two among such points,  lying on distinct straight lines $\gamma _ { 1, 2}$  of $\partial S$, and on opposite sides with respect to the arc $\gamma$. We label them so that $P \in \partial E \cap \gamma _1$ and $Q \in \partial (\Om_* \setminus E) \cap \gamma _2$.  Assume by contradiction that $P$ and $Q$ are not diametral points in $K_*$.

Let $P', Q'$ denote respectively the diametral points of $P, Q$.  Since the two diameters $[PP']$ and $[QQ']$ must intersect each other, necessarily one among $Q'$ and $P'$ lies inside the strip $S$, and the other one outside. Assume with no loss of generality that $Q' \in S$, so that $Q' \in \partial E\cap \partial \Om_*$.  Then, by a continuity argument, we can find an arc $\sigma$  of points in $\partial E\cap \partial \Om _*$ whose diametral points lie outside $S$. 
We consider the perturbation $\Om _ \e$ of $\Om_*$ defined by
$$\Om _\e =  {\rm co} ( ( I + \e V ) (\Om))  \,,$$ 
where $V ( x): = \phi ( x) \nu _K ( x)$ is a normal deformation, with $\phi \in {\mathcal C} ^ \infty _0 (\sigma; (0, 1))$, and ${\rm co} ( \cdot)$ denotes the convex envelope.

We observe that, provided $\e$ and ${\rm spt} \, \phi$ are sufficiently small,   the set $(\Om _\e \setminus \Om_*)$ is contained into an arbitrarily small tubular neighbourhood of the arc $\sigma$. 
(This follows from  a direct geometric argument
relying on the  fact that,  since $r _{K_*} \in (0, 1)$, the curvature of $\partial K_*$ is bounded below by $1$).
 
Such property, combined with the fact every point in the arc $\sigma$ has its diametral point outside $S$, ensures that, 
 provided $\e$ and ${\rm spt} \, \phi$ are sufficiently small,  the diameter of $\Om _\e$ is still equal to $1$. 
 
On the other hand, since $\Om _\e \supseteq \Om_*$, the inradius  of $\Om _\e$ is at least $r$. 
  Moreover,  the arc $\gamma$ splits $\Om _\e$ into a set $E _\e$ and its complement $\Om _\e \setminus E _\e$ such that, by construction 
$|E _ \e| > |E|$ and  $|\Om _\e \setminus E _\e|  =  |\Om_*  \setminus E|$.     
This contradicts the optimality of $\Om_*$. 
 \qed

\bigskip \bigskip 
\underbar{Conclusion of the proof of Theorem \ref{t:newexplicit}}.
For a fixed $r \in (0, \frac{1}{2})$,  let  $\Sigma ( r) $ 
be defined in  \eqref{f:defsigma},  and let us prove \eqref{f:aimshow}. 
We consider  an optimal set for $\Sigma ( r)$  of the form $\Om= K \cap S \in \mathcal X ( r)$ as given  by Lemma \ref{l:cw}, and we distinguish two cases:

\smallskip
{\it Case 1}: $\Om$ is improper, i.e. $\Om = K \in \mathcal W _1$. Then we are done by applying Theorem \ref{t:constwidth}. Indeed: 
$$\frac{\Sigma ( r) -2} {r ^2} =  \frac{\sigma_ 1 (\Om) -2} {r ^2} \geq  \frac{\sigma_ 1 ( B _ {\frac 1 2 } ) -2} {r ^2} \geq  \frac{\sigma_ 1 ( B _ {\frac 1 2 } ) -2} {  \frac{1}{4}  }   =  4 \big ( {\frac{8}{\pi} -2}   \big )  >      \frac 4 3\,.$$

\medskip

\smallskip
{\it Case 2}: $\Om$ is proper, i.e. $\Om = K \cap S$, with $\Om \neq K$. Let us prove that $\Om$ has the following property:
\begin{equation}\label{f:parallel} 
\text{ the straight lines bounding $S$ are parallel.}
\end{equation} 
Let us assume by contradiction that this is not the case. Then we claim that  it is possible to construct 
a perturbation $(\Om _\e, E _\e)$ of $(\Om, E)$ such that
\begin{equation}\label{f:willgive} 
{\rm Per} ( E_\e, \Om _ \e) - {\rm Per} ( E, \Om  ) = o ( \e) \, \qquad  0 < |E _\e| - |E| \sim \e  \, , \qquad 
 0 < |\Om _\e \setminus E _ \e| - |\Om \setminus E| \sim   \e  \,. \end{equation}  
This will give 
$$
\frac{d}{d\e}  \big ( {\rm Per} ( E_\e, \Om _ \e)  F (\Om _\e, E _\e) \big)  \Big | _{ \e = 0} <0\,,  
$$ 
which, for $\e$ small enough impliese 
$\sigma _1 (\Om_\e) \leq {\rm Per} ( E_\e, \Om _ \e)  F (\Om _\e, E _\e) <  {\rm Per} ( E, \Om )  F (\Om , E ) = \sigma _ 1  (\Om)$, leading to the required contradiction.

Let us show how the perturbation can be constructed.

\bigskip\bigskip\bigskip\bigskip
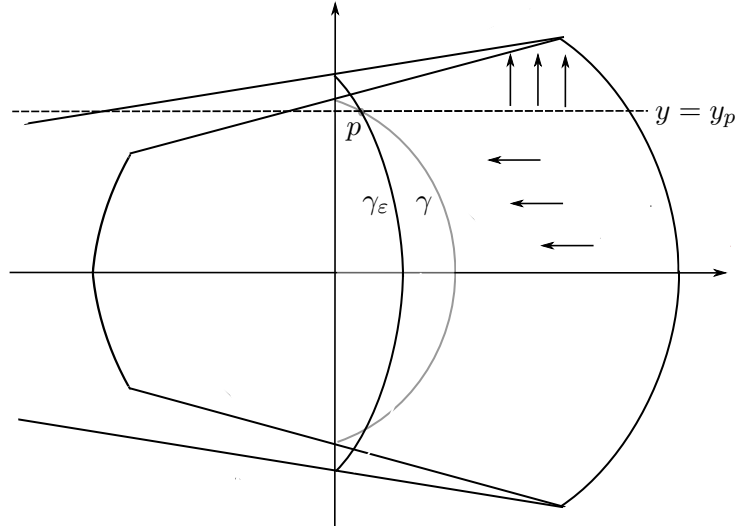
\begin{figure} [h] 
\centering   
\def\svgwidth{7cm}   
\begingroup%
  \makeatletter%
  \providecommand\color[2][]{%
    \errmessage{(Inkscape) Color is used for the text in Inkscape, but the package 'color.sty' is not loaded}%
    \renewcommand\color[2][]{}%
  }%
  \providecommand\transparent[1]{%
    \errmessage{(Inkscape) Transparency is used (non-zero) for the text in Inkscape, but the package 'transparent.sty' is not loaded}%
    \renewcommand\transparent[1]{}%
  }%
  \providecommand\rotatebox[2]{#2}%
  \ifx\svgwidth\undefined%
    \setlength{\unitlength}{265.53bp}%
    \ifx\svgscale\undefined%
      \relax%
    \else%
      \setlength{\unitlength}{\unitlength * \real{\svgscale}}%
    \fi%
  \else%
    \setlength{\unitlength}{\svgwidth}%
  \fi%
  \global\let\svgwidth\undefined%
  \global\let\svgscale\undefined%
  \makeatother%
  \begin{picture}(4,.8)%
    \put(0.6, 0){
     \includegraphics[height=7cm]{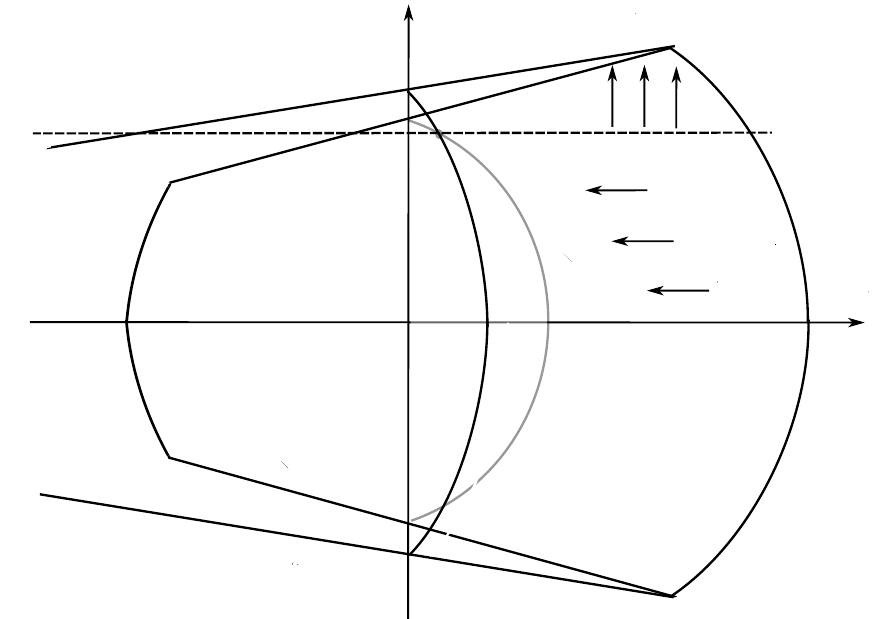}} %
         \put(1.3, 0.74){\color[rgb]{0,0,0}\makebox(0,0)[lb]{\smash{$p$}}}%
         \put(1.88, 0.77){\color[rgb]{0,0,0}\makebox(0,0)[lb]{\smash{$y = y _p$}}}%
         \put(1.33, 0.6){\color[rgb]{0,0,0}\makebox(0,0)[lb]{\smash{$\gamma _\e$}}}%
         \put(1.43, 0.6){\color[rgb]{0,0,0}\makebox(0,0)[lb]{\smash{$\gamma $}}}%

 \end{picture}%
\endgroup%
\caption{The geometry of Case 2}
\label{fig:2}   
\end{figure}

Let us fix for convenience a system of coordinates such that the vertical axis $\vec y$ is directed along  
the straight line passing through the endpoints of $\gamma:= \partial E \cap \Om$, 
the horizontal axis  $\vec x$  is directed along the bisector of $S$, and  
$(\partial E \cap \Om)  \subset \{   (x, y) \, :\,  x \geq 0 \}$. 

We consider an arc $\gamma _\e$ which is a small perturbation of $\gamma $, symmetric about the bisector of $S$, 
with the following properties: 
\begin{itemize} 
\item[--] $\gamma _\e$ has endpoints  on  $\vec y$, each one at distance $\e$ from an endpoint of  $\gamma$, 
in such way that   the segment joining the endpoints of $\gamma$ is contained into the segment joining 
the endpoints of $\gamma _\e$; 

\smallskip 
\item[--] if $\Gamma _\e$ is the subset delimited by $\vec y$ and $\gamma _\e$, it holds 
\begin{equation}\label{f:zerder1} \frac{d}{d\e} \big  | \Gamma _\e |    \Big | _{ \e = 0} =0
\end{equation} 
and  
\begin{equation}\label{f:zerder2} 
 \frac{d}{d\e}  \mathcal H ^ 1 (\gamma _\e) _{ \e = 0} =0\,.
 \end{equation}  
  \end{itemize} 
Let us describe how such arc $\gamma _\e$ can be constructed. By symmetry, we can limit ourselves to 
show that \eqref{f:zerder1}-\eqref{f:zerder2} are fulfilled replacing $\gamma _\e$ and $\Gamma _\e$ respectively by 
 $\gamma _\e ^+ := \gamma _\e \cap \{ (x, y)  :  y \geq 0 \}$ and $\Gamma _\e ^ +:= \Gamma _\e \cap \{ (x, y)  :  y \geq 0 \} $. 
 We fix a point $p = (x_p, y _p) 
\in \gamma$, with $x_p>0$, $y _p >0$, and we define 
  $$\gamma _\e  ^+    :=  ( Id + \e V)   (\gamma) \,,$$ 
  where $V$ is the vector field   
  \begin{equation}\label{f:defV} V (x, y):= \begin{cases} 
  (y - y _p) (0, 1) & \text{ for  } y \geq y _p
  \\ 
  c  x (y_p - y) (-1, 0) & \text{ for  } y < y _p\,,
  \end{cases} 
  \end{equation} 
 where $c$ is a constant chosen so that, denoting by $\nu$ the unit outer normal to the region  $\Gamma ^+$   delimited by $\vec y$, $\vec x$, and 
 $\gamma  ^+:= \gamma \cap \{ (x, y)  \,: \, y \geq 0 \} $, it holds  
 \begin{equation}\label{f:richiesta} 
 \int _{\gamma  ^+} V \cdot \nu = 0\,. 
 \end{equation}  
Notice that this is possible because  by construction it amounts to solve the following linear equation in the unknown $c$: 
$$\begin{aligned} 0 = \int _{\gamma  ^+} V \cdot \nu 
& = \int _{\gamma  ^+\cap \{ y \geq y _p \}} V \cdot \nu  + \int _{\gamma  ^+\cap \{ y < y _p \}} V \cdot \nu
\\ 
& = \int _{\gamma  ^+\cap \{ y \geq y _p \}}   (y - y _p) (0, 1)   \cdot \nu  + c \int _{\gamma  ^+\cap \{ y < y _p \}}   x (y_p - y) (-1, 0) \,. 
\end{aligned} 
$$ 
Then the requirement \eqref{f:zerder1} is fulfilled since
$$ \frac{d}{d\e} \big  | \Gamma _\e ^+ |    \Big | _{ \e = 0}  = \int _{\partial \Gamma ^+ } V \cdot \nu  = 
\int _{\partial \Gamma ^+ \cap \vec y } V \cdot  (-1, 0)  + \int _{\partial \Gamma ^+ \cap \vec x } V \cdot  (0, -1)   + \int _{\gamma ^+ } V \cdot \nu \,,  
$$ 
and the first two integrals at the right hand side vanish by \eqref{f:defV}, while the third one vanishes by the choice of $c$ ensuring \eqref{f:richiesta}. 

As a consequence of \eqref{f:richiesta}, and of the fact that  $\gamma $ is an arc of circle, we have that also the requirement \eqref{f:zerder2} is fulfilled: indeed,  denoting by $r$ the radius of curvature of $\gamma$, we have 
$$
 \frac{d}{d\e}  \mathcal H ^ 1 (\gamma _\e ^+) _{ \e = 0} =\int _{ \gamma ^+ } \frac{1}{r}  \,  V \cdot \nu  = 0  \,.
$$  

We are now in a position to construct a perturbation $(\Om _\e, E _\e)$ of $(\Om, E)$ such that
conditions \eqref{f:willgive}  hold. 

We replace $\Om = K \cap S$ by $\Om _\e = K \cap S _\e$, where  the new strip $S _\e$ 
is a small perturbation of $S$ such that the two straight lines bounding $S_\e$ pass each one by 
 an endpoint of $\gamma _\e$ and by an  intersection point  between $\partial (\Om \setminus E)$ and $\partial S$. 

We replace $E$ by a small perturbation of $E$ obtained as one of the two subsets 
of $\Om _\e$ in which the arc $\gamma _\e$ splits $\Om _\e$.

 Then it is readily seen that all the conditions in \eqref{f:willgive} are satisfied. Indeed, 
the condition ${\rm Per} ( E_\e, \Om _ \e) - {\rm Per} ( E, \Om  ) = o ( \e)$ is fulfilled thanks to \eqref{f:zerder2}.  
 The second and third conditions   in \eqref{f:willgive} are  satisfied thanks to \eqref{f:zerder1}, since  we have   
 $$   \begin{aligned}  |E _\e|   - |E| & =   
 \big |(\Om _ \e \setminus \Om  )  \cap \{ (x, y) \, : \,  x \leq 0 \} \big  |   +  |\Gamma _\e|   -| \Gamma |   
 \\ 
& = 
  \mathcal H ^ 1 (\partial S \cap \partial \Om\cap \{ (x, y) \, : \,  x \leq 0 \} )  \, \e  + 
 o ( \e) \,,
  \end{aligned} $$ 
 and in a similar way
 $$   \begin{aligned}  |\Om _ \e \setminus E _\e|   - |\Om \setminus E| & =  
  \big |(\Om _ \e \setminus \Om  )  \cap \{ (x, y) \, : \,  x \geq 0 \} \big  |   - | \Gamma _\e|  +  |\Gamma|   
 \\ 
& = 
  \mathcal H ^ 1 (\partial S \cap \partial \Om\cap \{ (x, y) \, : \,  x \geq 0 \} )  \, \e  + 
 o ( \e) \,.
  \end{aligned} $$ 
 
 This concludes the proof claim \eqref{f:parallel}.  
 In view of the claim just proved, we are reduced to show that, if  $\Om\in \mathcal K ( r)$ is the intersection of a ball of radius $\frac 1 2$ with a strip delimited by two parallel straight lines symmetric with respect to the  centre of the ball, it holds
  $$\frac{ \sigma _ 1 (\Om) -  2}{r ^ 2}  \geq \frac{4}{3}\,.$$
  Clearly, since $\rho _\Om \geq r$ for all $\Om \in \mathcal K ( r)$, it is enough to show that 
    $$\frac{ \sigma _ 1 (\Om) -  2}{\rho _\Om  ^ 2}  \geq \frac{4}{3}\,.$$ 
Moreover, we are going to show that the lower bound is asymptotically attained in the limit when $\rho _\Om$ tends to $0$.
   
 We observe that, if $E$ is  an optimal set  for $\sigma _1 (\Om)$, we have that
 $\partial E \cap \Om$ is a line segment of length $2 \rho _\Om$ orthogonal to the boundary of the strip,  and $|E| = |\Om \setminus E|$. Thus, writing 
for  brevity $\rho$ in place of $\rho _\Om$, the above inequality is equivalent to 
 \begin{equation}\label{f:equivto}
 \frac{1}{\rho ^ 2 } \Big [  \frac{1}{2} {{\rm Per} (E, \Om)} {\Big(  \frac{1}{|E| } + \frac{1}{|\Om \setminus E| } \Big ) }  - 2  \Big ]  = 
\frac { 2 (  \rho - |E|)} {\rho ^ 2 |E| } 
 \geq \frac{4}{3}\,. 
 \end{equation}
In terms of the angle
$$\theta := \arcsin ( 2 \rho) \in \big ( 0  , \frac \pi 2 \big  ) \,,$$   
we have
$$|E| =   \frac{1}{2}  \cos \theta \sin \theta + \frac{\theta}{4} - \frac{1 }{8}   \sin ( 2 \theta)  \,.  $$ 
Hence 
$$
\frac { 2 (  \rho - |E|)} {\rho ^ 2 |E| } 
  =  8 \, \frac{ 2 \sin \theta - \sin \theta \cos \theta - {\theta}  }{ \sin ^ 2 \theta \big ( \sin \theta \cos \theta + \theta  \big )   }  \,,
$$
 so that the inequality \eqref{f:equivto} can be rephrased as 
 \begin{equation}\label{f:entheta}
 \mathcal E (\theta): =  8 \, \frac{ 2 \sin \theta - \sin \theta \cos \theta - {\theta}  }{ \sin ^ 2 \theta \big ( \sin \theta \cos \theta + \theta  \big )   } \geq \frac{4}{3}\,,
 \qquad \forall \theta \in  \big ( 0  , \frac \pi 2 \big  )\,,
\end{equation} 
 or
 $$
f (\theta) := 12 \sin \theta - 6 \sin \theta \cos \theta - 6 \theta - \sin ^ 3 \theta \cos \theta - \theta \sin ^ 2 \theta \geq 0\,,  \qquad \forall \theta \in  \big ( 0  , \frac \pi 2 \big  )\,.
$$
Since $f ( 0) = 0$, in order to prove that $f (\theta) \geq 0$, it is enough to show that $f'(\theta ) \geq 0$.  We have
$$f' (\theta) = 12 \cos \theta - 12 \cos ^ 2 \theta - 4 \sin ^ 2 \theta \cos ^ 2 \theta - 2 \theta \sin \theta \cos \theta \,,$$ 
so that the inequality $f' ( \theta ) \geq 0$ is in turn equivalent to 
$$
g (\theta )
:=  6 \frac{1 -  \cos  \theta  }{\sin \theta} - 2  \sin \theta \cos \theta - \theta 
  \geq 0\,, \qquad  \forall \theta \in  \big ( 0  , \frac \pi 2 \big  )\,.
  $$ 
Such inequality is satisfied because $\lim _{\theta \to 0 ^ +} g (\theta ) = 0$, and  
$$g' ( \theta) = \frac{3 ( 1 - \cos \theta ) ^ 2 + 4 \sin ^ 4 \theta}{\sin ^ 2 \theta}  \geq 0\,.$$ 
Finally, we observe that the lower bound  \eqref{f:equivto} is asymptotically attained in the limit as $\rho \to 0 ^+$. Indeed, 
as $\theta \to 0 ^+$,
the elementary asymptotic expansions 
$$2 \sin \theta - \sin \theta \cos \theta - {\theta}  = \frac{\theta ^ 3 }{3} + o ( \theta  ^ 3) \qquad \text{ and } 
\qquad \sin ^ 2 \theta \big ( \sin \theta \cos \theta + \theta  \big )  = 2 \theta ^ 3 + o ( \theta ^ 3)\,,$$
show that the map $\theta \mapsto \mathcal E (\theta )$ in \eqref{f:entheta} satisfies
$$\lim _ { \theta \to 0 ^+} \mathcal E (\theta ) = \frac{4}{3} \,.$$ 
\qed

\section {Appendix}\label{sec:appendix}

\subsection{Asymptotics of the constant $\ka_0$ as $p \to 1^+$} 

Below we give a proof  of the claim  made in the Introduction that 
the constant    $\ka _0$ appearing in \eqref{f:ABF2}, as given by Theorem 1 in \cite{ABF2} for $N = 2$, tends to $0$ as 
$p \to 1 ^+$.

\medskip Adopting the same  notation as in \cite{ABF2}, we have that 
$$\ka _0=  
 \frac{ \ka_1 (\ka_2) ^ 2 }{6(7 \cdot 16 \cdot 256  )^2}\,,$$  being $\ka _1$ and $\ka _2$ given respectively   in eqs. (39) and (64) in \cite{ABF2}. 
We have
$$\ka_2 =      \frac{ 1}{2^{N+m}}   \frac{1}{  \ka_\infty  ^ p   \cdot \ka _{kroger}  ^ {N+m} }  \,,$$ 
where $\ka _\infty$ is independent of $p$ (see \cite[eq. (17)]{ABF2}) and  $  \ka _{kroger}   = (p+1)2^{\frac{N+  3 p}{2}}     N ^ { N+p} $.   Hence, in the limit as $p \to 1 ^+$, $\ka _2$ tends to a positive and finite limit, so that  our claim  $\ka _0 \to 0$ follows if we show that  $\ka _1 \to 0$. We have (see \cite[Section 4.1]{ABF2})
$$\ka _1 = 
 \frac{1}{2 ^ {\frac{1}{p-1}}}
  \frac{e ^  { -  (m+1)  2  ^ {m+3}  \gamma ^ { p(m+2)}    }   }{  \gamma  ^  {\frac{2p^2-p}{p-1} }} 
M  \Big (   p,  \frac{1}{4 \gamma ^ p } ,  \frac{1}{(4 \gamma ) ^ {\frac{1}{p-1}} }\,  ,
 ( 4 \gamma  ^ p) ^ { m+1}   \big ( (\pi _p ) ^ p + 1 \big ) \Big )  \,, $$ 
where
$$\gamma:=  \ka _\infty(\pi _ p ^ p + 1 )^{\frac{ m +1 }{p}} \frac{1}{{\left (\int_{I} 
( \min \{x , 1 - x \} ^ {m+1} ) \right )}^{\frac 1 p} }\,,  $$ 
$$M ( p, a, s, t) :=   \Big ( \frac{p-1}{p}  \Big ) ^ { p-1} \cdot \min \Big \{ \frac{b _0 ^2}{2}, \frac{b_0}{2}    s ^{  \frac{p}{p-1} } \Big \}\, , \qquad b _0 ( p, a, s , t) := \min \Big \{ \frac{a}{2},  \Big ( \frac{1}{2} \frac {p-1}{p} \frac{s}  { t}  \Big )  ^ {\frac{p}{p-1}  }  \Big \}\,.$$
Then the fact that  $\ka _1\to 0$ follows from by observing that $$ \frac{1}{2 ^ {\frac{1}{p-1}}}
  \frac{e ^  { -  (m+1)  2  ^ {m+3}  \gamma ^ { p(m+2)}    }   }{  \gamma  ^  {\frac{2p^2-p}{p-1} }}  \to 0  \qquad \text{ and } \qquad 
  M  \Big (   p,  \frac{1}{4 \gamma ^ p } ,  \frac{1}{(4 \gamma ) ^ {\frac{1}{p-1}} }\,  ,
 ( 4 \gamma  ^ p) ^ { m+1}   \big ( (\pi _p ) ^ p + 1 \big ) \Big )  \to 0 \,.$$ 
Indeed, the first convergence above holds because  $\gamma $ converges to a positive and finite limit,  while the second one holds because  
$$M  \Big (   p,  \frac{1}{4 \gamma ^ p } ,  \frac{1}{(4 \gamma ) ^ {\frac{1}{p-1}} }\,  ,
 ( 4 \gamma  ^ p) ^ { m+1}   \big ( (\pi _p ) ^ p + 1 \big ) \Big )  \sim\min \Big \{ \frac{b _0 ^2}{2}, \frac{b_0}{2} 
 \frac{1}{(4 \gamma ) ^ {\frac{p}{(p-1) ^ 2}} }  \Big \}  \leq \frac{b_0}{2} 
{(4 \gamma ) ^ {- \frac{p}{(p-1) ^ 2}} }  \,,$$ 
and the last term is infinitesimal because ${(4 \gamma ) ^ {- \frac{p}{(p-1) ^ 2}} } \to 0$, while 
 $$b _0 \Big (   p,  \frac{1}{4 \gamma ^ p } ,  \frac{1}{(4 \gamma ) ^ {\frac{1}{p-1}} }\,  ,
 ( 4 \gamma  ^ p) ^ { m+1}   \big ( (\pi _p ) ^ p + 1 \big ) \Big )    \leq \frac{1}{4 \gamma ^ p } \,. $$ 

\bigskip
\subsection{Proofs of the results of Section \ref{sec:fence}}
In order to prove Propositions \ref{p:k} and \ref{p:ad},  below we provide  extensions 
to the weighted  space $BV _\p (\Om)$ of 
several facts which are well-known in the classical space $BV (\Om)$
(see  e.g.\ \cite{AFP}). 
 
 \begin{proposition}\label{p:BV}  
 Let $\Om$ be  an open bounded set    in $\R^N$ and let   $\p$ be a positive function in $L ^ 1 (\Om)$.  
 
\smallskip 
(i) For every $u \in BV _\p  (\Om)$, if  $\{ u _ h \} \subset \mathcal D ( \R ^N)$ is any sequence of functions such that $u _ h \to u $ in $L ^ 1  _\p (\Om)$, it holds
$ \liminf _h \int _\Om |\nabla u _h |  \phi \geq    |D _\p u| (\Om) $; moreover, there exists a sequence $u _h  \subset \mathcal D ( \R ^N)$ such that $u _ h \to u $ in $L ^ 1  _\p (\Om)$  and $\lim _h \int _\Om |\nabla u _h |  \phi =  |D _\p u| (\Om) $. 
 
 \smallskip
 (ii) Statement (i) remains true replacing $\mathcal D ( \R ^N) $ by $W ^ { 1, 1} _\p (\Om)$. 
 
 \smallskip
 (iii)  For every  $u \in W^{1, 1} _ \p (\Om)$,   we have $u \in BV _\phi (\Om)$ and 
$|D_\p u | (\Om) = \int _\Om | \nabla u| \p $.

\smallskip
 
 (iv) For every $u \in BV _\p  (\Om)$, if  $u ^ + := \max \{ u, 0 \}$ and $u ^ -:= \min \{ u, 0\}$, it holds
$|D _\phi u |  (\Om) = |D _\phi u ^+|  (\Om)  + |D _\phi u ^-|  (\Om)$.

 \smallskip
 (v) For every $u \in BV _\p  (\Om)$, the level sets $\{ u >t \}$ have finite $\p$-perimeter in $\Om$ for a.e.\ $t \in \R$, and it holds 
 $$|D_\p u| (\Om) = \int _ \R  \Per_ \p \big ( \{ u > t \}, \Om \big ) \, dt\,.$$
 
 \smallskip
 (vi)  Assume in addition that $\p$ is 
 $( \frac{1}{m})$-concave, for some $m \in \N \setminus \{ 0 \}$. 
 Then, 
  if $\{ u _h \}$ is a bounded sequence in $W ^ { 1, 1} _\phi (\Om)$, up to subsequences we have
$u _h \to u$  in $L ^ 1 _\p (\Om)$, for some function $u$ in  $BV  _\p (\Omega)$.

 \end{proposition} 
 
 \proof 
  (i)  See \cite[Theorem 5.1]{BBF99}. 

 \smallskip 
 (ii)  See \cite[Remark 5.2]{BBF99}. 

\smallskip 
 (iii) See  \cite[Theorem 5.7]{BBF99}.  
 
 \smallskip

\smallskip
(iv)  Let $u \in BV _\p (\Om)$. If  $\{ u _ h \} \subset  W ^ { 1, 1} _\p (\Om)$ is any sequence  such that $u _ h \to u $ in $L ^ 1  _\p (\Om)$, it holds 
$u _ h  ^ \pm \to u  ^ \pm $ in $L ^ 1  _\p (\Om)$. Then, since
$$ \liminf _h \int _\Om |\nabla u _h | \p \geq  \liminf _h  \int _\Om |\nabla u _h ^+  | \p  + \liminf _h   \int _\Om |\nabla u _h ^- | \p   \,$$
we have  $$|D_\phi u |  (\Om)  = \inf \Big \{\liminf \int _\Om |\nabla u _h | \p \ :\ u _ h \to u  \text{  in } L ^ 1  _\p (\Om) \Big \}  \geq |D_\p u ^ +  | (\Om) +  |D_\p u ^ -  | (\Om)\,.$$ 
On the other hand, by statement (ii), there exist sequences  $u _h ^ \pm  \subset W ^ { 1, 1} _\p (\Om)$  such that $u _ h ^ \pm \to u^ \pm  $ in $L ^ 1  _\p (\Om)$  and $\lim _h \int _\Om |\nabla u _h ^\pm |  \phi =  |D _\p u ^\pm | (\Om) $. Since the sequence $u _h:= u _h ^ + +   u _h ^-$ is in $W ^ { 1, 1} _\p (\Om) $ and converges to $ u $ in $L ^ 1  _\p (\Om)$, we have 
$$ |D_\p u | (\Om) \leq \lim_h  \int _\Om |\nabla u _h   |  \leq \lim _h \Big (   \int _\Om |\nabla u _h ^+  | \p +   \int _\Om |\nabla u _h ^-  | \Big ) =  |D _\p u ^ +  | (\Om) + |D _\p u ^ -  | (\Om) \,.$$ 

 \smallskip
 (v) See \cite[Theorem 4.1]{BBF99}.  
 
 \smallskip
 (vi) 	Let  $\widetilde{\Omega} \sq \R^{N+m}$  be the open convex set defined by
	 $$
	 \widetilde{\Omega } = \left\{ (x,y) \in \R^{N}\times \R^{m} : \, x \in \Omega ,\,   \norma{y}_{\R^m} <     \omega _m ^ {- \frac{ 1}{m}}   \p^{1/m}(x)\right\}\,,
	 $$
	  where $\om _m$ is the Lebesgue measure of the unit ball in $\R ^m$,    
	 and consider the sequence of functions 
	 $$
	\tilde{u}_h(x,y):= u_h(x) \qquad \forall (x, y ) \in \widetilde \Om \,.$$
 We have 
	 \begin{gather*}
	\int_{\tilde{\Omega}} | \tilde {u}_h| (x,y) \, dxdy = \int_{\Omega}  |u_h (x)| \p(x) \, dx ,\\
		\int_{\tilde{\Omega}} \abs{\nabla_{x,y}\tilde{u}_h} (x,y) \, dxdy = \int_{\Omega} \abs{\nabla_x u_h} (x)  \p(x) \, dx.
	 \end{gather*}
 	Thus the sequence $\{\tilde u _h\}$ is bounded in $BV (\widetilde \Om)$. From the compact embedding 
 $
 	BV (\widetilde\Omega) \hookrightarrow L^1(\widetilde{\Omega})
 $
we infer that, up to a subsequence, there exists a function $\tilde u \in BV (\widetilde \Om) $ such that
$\tilde u _h  \to \tilde u$  in  $L ^ 1 (\widetilde \Om)$.
Since  the function $\tilde u$ is independent of the variable $y$, setting $u(x): = \tilde{u}(x,y)$ we get 	  
 $$\int _\Om  |u _ h - u | \p = \int _{\widetilde \Om} | \tilde u _h - \tilde u| \to 0 \,,$$
 namely $u _h\to u \in L ^ 1 _\phi (\Om)$. Since the sequence $\{u _h\}$ is bounded in 
 $W ^ { 1, 1} _\phi (\Om)$,  by statement  (ii) we infer that
 $|D_\phi u | (\Om) \leq \liminf \int _\Om |\nabla u  _h |\p < + \infty$, and hence $u \in BV _\p (\Om )$.  
  \qed

 \bigskip
 \bigskip 
\underbar {Proof of Proposition \ref{p:k}}. 
Throughout the proof, we write for brevity $\m$, $\overline\m$, and $\m ^*$ in place of $\m   (\Om,\p)$, $\overline \m  (\Om,\p)$,  and $\m ^*  (\Om,\p)$.  
 Let us show first that $\overline \m =\m  ^*$.   It is immediate that $\overline \m  \leq \m  ^*$, since  the  minimization problem which defines $\m  ^*$ is obtained by  restricting the class of admissible functions to the family of characteristic functions.  The converse inequality follows by using the coarea formula
stated in Proposition \ref{p:BV} (v): for any  $u \in BV  _\p (\Om; \R ^+)$ with $|\{ u >0 \} | _\p  \leq \frac{|\Om | _\p} {2}$, we have
$$ \frac {|D _\p u| (\Om) } {\int_ \Om |u|  \p }  = \frac{ \int _ 0 ^ {+ \infty} \Per  _\p ( \{ u > t \}, \Om ) \, dt }{\int _ 0 ^ {+ \infty} |  \{ u > t \} | _\p  \, dt }  \,. $$ 
Since, by definition of $\m ^*$, we have that  the inequality $\Per  _\p ( \{ u > t \} , \Om ) \geq  \m  ^* |  \{ u  > t \} |  _\p $ holds for every $t \in (0, + \infty)$, by the arbitrariness of $u$ it follows that 
 $\overline \m  \geq \m  ^*$. Moreover  we see that, if  $u$  is optimal for $\overline \m $, almost all its level sets
 $\{ u > t\} $ are optimal for $\m  ^*$.

We now prove that $ \m \geq \overline \m$. Let $\{ u _h \} \subset W ^ { 1, 1}_\phi (\Om)$ be a sequence with  $\int_\Om (\sign  u _h )\phi = 0$ and  $\frac {\int _\Om |\nabla u_h|  \phi} {\int_ \Om |u_h| \phi } \to \m $. 
  By Proposition \ref{p:BV} (vi), up to subsequences we have  $u _h \to u$  in $L ^ 1 _\p (\Om)$. 
  Then $u _h ^ \pm \to u ^ \pm $  in $L ^ 1 _\p (\Om)$, and applying  
 Proposition \ref{p:BV}  (ii) and (iv), we have
 $$ \m  = \lim _h   \frac {\int _\Om |\nabla u_h|  \p } {\int_ \Om |u_h| \p  } \geq  \frac { |D _\p u|  (\Om)  } {\int_ \Om |u| \p }    =  
  \frac { |D _\p u^+| (\Om)  + \int _\Om |D _\p u^-|  (\Om) } {\int_ \Om |u^+| \p  + \int_ \Om |u^-| \p } \,.$$

 From the condition  $\int_\Om ( \sign  u _h )\p = 0$ and the convergence of $u _ h$ to $u$ in $L ^ 1 _\p (\Om)$, it follows that
 $ | \{ u >0  \} |   _\p \leq  \frac{|\Om| _\p }{2}$ and  $| \{ u <0  \} |  _\p  \leq  \frac{|\Om| _\p}{2}$.  Moreover, we observe that 
 $$  \frac {|Du^+|   _\p  (\Om)  } {\int_ \Om |u^+| \p  } =   \frac {  |Du^-| _\p  (\Om) } {\int_ \Om |u^-|  \p } \,.$$  
 Indeed, if by contradiction one of the two terms would be smaller than the other, e.g. the left hand side,  then   for some $\alpha \in (0, 1)$, letting
 $v:= u ^+ - \alpha u ^ -$,  and  $v_h:= u_h ^+ - \alpha u_h ^ -$, it would be 
$$ \lim _h   \frac {\int _\Om |\nabla v_h|\p } {\int_ \Om |v_h| \p} =  \frac { |D_\p v| (\Om) } {\int_ \Om |v|\p  }  < \frac { |D_\p u| (\Om) } {\int_ \Om |u|\p  } \leq \m  \,,$$ 
against the definition of $\m $. 

We conclude that 
$$ \m  \geq   \frac { |D _\p u^+| (\Om)    } {\int_ \Om |u^+| \p   }  \geq  \min _ {\substack{u \in BV _\p (\Om; \R^+)\\  | \{ u >0  \} |  _\p  \leq  \frac{|\Om|_\p }{2} }}  \frac { |D_\p u|  (\Om) } {\int_ \Om |u| \p }  = \overline \m \,,$$

It remains to show that 
$\overline \m \geq \m$, or equivalently that $\m ^*\geq \m$. (Notice that the  inequality $\overline \m \geq \m$, combined with the above proof of the converse inequality, shows also that $\overline \m$ is attained).  
In order to prove that $\m ^*\geq \m$,  it is enough to show that, 
 for any set  $E \subset \Om$, with $|E|_\p \leq \frac{|\Om | _\p }{2}$, there exists a sequence $\{ u _h \} \subset W ^ { 1, 1} _\p (\Om)$, with 
$ \int_\Om  (\sign  u_h )  \p = 0$  such that
\begin{equation}\label{f:gollo} 
\frac{  \int _\Om |\nabla u _h |_ \p  }{\int _\Om {|u _h | _\p}}  \to \frac{  |D _\p \chi _E| (\Om)}{|E|_\p}  \,.
\end{equation} 

Let us first consider the case when $|E|_\p < \frac{|\Om | _\p }{2}$. In this case we can find a set $F \subset \Om$, with
$\Per_\p(F, \Om)  < + \infty$, such that $F \cap E = \emptyset$ and $|E|_\p  = |F|_\p$ (for instance by cutting  $\Om$ with a moving hyperplane).
Notice that by construction the set $H:= \Om \setminus ( E \cup F)$ turns out to be of strictly positive measure.  

For any $\delta  >0$,  we define the function $g _\delta := - \delta \chi _ F + \chi _E$. We have $ \int_\Om  (\sign  g_\d)  \p = 0$  and, in the limit as $\d \to 0$, 
 $$
g _\d \to \chi _ E  \text{  in } L ^ 1  _\p (\Om ) \qquad \text{ and } \qquad   |D _\p g _\d | (\Om)   \to \Per_\p(E, \Om)\,.$$ 
 
Therefore, by a diagonal argument, it is enough to show that \eqref{f:gollo} holds with $\chi _E$ replaced by $g _\d$, for a fixed $\delta >0$.
Since  $g _\delta$ belongs to $BV  _\p (\Om)$, there exists a sequence $\{ v _h \} \subset W ^ {1, 1} _\p (\Om)$ such that
$v _h \to g _\d$ in $L ^ 1 _\p (\Om)$, and $ \int _\Om |\nabla v _h | \p  \to |D_\p g _\d | (\Om)$.   Moreover, by a standard approximation argument, we can assume without loss of generality that every function $v_h$ has no plateau in its graph, namely that the level sets $\{ v _h = t \}$ are Lebesgue negligible. Then we can select $\e _h >0$ such that $| \{ v _h > \e _h \} | _\p  = \frac{|\Om|_\p}{2}$, so that 
 the sequence $u _h:= v _h - \e _h $  satisfies the condition $\int_\Om  (\sign  u_h )  \p = 0$.  Moreover, 
$u _h$ 
 converges in $L ^ 1_\p (\Om)$ to 
$v:=g _ \d - \e _0$, being $\e _0:= \lim_h \e _h$.  We notice that $\e _0 = 0$. Otherwise  the limit function $v$, which is strictly positive on $E$,  would be strictly  negative on the set $F \cup H$, which satisfies $|F \cup H | _\p > |E |_\p$. We conclude that 
$u _h \to g _\d$ in $L ^ 1 _\p (\Om)$, and 
$\int _\Om |\nabla u _h | \p = \int _\Om |\nabla v _h | \p  \to |D_\p g _\d | (\Om)$ as required. 

In the case when $|E|_\p = \frac{|\Om | _\p }{2}$, we consider the function $g = \chi _E - \chi _{\Om \setminus E}$, which satisfies 
\begin{equation}\label{f:gollodue} 
 \frac{  |D _\p \chi _E| (\Om)}{|E|_\p} =  \frac{  |D _\p g| (\Om)}{\int _\Om {|g | _\p} }  \,.\end{equation} 
Then, similarly as in the previous case, we take a sequence  of functions $\{ v _h \} \subset W ^ {1, 1} _\p (\Om)$, with Lebesgue-negligible level sets, such that
$v _h \to g $ in $L ^ 1 _\p (\Om)$, and $ \int _\Om |\nabla v _h | \p  \to |D_\p g | (\Om)$. 
We choose $\e _h >0$ such that $| \{ v _h > \e _h \} | _\p  = \frac{|\Om|_\p}{2}$, so that 
 the sequence $u _h:= v _h - \e _h $  satisfies  $\int_\Om  (\sign  u_h )  \p = 0$.  
 Moreover it converges in $L ^ 1_\p (\Om)$  to the function 
$v:=g  - \e _0$, with    $\e _0:= \lim_h \e _h$. Since the equality \eqref{f:gollodue} continues to hold with $g$ replaced by $v$, we conclude that
the sequence $\{u _h \}$ satisfies \eqref{f:gollo} as required. 
  \qed

 \bigskip

\underbar{Proof of Proposition \ref{p:ad}}. Throughout the proof, we write for brevity $\s$, $\overline \s$,  $\widetilde \s$,   and $\s ^*$ in place of 
$\s  (\Om,\p)$, $\overline \s   (\Om,\p)$,  $\widetilde \s (\Om, \p)$,    and $\s ^*   (\Om, \p)$.  

Let us show first that $  \widetilde  \s    =\s  ^*$.    
The inequality $\widetilde \s \leq \s ^*$ is immediate by choosing competitors for 
 $\widetilde \s$ 
 in which $B= \Om \setminus A$.   To show the converse, let $(A, B)$, $(a, b)$ be optimal for $\s ^*$; we can 
 assume without loss of generality that $a = |B|_\p$ and $b = |A|_\p$.  Then from the inequality $\widetilde \s \leq \s ^*$, taking respectively $A$ or  $B$ as competitors for $\s ^*$, we get 
 $$ 
\frac {|B| _\p  \Per _\p (A, \Om) + |A|_\p  \Per_\p (B, \Om)} {2|A| _\p |B| _\p} \leq  \min \Big \{ 
\frac{1}{2} \frac{|\Om|_\p \Per _\p (A, \Om) }{|A| _\p (|\Om | _\p- |A|_\p )}\, , \frac{1}{2} \frac{|\Om|_\p \Per_\p (B, \Om) }{|B|_\p  (|\Om |_\p - |B| _\p)} \Big \} \,.
$$ 
By writing the above inequality for the first of the two terms at the right hand side,  
we infer that
\begin{equation}\label{f:tre}
(|\Om |_\p - |A| _\p) \Per_\p (B, \Om) \leq |B|_\p \Per_\p (A, \Om) \,,
\end{equation} 
which since $|\Om |_\p - |A| _\p\geq |B|_\p$ implies $ \Per _\p (B, \Om) \leq  \Per_\p (A, \Om)$;
 similarly, by the inequality for the second term, we infer that $ \Per _\p(A, \Om) \leq  \Per _\p(B, \Om)$. 
But then  by \eqref{f:tre}  we have $|\Om |_\p = |A|_\p + |B|_\p$, namely $B = \Om \setminus A$,  
and hence $\widetilde \s = \s ^*$.

We now prove that $\overline\s \geq \s$. Let $u \in BV_\p (\Om; \R^+)$ with $\int _\Om u  \p = 0$.  
By Proposition \ref{p:BV} (i) there exists a sequence of functions $\{ u _ h \} \subseteq \mathcal D ( \R ^N) $ such that $u _ h \to u $ in $L ^ 1 _\p  (\Om)$  and $\int _\Om |\nabla u _h |  \p \to  |D _\p u| ( \Om ) $.    
Then $\int _\Om u _h \p  \to 0$, so that the sequence $v _h:= u _h -  \frac{1}{|\Om | _\p } \int _\Om u _h \p $ 
still converges to $u$ in $L ^ 1  _\p (\Om)$ and satisfies 
 $\int _\Om |\nabla v _h |_\p  \to |D_\p u|  (\Om )$. Since  $\int _\Om v _h \p = 0$, we have 
 $$\frac {|D_\p u| (\Om) } {\int_ \Om |u| \p }  = \lim _h  \frac {\int _\Om |\nabla v_h| \p } {\int_ \Om |v_h|\p } \geq \s\,,$$
which  implies  $\overline \s \geq \s$.
Viceversa,  let $\{ u _h \} \subset W ^ { 1, 1} _\p (\Om)$ be a sequence with  $\int_\Om  u _h \p= 0$ and  
$\frac {\int _\Om |\nabla u_h| \p } {\int_ \Om |u_h| \p } \to \s$. 
By Proposition \ref{p:BV} (vi), up to subsequences we have  $u _h \to u$  in $L ^ 1 _\p (\Om)$.
   Then we have $\int _\Om u \p = 0$,  and by Proposition \ref{p:BV} (ii) 
 $$ \s = \lim _h   \frac {\int _\Om |\nabla u_h| \p } {\int_ \Om |u_h| \p } \geq  \frac {|D_\p u| (\Om) } {\int_ \Om |u|\p  }  \geq \overline \s \,,$$
 
 which also shows that $\overline \s$ is attained.   

It remains to compare $\overline \s$ with $\s ^*$. The inequality  $\overline \s \leq \s ^*$ is immediate, since   the  minimization problem which defines $\s ^*$ is obtained by  restricting the class of admissible functions to the family of functions $u \in BV _\p  (\Omega)$ 
of the type $u =  \chi _E - \frac{|E| _\p }{|\Om \setminus E| _\p}  \chi _{\Om \setminus E}$. 

To show the converse inequality, we closely follow the argument in \cite[Lemma 1]{cianchi89}.   Let $u \in BV_\p  (\Om; \R ^+)$.  By the coarea formula
stated in Proposition \ref{p:BV} (v), we have
$$|D_\p u | (\Om) =  \int _ 0 ^ {+ \infty} \Per _\p  ( \{ u > t \}, \Om ) \, dt \,.$$
On the other hand, if we set for brevity $\overline u := \frac{1}{|\Om| _\p } \int _\Om  u \p$ and $\eta (t):= | \{ u > t \}|_\p$, we have 
  $$\overline u := \frac{1}{|\Om | _\p } \int_0 ^ {+ \infty} \eta ( t) \, dt \, , \qquad u ( x) = \int _0 ^ {+ \infty} \chi _{\{u > t\}} \, dt \,,$$
and hence 
$$\begin{aligned} 
\int _\Om |u - \overline u |  \phi  & = \int _{\Om} \Big | \int _0 ^ { + \infty} \chi _{\{u > t\}} - \frac{\eta ( t)} {|\Om|  _\p}  \, dt\Big |  \phi  \\ 
& \leq   \int _0 ^ { + \infty} \Big [ \int _{\Om}   \big |   \chi _{\{u > t\}} \p  - \frac{\eta ( t)} {|\Om| _\p  }   \p \big |  \Big ]  \, dt  \\
& =   \int _0 ^ { + \infty} \frac{1} {|\Om| _\p }  \Big [   \int _{\Om}  \big |   |\Om|  _\p \p  \chi _{\{u > t\}} - \eta ( t)  \p \big | \Big ]  \, dt   \\
& =  \int _0 ^ { + \infty} \frac{1} {|\Om|  _\p}  \Big [   \int _{ \{ u > t \}}  \big (  |\Om| _\p - \eta ( t)  \big ) \p   +  
\int _{ \{ u  \leq t \} }   \eta(t)   \p      \Big ] 
& =  \int _0 ^ { + \infty} \frac{1} {|\Om|  _\p }  2 \eta ( t) \big ( |\Om| _\p  - \eta ( t) \big ) \, dt   \\ 
& =  \int _0 ^ { + \infty}  \Per_\p ( \{ u > t\}, \Om)  \frac{ 2 |\{ u > t \} | _\p |\Om   \setminus  \{ u > t \}| _\p  }{ |\Om| _\p \Per_\p ( \{ u > t\}, \Om) }\, dt  \,.
 \end{aligned}
$$ 
Therefore, 
$$\begin{aligned} \frac {|D_\p u| (\Om) } {\int_ \Om |u - \overline u |  \p } = \frac { \int _0 ^ { + \infty}  \Per_\p ( \{ u > t \}, \Om ) } { \int _0 ^ { + \infty}  \Per
_\p ( \{ u > t\}, \Om)  \frac{ 2 |\{ u > t \} | _\p  |\Om   \setminus  \{ u > t \}| _\p  }{ |\Om| _\p  \Per _\p ( \{ u > t\}, \Om) }\, dt   }   
& \geq \frac{1}{ \sup _ {E\subset \Om}  
\frac {2 |E| _\p|\Om \setminus E|_\p}  { |\Om |_\p \Per_\p (E, \Om) } }\\ \noalign{\medskip}  &=   \inf _ {{E\subset \Om}}  
 \frac { |\Om |_\p \Per_\p (E, \Om) }  {2 |E|_\p |\Om \setminus E|_\p}  = \s ^* 
 \,.  \end{aligned}$$ 
The same argument applies for $u \in BV_\p  (\Om; \R ^-)$. In view of the equality stated in Proposition \ref{p:BV}    (iv), we infer that, for any 
$u \in BV_\p  (\Om)$, it holds
$$\begin{aligned} \frac {|D_\p u| (\Om) } {\int_ \Om |u - \overline u |  \p }   & = \frac {|D_\p u^+| (\Om)  + |D_\p u^-| (\Om)} {\int_ \Om |u ^+ +  u ^ -  - ( \overline{ u ^+}  +  
 \overline{ u ^-} )  |  \p }  \geq \frac {|D_\p u^+| (\Om)  + |D_\p u^-| (\Om)} {\int_ \Om |u ^+ -   \overline{ u ^+}  | \p  + \int _\Om   
  |u ^- -   \overline{ u ^-}  |  \p } \\ \noalign{\medskip} & \geq
   \min \Big \{ \frac {|D_\p u^+| (\Om)}{\int_ \Om |u ^+ -  \overline{ u ^+}  | \p }, \frac{|D_\p u^-| (\Om)}{\int _\Om   
  |u ^- -   \overline{ u ^-}  | \p }  \Big \}  \geq \s ^* \,. \end{aligned} 
  $$ 
We conclude, by the arbitrariness of $u \in BV_\p  (\Om)$, that 
$\overline \s \geq \s ^*$.  In addition, if $u$ is optimal for $\overline \s$, we see that  also $u ^ \pm$ are optimal for $\s ^*$ and 
almost all its level sets
 $\{ u > t\} $ (namely, their characteristic functions) are optimal for $\s ^*$. 
  \qed

\bigskip

\textbf{Data Availability Statement.}
Data sharing not applicable to this article as no datasets were generated or analysed during the current study.

\bigskip
\textbf{Conflict of Interest Statement.} 
The Authors declare that this article does not involve any  conflict of interest.  

\bigskip 

{\bf Acknowledgments.}  DB was supported by the  ANR STOIQUES research programme.

\bigskip \bigskip 
\bibliographystyle{mybst}


\end{document}